\documentclass[11pt, lmargin=1.25in, rmargin=1.25in]{article}

\RequirePackage{amsthm,amsmath,amsfonts,amssymb}
\RequirePackage[authoryear]{natbib} 
\RequirePackage[colorlinks,citecolor=blue,urlcolor=blue,breaklinks]{hyperref}

\usepackage{bm}
\usepackage{mathrsfs}  
\usepackage{appendix}
\usepackage{commath}
\usepackage{dsfont}
\usepackage{subfig}
\usepackage{xcolor}
\usepackage{mathtools}
\usepackage{enumitem}
\usepackage{fancyhdr}
\usepackage{mathrsfs} 
\usepackage{geometry}
\usepackage{caption}
\usepackage{pdflscape}

\usepackage{xcolor}
\usepackage{todonotes}
\usepackage{mdframed} 

\usepackage{xcolor}
\usepackage{todonotes}

\newcommand{\cop}[2][]{%
\textbf{\textcolor{red!50}{#2}}%
\if\relax\detokenize{#1}\relax\else%
\todo[color=red!50]{\footnotesize #1}%
\fi
}

\usepackage{setspace}
\newcommand{\N}{\mathbb{N}}
\newcommand{\R}{\mathbb{R}}
\newcommand{\eps}{\varepsilon}
\renewcommand{\P}{\mathbb{P}}
\newcommand\E{\mathbb{E}}
\newcommand{\mc}{\mathcal}

\theoremstyle{plain}
\newtheorem{theorem}{Theorem}[section]
\newtheorem{lemma}[theorem]{Lemma}

\theoremstyle{remark}
\newtheorem{remark}{Remark}[section]

\newtheorem{assumption}{Assumption}[section]

\begin{document}

\title{High-dimensional Gaussian and bootstrap approximations for robust means}

\author{
\begin{tabular}{c}
Anders Bredahl Kock\footnote{Kock's research was supported by the European Research Council (ERC) grant number 101124535 -- HIDI (UKRI EP/Z002222/1).  He is also a member of, and grateful for support from, i) the Aarhus Center for Econometrics (ACE), funded by the Danish National Research Foundation grant number DNRF186,  and ii) the Center for Research in Energy: Economics and Markets (CoRe).} \\ 
\footnotesize	University of Oxford \\
\footnotesize Department of Economics\\
\footnotesize	10 Manor Rd, Oxford OX1 3UQ
\\
\footnotesize	{\footnotesize	\href{mailto:anders.kock@economics.ox.ac.uk}{anders.kock@economics.ox.ac.uk}} 
\end{tabular}
\begin{tabular}{c}
David Preinerstorfer \\ 
{\footnotesize WU Vienna University of Economics and Business} \\
{\footnotesize Institute for Statistics and Mathematics} \\
{\footnotesize Welthandelsplatz 1, 1020 Vienna} \\ 
{\footnotesize	 \href{mailto:david.preinerstorfer@wu.ac.at}{david.preinerstorfer@wu.ac.at}}
\end{tabular}
}

\date{First version: April 9, 2025 \\ This version: June, 2026}

\maketitle	

\begin{abstract}
Recent years have witnessed much progress on Gaussian and bootstrap approximations to the distribution of sums of independent random vectors with dimension~$d$ large relative to the sample size~$n$. However, for any number of moments~$m>2$ that the summands may possess, there exist distributions such that these approximations break down if~$d$ grows faster than the polynomial barrier~$n^{\frac{m}{2}-1}$. In this paper, we establish Gaussian and bootstrap approximations to the distributions of winsorized and trimmed means that allow~$d$ to grow at an exponential rate in~$n$ as long as~$m>2$ moments exist. The approximations remain valid under some amount of adversarial contamination. Our implementations of the winsorized and trimmed means do not require knowledge of~$m$. As a consequence, the approximation guarantees ``adapt'' to~$m$. 
\end{abstract}

\section{Introduction}\label{sec:Intro}
Let~$X_1,\hdots,X_n$ be a sample of i.i.d.~random vectors in~$\R^d$ with mean vector~$\mu$ and covariance matrix~$\Sigma$. Furthermore, let~$S_n=n^{-1/2}\sum_{i=1}^n \del[0]{X_i-\mu}$. Since the seminal paper of~\cite{chernozhukov2013gaussian} there has been much renewed interest in Gaussian approximations to the distribution of~$S_{n}$ when~$d=d(n)$ is large relative to $n$. Letting~$Z\sim\mathsf{N}_d(0,\Sigma)$ and $\mc{H}$ be the class of (generalized) hyperrectangles in~$\R^d$, that is the class of all sets of the form
\begin{equation*}
	H=\cbr[1]{x\in\R^d:a_j\leq x_j\leq b_j\text{ for all }j=1,\hdots, d},
\end{equation*}
where~$-\infty\leq a_j\leq b_j\leq \infty$ for all~$j = 1, \hdots, d$, increasingly refined upper bounds have been established on
\begin{equation}\label{eq:approxerror}
	\rho_n:=\sup_{H\in\mc{H}}\envert[2]{\P\del[1]{S_{n}\in H}-\P\del[1]{Z\in H}}
\end{equation} 	
and related quantities, cf., e.g., \cite{chernozhukov2017central, deng2020beyond, lopes2020bootstrapping, kuchibhotla2020high, das2021central, koike2021notes, kuchibhotla2021high, lopes2022central, chernozhuokov2022improved, fang2023high, chernozhukov2023nearly, koike2024high}. We refer to the review in~\cite{chernozhukov2023high} for further references. For example, when the entries of~$X_i=(X_{i,1},\hdots,X_{i,d})'$ are (uniformly) sub-exponential,~$\rho_n\to 0$ if~$\log(d)=o\del[1]{n^{1/5}}$, cf.~\cite{chernozhuokov2022improved}. Thus,~$d$ can grow exponentially fast with~$n$ and this rate can be further improved under additional assumptions on the distribution of the~$X_i$ such as, e.g.,~variance decay conditions on~$\Sigma$ as in~\cite{lopes2020bootstrapping} or eigenvalue conditions as in~\cite{fang2021high, kuchibhotla2020high, chernozhukov2023nearly}.

\begin{figure}[t]
	\begin{center}
		\hspace{0.5cm}\footnotesize $n=200$ and~$d=5{,}000$\hspace{3.5cm} 	$n=1{,}000$ and~$d=10{,}000$
	\end{center}
	\vspace{-1cm}
	\begin{center}
		\includegraphics[width=6.5cm]{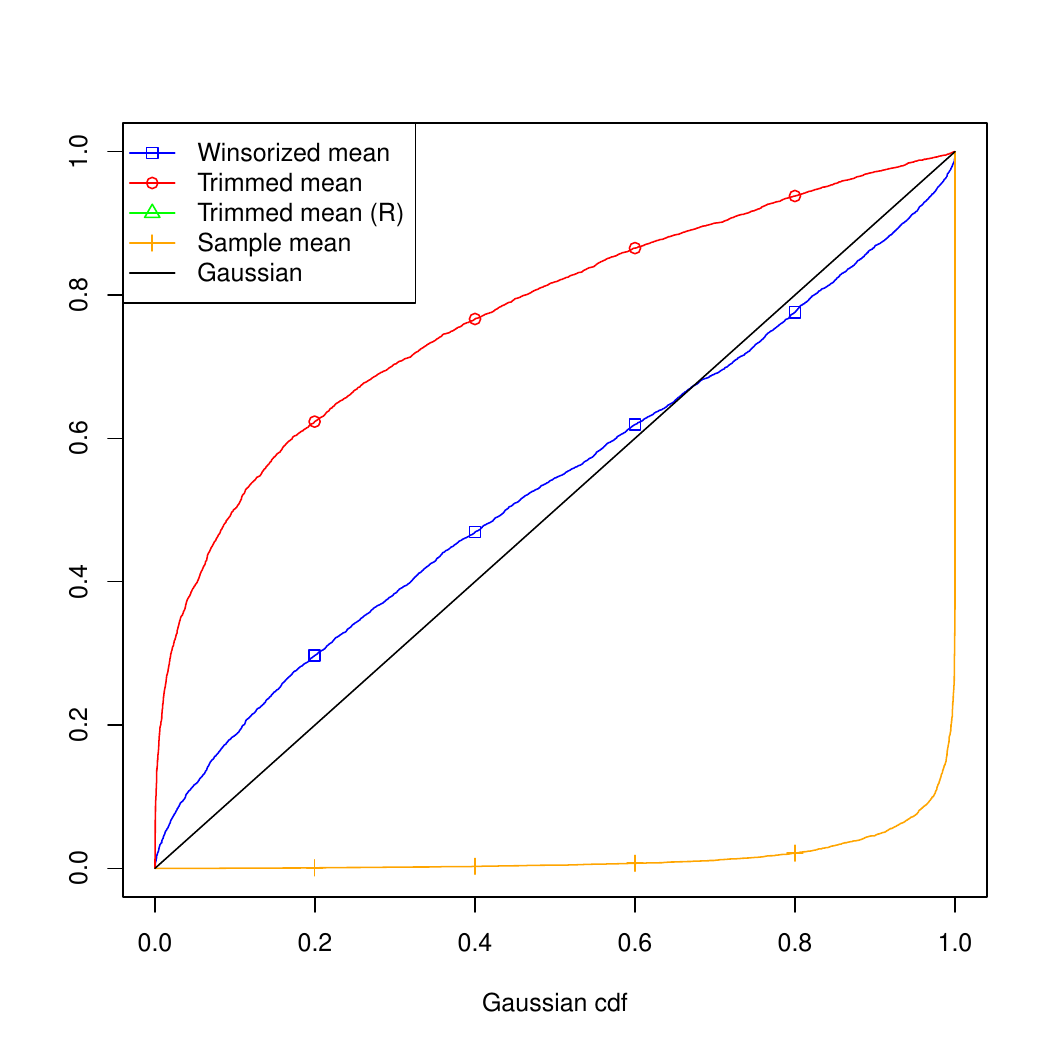}
		\hspace{-0.2cm}
		\includegraphics[width=6.5cm]{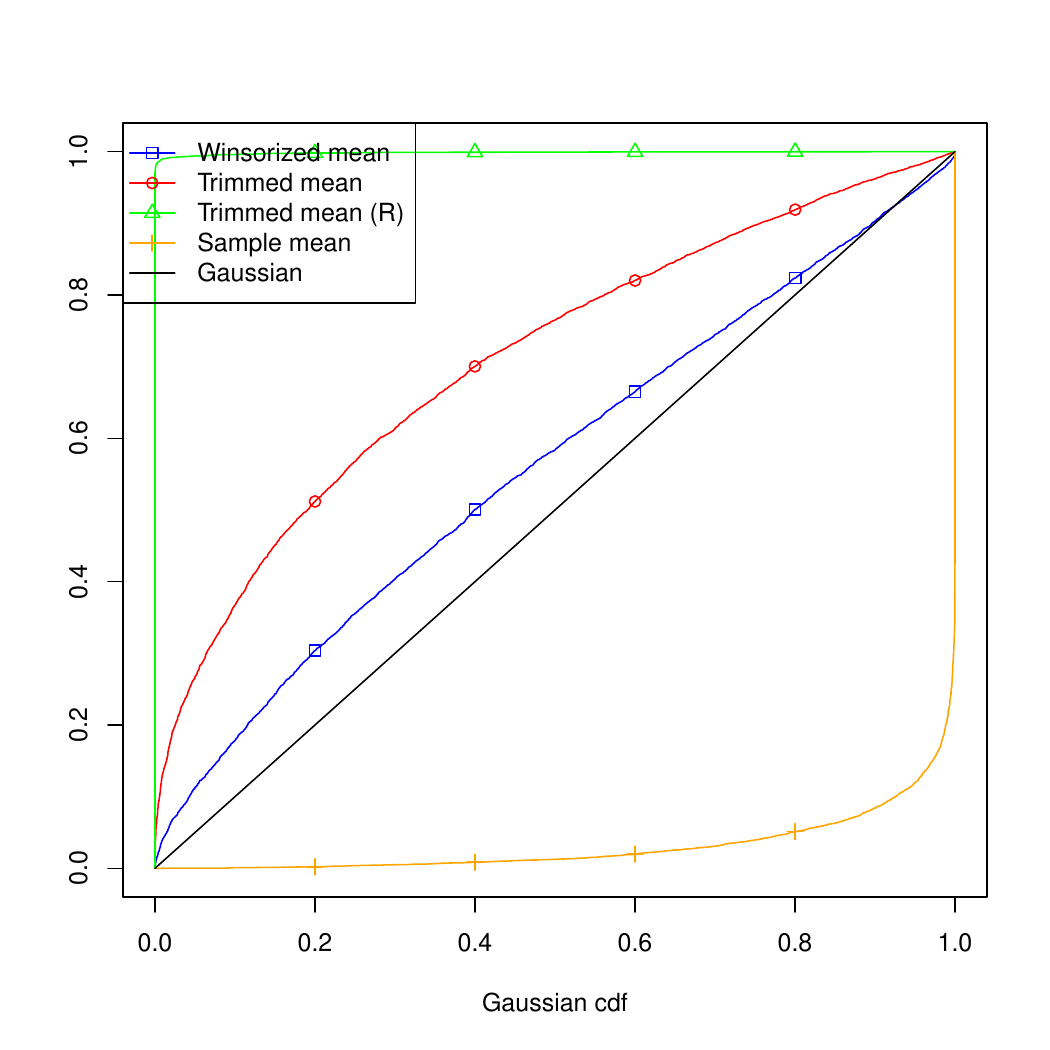}
	\end{center}
	\vspace{-0.5cm}
	\caption{\footnotesize P-P (probability-probability) plots: $X_{i,j}\sim$ i.i.d.~$t(3.01)$ implying~$\mu=0$. For each pair~$(n,d)$, the plots illustrate the Gaussian approximations to~$T_n$ for~$T_n\in\cbr[0]{S_{n,W},S_{n,T},S_{n,T,R},S_n}$; where~$S_{n,W}$ is our winsorized mean ({\color{blue} blue}) defined in~\eqref{eqn:winsmeandef}--\eqref{eq:epsfam},~$S_{n,T}$ is our trimmed mean ({\color{red} red}) defined in~\eqref{eq:trimean},~$S_{n,T,R}$ ({\color{green} green}) is the trimmed mean of~\cite{resende2024robust}, and~$S_n=n^{-1/2}\sum_{i=1}^nX_i$ ({\color{orange} orange}).
		The plots compare~$P(||T_{n}||_\infty\leq t)$,~$t\in\R$ to~$P(||Z||_\infty\leq t)$ (black) for~$Z\sim\mathsf{N}(0,\Sigma)$,~$\Sigma=\text{diag}(2.980,\hdots,2.980)$ the covariance matrix of~$X_1$ and recalling that for~$x\in\R^d$,~$||x||_\infty=\max_{1\leq j\leq d}|x_j|$. The absence of the trimmed mean of~\cite{resende2024robust} in the left panel indicates that it requires trimming more than~$n$ observations and is hence not possible to implement for the studied values of~$n$ and~$d$ (and the number of moments~$m<3.01$ that this procedure needs to specify in order to be implemented). Full implementation details are provided in Section~\ref{sec:sims}.}
	\label{fig:intro}
\end{figure}

Despite the progress on such high-dimensional Gaussian approximations for~$S_n$,~it follows from Remark~2 in~\cite{zhangwu2017} and Theorem 2.1 in~\cite{kock2024remark} that for every~$m\in(2,\infty)$ there exist~i.i.d.~random vectors~$X_1,\hdots,X_n$ with independent entries~$X_{i,j}\sim P_m$, and~$P_m$ depending neither on~$n$ nor~$d$, having mean zero, variance one, and finite~$m$th absolute moment, such that if for some~$\xi\in(0,\infty)$ it holds that~$\limsup_{n\to\infty}\frac{d}{n^{m/2-1+\xi}}>0$, then~$\limsup_{n\to\infty}\rho_n=1$. In particular, for any given~$m \in (2, \infty)$, the Gaussian approximation~$\rho_n \to 0$ does not hold uniformly over all distributions with bounded~$m$th moments when~$d$ grows exponentially in~$n$. Figure~\ref{fig:intro} illustrates how the Gaussian approximation to~$S_n$ can fail severely over~$\mc{H}$ when the~$X_{i,j}$ have relatively heavy tails. In particular, letting~$||x||_\infty=\max_{1\leq j\leq d}|x_j|$ for~$x\in\R^d$, it is seen that~$||S_n||_\infty$ has much heavier tails than~$||Z||_\infty$. 


Conversely, it is a simple consequence of Theorem~2 in~\cite{chernozhukov2023high}, cf.~Theorem 2.2 in~\cite{kock2024remark}, that~$\rho_n\to 0$ uniformly over a large class of distributions with bounded~$m$th moments if there exists a~$\xi\in(0,\infty)$ such that $\lim_{n\to\infty}\frac{d}{n^{m/2-1-\xi}}=0$. Hence, a critical phase transition occurs for the asymptotic behaviour of~$\rho_n$ at~$d=n^{m/2-1}$. As~$d$ passes this threshold from below, the limit of~$\rho_n$ jumps from zero to one. 

Motivated by this phase transition,~\cite{resende2024robust} recently  studied the case where~$S_n$ is replaced by a suitably \emph{trimmed} mean and~$\mc{H}$ is replaced by the subfamily~$\mc{R}$ consisting of all sets of the form~$R=\cbr[0]{x\in\R^d:x_j\leq t \text{ for all }j=1,\hdots,d}$, where~$t\in\R$. He obtained Gaussian approximations that are informative even when the~$X_i$ only possess~$m>2$ moments and~$d$ grows exponentially fast in~$n$. This is of fundamental importance, as it shows that one can break through the barrier~$d = n^{m/2 - 1}$ faced by~$\rho_n$, which is based on~$S_n$. The exact permitted growth rate of~$d$ depends on~$m$. As~$m\to \infty$, his result allows~$d$ to grow almost as fast as~$\exp(n^{1/6})$ for the Gaussian approximation as well as an empirical bootstrap and as fast as~$\exp(n^{1/8})$ for a multiplier bootstrap. Furthermore, the trimming ensures that these approximations remain valid even when some of the~$X_i$ have been adversarially contaminated prior to being given to the statistician. This is in stark contrast to statistics based on the sample mean~$S_n$, which have a breakdown point of~$1/n$ ($S_n$ can be changed to any value by manipulating only one of the vectors~$X_i$).

We also mention the work of~\cite{liu2024robust} who, motivated by the poor performance of~$S_n$ in the presence of heavy tails, even established a dimension-independent bootstrap approximation over~$\mc{R}$ for certain \emph{robust max-statistics} related to the winsorized means we study (but working with a ``variance-based'' instead of a ``quantile-based'' winsorization) under the conditions of~$L^4$-$L^2$ moment equivalence, a variance decay condition on~$\Sigma$, and restrictions on the Frobenius norm of certain submatrices of the correlation matrix of the~$X_i$. Robustness to outliers or other sources of (adversarial) data contamination were not investigated. 

In the present paper, we obtain Gaussian approximations to the distributions of winsorized and trimmed means over~$\mc{H}$ without imposing any structural assumptions on~$\Sigma$ (apart from positive variances). Our Gaussian approximations are valid over all of~$\mc{H}$ rather than the smaller family~$\mc{R}$ studied in~\cite{resende2024robust}, which is crucial for the study of power properties of max-type tests based on our winsorized and trimmed means; cf. the results in the companion paper~\cite{kp26rff}. In contrast to the present paper, which focuses exclusively on the canonical problem of Gaussian approximations in~$\R^d$,~\cite{resende2024robust} also considers Gaussian approximations over VC-subgraph classes of functions and applies his results to vector mean estimation under general norms.

One limitation of the trimmed mean studied in Theorem 2 in~\cite{resende2024robust} is that the amount of trimming needed depends on~$m$, which is typically unknown in practice. If one implements the trimmed mean studied in~\cite{resende2024robust} with an~$m$ higher than the actual number of moments that the~$X_i$ possess, one does not have any approximation guarantees, whereas the guarantees one obtains may be suboptimal if the~$X_i$ possess more moments than the~$m$ used in the implementation of the trimmed mean. What is more, \emph{even when~$m$ is known}, it is often the case that the required number of trimmed observations for the estimator studied in \cite{resende2024robust}
\begin{itemize}
	\item[i)] either exceeds the sample size~$n$, making that trimmed mean unimplementable (cf.~the left panel of Figure~\ref{fig:intro} where the  trimmed mean studied in \cite{resende2024robust} is absent);
	\item[ii)] or is smaller than~$n$ but so large that the Gaussian approximation to the distribution of that trimmed mean is not practically useful (cf.~the right panel of Figure~\ref{fig:intro} showing that the tails of the trimmed mean are much lighter than those of the approximating Gaussian due to excessive trimming).
\end{itemize}

In contrast, to implement the winsorized and trimmed means we investigate, one does not need to know the number of moments~$m$ that the~$X_i$ possess. The winsorization/trimming used in the present article draws from ideas in~\cite{LM21} and~\cite{Wins1}, is quantile-based, and ``adapts'' to~$m$. Figure~\ref{fig:intro} illustrates, in particular, that the winsorized mean we study obeys a tight Gaussian approximation even for rather large~$d$.

Furthermore, our results only require~$\log(d)=o\del[1]{n^{\frac{m-2}{5m-2}}}$, which is exponential for all~$m>2$, albeit with a small exponent for~$m$ close to two. As~$m\to\infty$,~$d$ is allowed to grow almost as fast as~$\exp(n^{1/5})$, improving on the rate in~\cite{resende2024robust}, and ``recovering'' the best known rate  for Gaussian approximations based on the sample mean of~$X_i$ with \emph{sub-exponential} entries (cf.~Remark 2 in~\cite{chernozhukov2023high}). This rate remains valid for the bootstrap procedures that we consider and all results are robust to some adversarial contamination. 

Our bootstrap approximations are based on a novel ``winsorized'' covariance matrix estimator, for which we establish performance guarantees in the presence of adversarial contamination and only requiring~$m > 2$ moments. Because the winsorization mechanism used in the construction of that covariance estimator is quantile-based, the estimator does not require knowledge of any unknown population quantities. In particular, its performance guarantees adapt to the unknown~$m$, which may be of independent interest.

Gaussian approximations for winsorized means are obtained in Section~\ref{sec:gaussapproxwins}, whereas bootstrap approximations are studied in Section~\ref{sec:Bootstrap}. Some results on trimmed means are established in Section~\ref{sec:trimmedmean}. Numerical results in addition to the ones already discussed in Figure~\ref{fig:intro} above are presented in Section~\ref{sec:sims}.

\section{Gaussian approximations for winsorized means}\label{sec:gaussapproxwins}

We first present our approximations to the distributions of high-dimensional winsorized means. Section~\ref{sec:trimmedmean} outlines the corresponding results for the version of the trimmed means we study.

Recall that~$X_1,\hdots,X_n$ is a sample of i.i.d.~random vectors in~$\R^d$, where~$X_i=(X_{i,1},\hdots,X_{i,d})'$ for~$i=1,\hdots,n$. Let~$\mu=(\mu_{1},\hdots,\mu_d)'=\E X_1$, $\Sigma$ be the covariance matrix of~$X_1$, and for~$m\in[2,\infty)$ let~$\sigma_{m,j}^m:=\E|X_{1,j}-\mu_j|^m$, all of which are well-defined under Assumption~\ref{ass:setting} below. We suppress the dependence of~$d=d(n)$ on~$n$ in our notation. 

In this section, our main focus is to establish Gaussian approximations for winsorized means that are valid for~$d$ growing exponentially in~$n$ imposing only that the~$X_{1,j}$ possess~$m>2$ moments,~$j=1,\hdots,d$. An added benefit of the winsorization is that the Gaussian approximations are robust to some amount of \emph{adversarial contamination}. Under such contamination an adversary inspects the sample and returns a corrupted sample~$\tilde{X}_1,\hdots,\tilde{X}_n$ to the statistician satisfying that
\begin{equation}\label{eq:contamfrac}
	\envert[1]{\cbr[1]{i\in\cbr[0]{1,\hdots,n}:\tilde{X}_i\neq X_i}}
	\leq
	\overline{\eta}_n n,
\end{equation}
where~$\overline{\eta}_n\in[0,1/2)$ is a non-random and known upper bound on the fraction of contaminated observations. Which of~$\tilde{X}_i$ differ from~$X_i$ as well as their values can depend on the uncontaminated sample~$X_1,\hdots, X_n$. Adversarial contamination has become a popular criterion to evaluate robustness of a statistic against, as it allows for many forms of data manipulation, cf.~\cite{lai2016agnostic}, \cite{cheng2019high}, \cite{diakonikolas2019robust}, \cite{hopkins2020robust}, \cite{LM21}, \cite{minsker2021robust}, \cite{bhatt2022minimax}, \cite{depersin2022robust}, \cite{dalalyan2022all}, \cite{minasyan2023statistically}, \cite{minsker2023efficient}, \cite{oliveira2025finite}. The recent book by~\cite{diakonikolas2023algorithmic} provides further references and discussion of various contamination settings. Since the sample mean has a breakdown point of~$1/n$, Gaussian approximations based on~$S_n$ are not robust to adversarial contamination (or large outliers). 

In all asymptotic statements~$n\to\infty$. Throughout, we impose the following assumption (for various values of~$m$).
\begin{assumption}\label{ass:setting}
	The~$X_1,\hdots,X_n$ are i.i.d.~random vectors in~$\R^d$ with~$\E |X_{1,j}|^m<\infty$ for some~$m\in(2,\infty)$ and all~$j=1,\hdots,d$. There exist strictly positive constants~$b_1$ and~$b_2$ such that~$\min_{j=1,\hdots,d}\sigma_{2,j}\geq b_1$ and $\sigma_m:=\max_{j=1,\hdots,d}\sigma_{m,j}\leq b_2$. The actually observed adversarially contaminated random vectors (in~$\R^d$) are denoted~$\tilde{X}_1,\hdots,\tilde{X}_{n}$ and satisfy~\eqref{eq:contamfrac}.	
\end{assumption}

Imposing lower and upper bounds~$b_1$ and~$b_2$ on moments of the~$X_{1,j}$ is commonplace when establishing upper bounds on~$\rho_n$ in~\eqref{eq:approxerror}, cf., e.g., the results in the overview~\cite{chernozhukov2023high}. Note that~$b_1$ and~$b_2$ are describing a statistical model and are \emph{not} intended to be optimized by the user (our procedures do not require knowledge of~$b_1$ or~$b_2$). Let us finally emphasize that all of our results are valid (in particular) \emph{absent} adversarial contamination, i.e., for~$\overline{\eta}_n=0$, which is the case studied in the literature on upper bounds on~$\rho_n$ summarized in the introduction.

\subsection{The winsorized means}
For real numbers~$x_1,\hdots,x_n$, denote by~$x_1^*\leq \hdots\leq x_n^*$ their non-decreasing rearrangement. Let~$-\infty<\alpha\leq\beta<\infty$ and 
\begin{align*}
	\phi_{\alpha,\beta}(x)
	=
	\begin{cases}
		\alpha\qquad \text{if }x<\alpha\\
		x\qquad \text{if }x\in[\alpha,\beta]\\
		\beta\qquad \text{if }x>\beta.
	\end{cases}
\end{align*} 
We establish Gaussian and bootstrap approximations to the distribution of centered winsorized means~$S_{n,W}\in\R^d$ where
\begin{equation}\label{eqn:winsmeandef}
	S_{n,W,j}
	=
	n^{-1/2}\sum_{i=1}^n\del[1]{\phi_{\hat\alpha_j,\hat\beta_j}(\tilde{X}_{i,j})-\mu_j},\qquad j=1,\hdots,d,
\end{equation}
with~$\hat{\alpha}_j=\tilde X_{\lceil \eps_n n \rceil,j}^*$  and $\hat{\beta}_j=\tilde X_{\lfloor(1-\eps_n )n\rfloor+1,j}^*$ for a carefully chosen~$\eps_n\in(0,1/2)$. Note that the choice of~$\eps_n$ determines the fraction of observations that are winsorized. For each coordinate~$j$, the winsorization points~$\hat{\alpha}_j$ and~$\hat{\beta}_j$ are the lower- and upper-$\eps_n$ order statistics of the contaminated data~$\tilde{X}_{1,j},\hdots,\tilde{X}_{n,j}$. In this sense, the winsorization mechanism employed is quantile-based.\footnote{For the purpose of constructing estimators of~$\mu\in\R$ with finite-sample sub-Gaussian concentration properties, related winsorized mean estimators were recently studied in~\cite{LM21} and~\cite{Wins1}. The winsorized estimators studied in the latter article remove some practical limitations of the ones studied in the former. We therefore mainly build on the results in~\cite{Wins1}.} Under adversarial contamination it is clear that even~$S_{n,W}$ can perform arbitrarily badly unless at least the smallest and largest~$\overline{\eta}_n n$ observations are winsorized. Thus, one must choose~$\eps_n\geq \overline{\eta}_n$. Our Gaussian approximations below apply for any $\eps_n\in(0,0.5)$ of the form
\begin{equation}\label{eq:epsfam}
	\eps_n=\lambda_1\cdot \overline{\eta}_n +\lambda_2\cdot \frac{\log(dn)}{n},\qquad \text{with }\lambda_1\in(1,\infty)\text{ and }\lambda_2\in (0,\infty).
\end{equation} 
We impose~$\lambda_1>1$, since otherwise even~$S_{n,W}$ can be arbitrarily manipulated by the adversary. To make efficient use of the data, we recommend choosing~$\lambda_1$ close to but greater than one, e.g.,~$\lambda_1 = 1.01$. In the leading special case of~$\overline{\eta}_n=0$, which is the setting in which high-dimensional Gaussian approximations for the sample mean based~$S_n$ have been studied (cf.~the literature summarized in the introduction), one has that~$\eps_n=\lambda_2\cdot \log(dn)/n$. The simulations in Section~\ref{sec:sims} suggest that~$\lambda_2\approx 0.1$ works well in practice.

\emph{Without further mentioning,~$\lambda_1\in(1,\infty)$ and~$\lambda_2\in(0,\infty)$ are taken as fixed throughout the remainder of the paper.}

\subsection{Gaussian approximation}\label{ss:gapp}

We now present a high-dimensional Gaussian approximation result for~$S_{n,W}$ over~$\mc{H}$ (defined prior to~\eqref{eq:approxerror}) in the form of an upper bound on
\begin{equation*}
	\rho_{n,W}:=\sup_{H\in\mc{H}}\envert[2]{\P\del[1]{S_{n,W}\in H}-\P\del[1]{Z\in H}}, \qquad\text{where }Z\sim \mathsf{N}_d(0,\Sigma).
\end{equation*}
We assume throughout that~$d \geq 2$ and~$n > 3$ (such that, e.g.,~$\sqrt{\log(d)}>0$). 
\begin{theorem}\label{thm:HDGauss}
	Let Assumption~\ref{ass:setting} be satisfied with~$m>2$. If~$\eps_n\in(0,1/2)$, with~$\eps_n$ as in~\eqref{eq:epsfam}, then, for a constant~$C$ depending only on~$b_1,b_2,\lambda_1,\lambda_2$, and~$m$, we have
	\begin{equation}\label{eq:Gaussapprox}
		\begin{aligned}
			\rho_{n,W}
			&\leq 	
			C\del[4]{\sbr[3]{\frac{\log^{5-\frac{2}{m}}(dn)}{n^{1-\frac{2}{m}}}}^{\frac{1}{4}}
				+
				\sbr[3]{\overline{\eta}_n^{1-\frac{1}{m}}+\sbr[2]{\frac{\log(dn)}{n}}^{1-\frac{1}{m}}}\sqrt{n\log(d)}}\\
			&\quad +
			C\del[4]{\log^2(d)\sbr[3]{\overline{\eta}_n^{1-\frac{2}{m}}+\sbr[2]{\frac{\log(dn)}{n}}^{1-\frac{2}{m}}}}^{1/2}.
		\end{aligned}
	\end{equation}
	In particular,~$\rho_{n,W}\to 0$ if~$\sqrt{n\log(d)}\overline{\eta}_n^{1-\frac{1}{m}}\to 0$ and~$\log(d)/n^{\frac{m-2}{5m-2}}\to 0$.
\end{theorem}
Consider the case of~$\overline{\eta}_n=0$. Theorem~\ref{thm:HDGauss} then shows that winsorized means can break through the \emph{polynomial} growth rate barrier~$d=n^{m/2-1}$ that~$d$ must obey for the Gaussian approximation error to the distribution of~$S_n$, i.e.,~$\rho_n$ in~\eqref{eq:approxerror}, to converge to zero. In particular,~$\rho_{n,W}\to 0$ if only~$\log(d)=o\del[1]{n^\frac{m-2}{5m-2}}$, which allows for~$d$ growing \emph{exponentially} in~$n$ for any~$m>2$, albeit with a small exponent for~$m$ close to two. Thus, the winsorized mean obeys a Gaussian approximation result over~$\mc{H}$ under heavy tails and adversarial contamination, which is in analogy to the Gaussian approximation result over~$\mc{R}\subsetneq\mc{H}$ in~\cite{resende2024robust} for the trimmed mean analyzed there. However, let us highlight the following key differences, which we already commented on in the introduction: The implementation of the trimmed mean estimator in~\cite{resende2024robust} requires knowledge of the (typically) unknown~$m$. Thus, if the~$X_i$ have fewer moments than the chosen~$m$, then there are no approximation guarantees. If, on the other hand, the~$X_i$ have more moments than the specified~$m$, then the guarantees may be suboptimal. In contrast, the implementation of our winsorized mean does not depend on~$m$ --- it ``adapts'' to it. Furthermore, we recover the rate~$d=o\del[1]{\exp(n^{1/5})}$ as~$m\to\infty$, which is currently the best available for Gaussian approximations based on~$S_n$ with sub-exponential~$X_i$, instead of~$d=o\del[1]{\exp(n^{1/6})}$ for the trimmed mean analyzed in~\cite{resende2024robust}. To make this more concrete, consider the case where the entries of the~$X_i$ are heavy-tailed in the sense of possessing exactly, e.g.,~$m=4$ moments. Then,~$\rho_{n,W}\to 0$ if~$\log(d)=o(n^{1/9})$, whereas the trimmed mean analyzed in Theorem 2 in \cite{resende2024robust} allows~$\log(d)=o(n^{3/35})$ if implemented with~$m=3$, $\log(d)=o(n^{4/35})$ if implemented with the true~$m=4$, and provides no guarantees if implemented with~$m>4$. Thus, unless one has reliable information on the number of moments the data possesses (in this case: knows~$m\approx 4$), our method adapting to~$m$ allows for larger~$d$. Finally, Section~\ref{sec:sims} (cf.~also Figure~\ref{fig:intro} in the Introduction) shows that the trimmed mean of~\cite{resende2024robust} i) often requires trimming more observations than the sample size~$n$, rendering it unimplementable, or ii) its distribution is very poorly approximated by that of~$Z$. In contrast, the distribution of~$S_{n,W}$ is closely approximated by that of~$Z$.

To prove Theorem~\ref{thm:HDGauss}, we first show that the order statistics~$\hat{\alpha}_j$ and~$\hat{\beta}_j$ can (essentially) be replaced by closely related population quantiles~$Q_{\eps_n,j}$ and~$Q_{1-\eps_n,j}$ of the~$X_{1,j}$, such that one can analyze the non-random winsorization functions~$\phi_{n,j}(x):=\phi_{Q_{\eps_n,j},Q_{1-\eps_n,j}}(x)$ instead of~$\phi_{\hat{\alpha}_j,\hat{\beta}_j}(x)$ for all~$j=1,\hdots,d$:\footnote{A similar control of the order statistics was also used in~\cite{LM21} and~\cite{Wins1}. \cite{resende2024robust} studied the trimmed mean by relating it to a winsorized mean.}
\begin{align*}
	S_{n,W,j}
	\approx 
	\underbrace{\frac{1}{n^{1/2}}\sum_{i=1}^n\sbr[1]{\phi_{n,j}(\tilde{X}_{i,j})-\mu_j}}_{T_{n,j}}
	&=
	\underbrace{\frac{1}{n^{1/2}}\sum_{i=1}^n\sbr[1]{\phi_{n,j}(\tilde{X}_{i,j})-\phi_{n,j}(X_{i,j})}}_{I_{n,j,1}}\\
	&+
	\underbrace{\frac{1}{n^{1/2}}\sum_{i=1}^n\sbr[1]{\phi_{n,j}(X_{i,j})-\E\phi_{n,j}(X_{i,j})}}_{I_{n,j,2}} \\
	&+\underbrace{\frac{1}{n^{1/2}}\sum_{i=1}^n\sbr[1]{\E\phi_{n,j}(X_{i,j})-\mu_j}}_{I_{n,j,3}}.
\end{align*}
The term~$\max_{j=1,\hdots,d}|I_{n,j,1}|$ isolates the effect of the adversarial contamination of the data and $\max_{j=1,\hdots,d}|I_{n,j,3}|$ quantifies how far the winsorized means~$\E\phi_{n,j}(X_{1,j})$ are from the population means~$\mu_j$ of interest. Next, note that~$I_{n,j,2}$ is a sample average of \emph{bounded} i.i.d.~random variables. Thus, letting~$Z_{n}\sim \mathsf{N}_d(0,\Sigma_{\phi_n})$ with~$\Sigma_{\phi_n}$ being the covariance matrix of~$I_{n,2}=(I_{n,1,2},\hdots,I_{n,d,2})'$, one can apply, e.g.,~the Gaussian approximation for sums of sub-exponential random vectors from~\cite{chernozhuokov2022improved} to this term to show that
\begin{equation*}
	\sup_{H\in\mc{H}}\envert[2]{\P\del[1]{I_{n,2}\in H}
		-
		\P\del[1]{Z_{n}\in H}}\qquad \text{is small}.
\end{equation*}
Furthermore, we show that~$\max_{1\leq j,k\leq d}\envert[0]{\Sigma_{\phi_n,j,k}-\Sigma_{j,k}}$ is sufficiently small for the Gaussian-to-Gaussian comparison inequality as stated in Proposition 2.1 in the previous reference (cf.~also Proposition 2 in~\cite{chernozhukov2023high}) to imply that
\begin{equation*}
	\P\del[1]{Z_{n}\in H}\qquad\text{in the previous display can be replaced by}\qquad 
	\P\del[1]{Z\in H}.	
\end{equation*}
Finally, we show that for~$l\in\cbr[0]{1,3}$ one has that~$\max_{j=1,\hdots,d}|I_{n,j,l}|$ are sufficiently small for a Gaussian anti-concentration inequality to imply that these  can be ``ignored". Therefore,~$\P\del[1]{I_{n,2}\in H}$ can be replaced by $\P\del[1]{T_n\in H}$ in the penultimate display, where~$T_n=(T_{n,1},\hdots,T_{n,d})'$.

\section{Bootstrap approximations for winsorized means}\label{sec:Bootstrap}

Similarly to Gaussian approximations for~$S_{n}$, cf.~the references in the introduction, the one for~$S_{n,W}$ in Theorem~\ref{thm:HDGauss} is not directly useful for statistical inference since the covariance matrix~$\Sigma$ of the approximating distribution~$\mathsf{N}_d(0,\Sigma)$ is typically unknown. Because the Gaussian-to-Gaussian comparison inequality as stated in Proposition 2.1 in~\cite{chernozhuokov2022improved} (cf.~also Proposition 2 in~\cite{chernozhukov2023high}) shows that for~$Z_1\sim\mathsf N_d(0,\Sigma^{(1)})$ and~$Z_2\sim\mathsf N_d(0,\Sigma^{(2)})$ with~$\min_{j=1,\hdots,d}\Sigma^{(2)}_{j,j}>b$ for some~$b>0$, it holds that
\begin{equation}\label{eq:Gaussiancomparison}
	\sup_{H\in\mc{H}}\envert[2]{\P\del[1]{Z_{1}\in H}-\P\del[1]{Z_{2}\in H}}
	\leq
	C\del[2]{\log^2(d)\max_{1\leq j,k\leq d}\envert[1]{\Sigma^{(1)}_{j,k}-\Sigma^{(2)}_{j,k}}}^{1/2},	
\end{equation}
for some~$C=C(b)$, one can approximate the unknown~$\P\del[0]{Z\in H}$ from Theorem~\ref{thm:HDGauss} if an estimator~$\hat{\Sigma}_n$ satisfying an upper bound on~$\log^2(d)\max_{1\leq j,k\leq d}\envert[1]{\hat{\Sigma}_{n,j,k}-\Sigma_{j,k}}$ can be exhibited. The sample covariance matrix can be used for~$d$ growing exponentially in~$n$ when i)~$\overline{\eta}_n=0$ and ii) the~$X_i$ have sub-exponential entries.\footnote{See Section 4.1 of~\cite{kuchibhotla2022moving} for properties of the sample covariance matrix when the~$X_i$ have sub-Weibull entries (generalizing sub-exponential distributions).} However, since we allow for~$\overline{\eta}_n>0$ and only impose the existence of~$m>2$ moments, the sample covariance matrix cannot be used in our context.

There has been a recent interest in constructing estimators of~$\Sigma$ which perform well under heavy tails (and frequently also adversarial contamination). For example, estimators with precision guarantees in the entrywise maximal distance~$\max_{1\leq j,k\leq d}\envert[1]{\hat{\Sigma}_{n,j,k}-\Sigma_{j,k}}$ needed in~\eqref{eq:Gaussiancomparison}, have been proposed in~\cite{ke2019user} in the setting of heavy-tailed~$X_i$. These estimators are based on, e.g., entrywise truncation or the median-of-means principle and the practical choice of the needed tuning parameters (which depend on unknown population quantities) is also discussed there. 

In the next section, we construct an estimator~$\tilde{\Sigma}_n$, say, which is based on suitably winsorized observations~$\tilde{X}_i$, and for which we establish performance guarantees even for~$\overline{\eta}_n > 0$ and when the $X_i$ possess only~$m > 2$ moments. Our estimator does not depend on unknown population quantities and its performance guarantees ``adapt'' to the unknown~$m$. We stress that Theorem~\ref{thm:HDBootstrap} below is modular in the sense that it remains valid for any other estimator~$\hat{\Sigma}_n$ satisfying a bound as in Theorem~\ref{thm:covestimGRAM}. 	

\subsection{Estimating~$\Sigma$}\label{sec:estsigma}
Imposing only~$m>2$ moments to exist we now construct an estimator of~$\Sigma$ with precision guarantees in the maximal entrywise norm. Set
\begin{equation}\label{eq:epsprime}
	\eps_n'= \lambda_1'\cdot\overline{\eta}_n+\lambda_2'\cdot\frac{\log(d^2n)}{n},\qquad \lambda_1'\in(1,\infty)\text{ and }\lambda_2'\in (0,\infty).
\end{equation}
\emph{Like for~$\lambda_1$ and~$\lambda_2$, we take~$\lambda_1'$ and~$\lambda_2'$ as fixed throughout}.

Writing~$\hat{a}_j=\tilde{X}_{\lceil \eps_n'n\rceil,j}^*$,~$\hat{b}_j=\tilde{X}_{\lfloor (1-\eps_n')n\rfloor+1,j}^*$, and~$\tilde{\mu}_{n,j}=n^{-1}\sum_{i=1}^n\phi_{\hat{a}_j,\hat{b}_j}(\tilde{X}_{i,j})$ for~$j=1,\hdots,d$, define~$\tilde{\Sigma}_n$ as the matrix with entries 
\begin{align}\label{eq:tildeSigma}
	\tilde{\Sigma}_{n,j,k}
	=
	\frac{1}{n}\sum_{i=1}^n\sbr[1]{\phi_{\hat{a}_j,\hat{b}_j}(\tilde{X}_{i,j})-\tilde{\mu}_{n,j}}\sbr[1]{\phi_{\hat{a}_k,\hat{b}_k}(\tilde{X}_{i,k})-\tilde{\mu}_{n,k}},\qquad 1\leq j,k\leq d.
\end{align}
The estimator~$\tilde{\Sigma}_n$, which is a winsorized analog of the sample covariance matrix, is positive semi-definite and symmetric by virtue of being a Gram matrix. Note that the choice of~$\eps_n'$ (which slightly differs from~$\eps_n$ as defined in Equation~\eqref{eq:epsfam}) determines the fraction of observations that are winsorized for the purpose of \emph{covariance} estimation. Solely for the purpose of stating the precision guarantee for~$\tilde{\Sigma}_n$ in Theorem~\ref{thm:covestimGRAM} below, we introduce the condition
\begin{align}\label{eq:epscond'}
	2\eps_n' +\frac{\log(d^2n)}{n}+\sqrt{\del[2]{\frac{\log(d^2n)}{n}}^2+4\frac{\log(d^2n)}{n}\eps_n'}<1.	
\end{align}
We stress that~\eqref{eq:epscond'} is \emph{not needed} for any of the subsequent Gaussian and bootstrap approximations of this paper to hold, as these will be argued to be trivially true when~\eqref{eq:epscond'} is violated, cf., e.g., the closing paragraph of the proof of Theorem~\ref{thm:HDBootstrap}. Note also that in the leading special case of~$\overline{\eta}_n=0$, which is always imposed in the literature when studying Gaussian and bootstrap approximations to the distribution of~$S_n$,~\eqref{eq:epscond'} reduces to
\begin{align*}
	\del[2]{2\lambda_2'+1+\sqrt{1+4\lambda_2'}}\frac{\log(d^2n)}{n} <1,		
\end{align*}
which is often satisfied for~$\lambda_2'$ small. 

Concerning the practical choice of~$\lambda_1'$ and~$\lambda_2'$, we suggest to choose~$\lambda_1'$ close to but greater than one (to make efficient use of the data), e.g.,~$\lambda_1' = 1.01$. The simulations in Section~\ref{sec:sims} suggest that~$\lambda_2'\approx 0.07$ works well in practice. Note that for this choice, the number of observations winsorized in implementing~$\tilde{\Sigma}_n$ is typically the same as in implementing~$S_{n,W}$.
\begin{remark}
	The condition in~\eqref{eq:epscond'} is actually just a conservative (simple) sufficient condition for the following milder condition: Writing~$A_+'=1-\lambda_1'^{-1}\mathds{1}\cbr[0]{\overline{\eta}_n>0}\in(0,1]$,~$A_-'=1+\lambda_1'^{-1}\mathds{1}\cbr[0]{\overline{\eta}_n>0}\in [1,\infty)$, and denoting by~$W_0$ and~$W_{-1}$ the principal and lower branch of Lambert's~$W$ function (cf., e.g., \cite{Corless1996}), respectively,~\eqref{eq:epscond'} can be replaced by
	\begin{align*}
		\eps_n'\del[3]{-A_+'W_0\del[2]{-e^{-(\frac{\log(d^2n)}{\eps_n' n}+A_+')/A_+'}}-A_-'W_{-1}\del[2]{-e^{-(\frac{\log(d^2n)}{\eps_n' n}+A_-')/A_-'}}}<1.
	\end{align*}
	By \eqref{eq:welldefined'} of Lemma~\ref{lem:cControl'} in the appendix, the left-hand side of the previous display is upper bounded by the left-hand side of~\eqref{eq:epscond'}, explaining the condition in~\eqref{eq:epscond'}. Note that~\eqref{eq:epscond'} also implies~$\frac{\log(d^2n)}{n} <1$, which we repeatedly use in the proofs.
\end{remark}
\begin{theorem}\label{thm:covestimGRAM}
	Let Assumption~\ref{ass:setting} be satisfied with~$m>2$. If~$\eps_n'$ is chosen as in~\eqref{eq:epsprime} and satisfies~\eqref{eq:epscond'}, then for a constant~$C$ depending only on~$b_2,\lambda_1',\lambda_2'$ and~$m$, we have
	\begin{align}\label{eq:covestimGram}
		\P\del[4]{\max_{1\leq j,k\leq d}\envert[1]{\tilde{\Sigma}_{n,j,k}-\Sigma_{j,k}}> C\sbr[3]{\overline{\eta}_n^{1-\frac{2}{m}}+\del[2]{\frac{\log(dn)}{n}}^{1-\frac{1}{(m/2)\wedge 2}}}}
		\leq \frac{24}{n}.
	\end{align}	
\end{theorem}
Theorem~\ref{thm:covestimGRAM} shows that for any~$m>2$ it is possible for~$\max_{1\leq j,k\leq d}\envert[0]{\tilde{\Sigma}_{n,j,k}-\Sigma_{j,k}}$ to converge to zero in probability even when~$d$ grows exponentially in~$n$. Finally, we reiterate that the approximation in Theorem~\ref{thm:HDBootstrap} below remains valid for \emph{any} covariance estimator satisfying~\eqref{eq:covestimGram}.

\subsection{Bootstrap consistency}

Equipped with the estimator~$\tilde{\Sigma}_n$, the following theorem justifies approximating~$\P\del[1]{S_{n,W}\in H}$ for~$H\in\mc{H}$ by sampling repeatedly from~$\mathsf{N}_d(0,\tilde{\Sigma}_n)$.
\begin{theorem}\label{thm:HDBootstrap}
	Let Assumption~\ref{ass:setting} be satisfied with~$m>2$. Let~$n > 24$.  If~$\eps_n,\eps_n'\in(0,1/2)$, with~$\eps_n$ as in~\eqref{eq:epsfam} and~$\eps_n'$ as in~\eqref{eq:epsprime}, and~$\tilde{Z}\sim \mathsf{N}_d(0,\tilde{\Sigma}_n)$ conditionally on~$\tilde{X}_1,\hdots,\tilde{X}_n$, it holds with probability at least~$1-\frac{24}{n}$ that
	\begin{align*}
		\tilde{\rho}_{n,W}:&=\sup_{H\in\mc{H}}\envert[2]{\P\del[1]{S_{n,W}\in H}-\P\del[1]{\tilde{Z}\in H\mid\tilde{X}_1,\hdots,\tilde{X}_n}}\\
		&\leq
		\mathfrak{A}_n
		+C\del[4]{\log^2(d)\sbr[3]{\overline{\eta}_n^{1-\frac{2}{m}}+\del[2]{\frac{\log(dn)}{n}}^{1-\frac{1}{(m/2)\wedge 2}}}}^{1/2},
	\end{align*}	
	where~$\mathfrak{A}_n$ is the upper bound on~$\rho_{n,W}$ in~\eqref{eq:Gaussapprox} of Theorem~\ref{thm:HDGauss} and~$C$ is a constant depending only on~$b_1,b_2,\lambda_1',\lambda_2'$ and~$m$. 
	
	In particular,~$\tilde{\rho}_{n,W}\to 0$ in probability if~$\sqrt{n\log(d)}\overline{\eta}_n^{1-\frac{1}{m}}\to 0$ and~$\log(d)/n^{\frac{m-2}{5m-2}}\to 0$.
\end{theorem}
Draws from~$\mathsf{N}_d(0,\tilde{\Sigma}_n)$ can be obtained efficiently by, e.g., the following multiplier bootstrap: Let~$\xi_1,\hdots,\xi_n$ be i.i.d.~$N_1(0,1)$, independent of~$\tilde{X}_1,\hdots,\tilde{X}_n$, and set
\begin{equation*}
	S_{n,\text{MB}}
	=
	\frac{1}{\sqrt{n}}\sum_{i=1}^n\xi_iV_i \quad\text{where}\quad V_i=\del[1]{\phi_{\hat{a}_1,\hat{b}_1}(\tilde{X}_{i,1})-\tilde{\mu}_{n,1},\hdots,\phi_{\hat{a}_d,\hat{b}_d}(\tilde{X}_{i,d})-\tilde{\mu}_{n,d}}'\in\R^d.
\end{equation*}
Observe that conditionally on~$\tilde{X}_1,\hdots,\tilde{X}_n$ the distribution of~$S_{n,\text{MB}}$ is~$\mathsf{N}_d(0,\tilde{\Sigma}_n)$. Generating a draw from~$\mathsf{N}_d(0,\tilde{\Sigma}_n)$ via this multiplier bootstrap may be numerically preferable to first drawing from~$\mathsf{N}_d(0,\mathsf{I}_d)$ and then premultiplying this by~$\tilde{\Sigma}_n^{1/2}$ since the calculation of the matrix square root may be costly for~$d$ large. 

\section{Normalized winsorized means}\label{sec:Normalized}
In practice one often normalizes the data by estimates of the~$\sigma_{2,j}$ to bring the variables on the same scale.  We now describe how the Gaussian and bootstrap approximations established so far remain valid upon normalization by the diagonal elements~$\tilde{\sigma}_{n,j}=\tilde{\Sigma}_{n,j,j}^{1/2}$ of the robust estimator~$\tilde{\Sigma}_n$ of~$\Sigma$, cf.~\eqref{eq:tildeSigma}. Writing~$D=\text{diag}(\sigma_{2,1},\hdots,\sigma_{2,d})$ and~$\Sigma_0=D^{-1}\Sigma D^{-1}$ for the correlation matrix of the (centered) uncontaminated~$X_1,\hdots,X_n$, one has the following Gaussian approximation for the vector~$S_{n,W,S}$ of normalized winsorized means with elements (we leave the quotients undefined if one of the variance estimators equals~$0$, cf.~also Remark~\ref{rem:rem0var} below)
\begin{equation}\label{eq:S_nWS}
	S_{n,W,S,j}=\frac{1}{\sqrt{n}\tilde{\sigma}_{n,j}}\sum_{i=1}^n\sbr[1]{\phi_{\hat\alpha_j,\hat\beta_j}(\tilde{X}_{i,j})-\mu_j},\qquad j=1,\hdots,d.
\end{equation}
\begin{theorem}\label{thm:HDGauss_studentized}
	Let Assumption~\ref{ass:setting} be satisfied with~$m>2$.  If~$\eps_n,\eps_n'\in(0,1/2)$, with~$\eps_n$ as in~\eqref{eq:epsfam} and~$\eps_n'$ as in~\eqref{eq:epsprime}, then for~$Z'\sim\mathsf{N}_d(0,\Sigma_0)$,
	\begin{align}\label{eq:HDGaussstudentized}
		\rho_{n,W,S}
		:&=
		\sup_{H\in\mc{H}}\envert[2]{\P\del[1]{S_{n,W,S}\in H}-\P\del[1]{Z'\in H}}\notag\\
		&\leq
		C\del[4]{\mathfrak{A}_n+\sqrt{\log(d)\log(dn)}\sbr[3]{\overline{\eta}_n^{1-\frac{2}{m}}+\del[2]{\frac{\log(dn)}{n}}^{1-\frac{1}{(m/2)\wedge 2}}}},
	\end{align}
	where~$\mathfrak{A}_n$ is the upper bound on~$\rho_{n,W}$ in~\eqref{eq:Gaussapprox} of Theorem~\ref{thm:HDGauss} and~$C$ is a constant depending only on~$b_1,b_2,\lambda_1,\lambda_2,\lambda_1',\lambda_2'$ and~$m$.
	
	In particular,~$\rho_{n,W,S}\to 0$ if~$\sqrt{n\log(d)}\overline{\eta}_n^{1-\frac{1}{m}}\to 0$ and~$\log(d)/n^{\frac{m-2}{5m-2}}\to 0$.
\end{theorem}
Theorem~\ref{thm:HDGauss_studentized} allows for the same (exponential) growth rate of~$d$ (for all~$m>2$) and contamination rate~$\overline{\eta}_n$ as Theorem~\ref{thm:HDGauss} for non-normalized data. In analogy to~$\Sigma$ in Theorem~\ref{thm:HDGauss},~$\Sigma_0$ is unknown in Theorem~\ref{thm:HDGauss_studentized}. Letting~$\tilde{D}_n=\text{diag}\del[0]{\tilde{\sigma}_{n,1},\hdots,\tilde{\sigma}_{n,d}}$ and~$\tilde{\Sigma}_{n,0}=\tilde{D}_n^{-1}\tilde{\Sigma}_n\tilde{D}_n^{-1}$ (we leave the expression undefined if one of the variance estimators equals~$0$, cf.~also Remark~\ref{rem:rem0var} below), we have the following analogue to the bootstrap approximation in Theorem~\ref{thm:HDBootstrap}. 
\begin{theorem}\label{thm:HDBootstrap_studentized}
	Let Assumption~\ref{ass:setting} be satisfied with~$m>2$. Let~$n > 24$. If~$\eps_n,\eps_n'\in(0,1/2)$, with~$\eps_n$ as in~\eqref{eq:epsfam} and~$\eps_n'$ as in~\eqref{eq:epsprime}, and~$\tilde{Z}'\sim \mathsf{N}_d(0,\tilde{\Sigma}_{n,0})$ conditionally on~$\tilde{X}_1,\hdots,\tilde{X}_n$, it holds with probability at least~$1-\frac{24}{n}$ that
	\begin{align}
		\tilde{\rho}_{n,W,S}
		:&=
		\sup_{H\in\mc{H}}\envert[3]{\P\del[1]{S_{n,W,S}\in H}-\P\del[1]{\tilde{Z}'\in H\mid\tilde{X}_1,\hdots,\tilde{X}_n}}\notag\\
		&\leq
		\mathfrak{B}_n+ C\del[4]{\log^2(d)\sbr[3]{\overline{\eta}_n^{1-\frac{2}{m}}+\del[2]{\frac{\log(dn)}{n}}^{1-\frac{1}{(m/2)\wedge 2}}}}^{1/2}	\label{eq:StudentizeBS},
	\end{align}
	where~$\mathfrak{B}_n$ is the upper bound on~$\rho_{n,W,S}$ in~\eqref{eq:HDGaussstudentized} of Theorem~\ref{thm:HDGauss_studentized} and~$C$ is a constant depending only on~$b_1,b_2,\lambda_1',\lambda_2'$ and~$m$. 
	
	In particular,~$\tilde{\rho}_{n,W,S}\to 0$ in probability if~$\sqrt{n\log(d)}\overline{\eta}_n^{1-\frac{1}{m}}\to 0$ and~$\log(d)/n^{\frac{m-2}{5m-2}}\to 0$.
\end{theorem}

\begin{remark}\label{rem:rem0var}
	We have left the quantities in~\eqref{eq:S_nWS} and~$\tilde{\Sigma}_{n,0}$ undefined on the event where one of the variance estimators~$\tilde{\sigma}_{n,j}$ ($j = 1, \hdots, d$) equals~$0$. The proofs of Theorems~\ref{thm:HDGauss_studentized} and~\ref{thm:HDBootstrap_studentized} exploit that the just-mentioned event is negligible in regimes where the upper bounds in these theorems
	are informative. We sketch the argument here: Theorem~\ref{thm:covestimGRAM} shows that with high probability
	\begin{equation}\label{discusseqnS}
		\max_{1\leq j,k\leq d}\envert[1]{\tilde{\Sigma}_{n,j,k}-\Sigma_{j,k}} \leq C\sbr[3]{\overline{\eta}_n^{1-\frac{2}{m}}+\del[2]{\frac{\log(dn)}{n}}^{1-\frac{1}{(m/2)\wedge 2}}}.
	\end{equation}
	In regimes where the upper bounds in Theorems~\ref{thm:HDGauss_studentized} and~\ref{thm:HDBootstrap_studentized}
	are small, the term in brackets in~\eqref{discusseqnS} is small, so that then~$\tilde{\sigma}_{n,j} \approx \sigma_{2,j} \geq b_1 > 0$, the inequalities being due to Assumption~\ref{ass:setting}. Thus, with high probability~$\tilde{\sigma}_{n,j} > 0$ for~$j = 1, \hdots, d$.\footnote{The asymptotic argument just given shows that in regimes where the bounds in Theorems~\ref{thm:HDGauss_studentized} and~\ref{thm:HDBootstrap_studentized} are informative, the event where a variance estimator vanishes is negligible. For practitioners it is relevant to know that typically more can be said in finite samples upon closer inspection: fix a~$j \in \{1, \hdots, d\}$ and deduce from~\eqref{eq:tildeSigma} that~$\tilde{\sigma}_{n,j} = 0$ if and only if $$\tilde{X}_{\lceil \eps_n'n\rceil,j}^* = \tilde{X}_{\lceil \eps_n'n\rceil + 1,j}^* = \hdots = \tilde{X}_{\lceil (1-\eps_n')n\rceil,j}^*,$$ so that $\lceil (1-\eps_n')n\rceil - \lceil \eps_n'n\rceil + 1 \geq (1-2\eps_n')n$ of the \emph{contaminated} observations coincide in their~$j$-th coordinate. Since at most~$\overline{\eta}_n n$ observations are contaminated, at least~$(1-2\eps_n' - \overline{\eta}_n) n$ of the \emph{uncontaminated} observations then coincide in their~$j$-th coordinate. The latter event has non-negligible probability only if~$2\eps_n' + \overline{\eta}_n$ is close to or larger than~$1$ or if the distribution of~$X_{i,j}$ has a dominant discrete component. Both scenarios do not appear very relevant for practical applications. For example, if~$n > 3$, $(2\eps_n' + \overline{\eta}_n) \leq 1-(2/n)$ and the uncontaminated observations have a distribution absolutely continuous~w.r.t.~Lebesgue measure on~$\R^d$, the event where one of the variance estimators~$\tilde{\sigma}_{n,j}$ ($j = 1, \hdots, d$) equals zero in fact has probability zero.}
\end{remark}

\section{Trimmed means}\label{sec:trimmedmean}
In this section we show that the Gaussian and bootstrap approximations established so far for \emph{winsorized} means carry over to (suitably) \emph{trimmed} means allowing for exactly the same growth rates of~$\overline{\eta}_n$ and~$d$. 
For~$\eps_n$ as in~\eqref{eq:epsfam}, let~$I_n:=\cbr[1]{\lceil\eps_nn\rceil,\hdots,\lfloor(1-\eps_n)n\rfloor+1}$ with cardinality~$$|I_n|=\lfloor(1-\eps_n)n\rfloor+1-\lceil\eps_nn\rceil+1
=
n-2\lceil\eps_nn\rceil+2,$$
and consider the vector of trimmed means~$S_{n,T}\in\R^d$ with entries ($j = 1, \hdots, d$)
\begin{align}\label{eq:trimean}
	S_{n,T,j}
	=
	\frac{\sqrt{n}}{|I_n|}\sum_{i\in I_n}\sbr[1]{\tilde{X}_{i,j}^*-\mu_j}
	=
	\frac{\sqrt{n}}{n-2\lceil\eps_nn\rceil+2}\sum_{i=\lceil \eps_nn\rceil}^{\lfloor(1-\eps_n)n\rfloor+1}\sbr[1]{\tilde{X}_{i,j}^*-\mu_j}.
\end{align}
Analogously to the winsorized means in~\eqref{eqn:winsmeandef}, the trimmed mean in~\eqref{eq:trimean} can be implemented in a fully data-driven way: Its amount of trimming is governed by~$\eps_n$ in~\eqref{eq:epsfam}, which does not depend on any unknown population quantities.\footnote{As is typical in the literature on inference under data contamination, an upper bound~$\overline{\eta}_n$ on the contamination rate must be supplied, however.} This sets this trimmed mean apart from the one studied in Theorem 2 in~\cite{resende2024robust}, where the amount of trimming depends on~$m$. Furthermore, the bound in Theorem~\ref{thm:HDGaussTrim} below is valid over the larger family of sets~$\mc{H}$ instead of~$\mc{R}$. In analogy to~$\rho_{n,W}$ define
\begin{align*}
	\rho_{n,T}
	=
	\sup_{H\in\mc{H}}\envert[2]{\P\del[1]{S_{n,T}\in H}-\P\del[1]{Z\in H}}.
\end{align*}
The following theorem is the trimmed mean counterpart to Theorem~\ref{thm:HDGauss}. 

\begin{theorem}\label{thm:HDGaussTrim}
	Let Assumption~\ref{ass:setting} be satisfied with~$m>2$. If~$\eps_n\in(0,1/2)$, with~$\eps_n$ as in~\eqref{eq:epsfam}, then 
	\begin{equation}
		\rho_{n,T}
		\leq 	
		\mathfrak{A}_n
		+
		C\sqrt{n\log(d)}\del[3]{\overline{\eta}_n^{1-\frac{1}{m}}+\sbr[2]{\frac{\log(dn)}{n}}^{1-\frac{1}{m}}},
		\label{eq:trimGauss}
	\end{equation}
	where~$\mathfrak{A}_n$ is the upper bound on~$\rho_{n,W}$ in~\eqref{eq:Gaussapprox} of Theorem~\ref{thm:HDGauss} and~$C$ is a constant depending only on~$b_1,b_2,\lambda_1,\lambda_2$ and~$m$. 
	
	In particular,~$\rho_{n,T}\to 0$ if~$\sqrt{n\log(d)}\overline{\eta}_n^{1-\frac{1}{m}}\to 0$ and~$\log(d)/n^{\frac{m-2}{5m-2}}\to 0$.
\end{theorem}
We establish Theorem~\ref{thm:HDGaussTrim} by verifying that~$S_{n,T}$ is sufficiently close to~$S_{n,W}$ in the supremum-norm, in order to use a Gaussian anti-concentration inequality to deduce~\eqref{eq:trimGauss} from~\eqref{eq:Gaussapprox} in Theorem~\ref{thm:HDGauss}. The latter theorem was established using the proof strategy outlined in Section~\ref{ss:gapp} (making use of a Gaussian approximation result applicable to sums of bounded~i.i.d.~random variables).  This explains the presence of the additional summand on the right-hand side in~\eqref{eq:trimGauss}. The same discussion as that following Theorem~\ref{thm:HDGauss} also applies to the trimmed mean in~\eqref{eq:trimean}.

The following theorem is the trimmed mean analogue to Theorem~\ref{thm:HDBootstrap} on bootstrap approximation to the distribution of winsorized means. 
\begin{theorem}\label{thm:HDBootstrapTrim}
	Let Assumption~\ref{ass:setting} be satisfied with~$m>2$. Let~$n > 24$. If~$\eps_n,\eps_n'\in(0,1/2)$, with~$\eps_n$ as in~\eqref{eq:epsfam} and~$\eps_n'$ as in~\eqref{eq:epsprime}, and~$\tilde{Z}\sim \mathsf{N}_d(0,\tilde{\Sigma}_n)$ conditionally on~$\tilde{X}_1,\hdots,\tilde{X}_n$, it holds with probability at least~$1-\frac{24}{n}$ that
	\begin{align*}
		\tilde{\rho}_{n,T}:&=\sup_{H\in\mc{H}}\envert[2]{\P\del[1]{S_{n,T}\in H}-\P\del[1]{\tilde{Z}\in H\mid\tilde{X}_1,\hdots,\tilde{X}_n}}\\
		&\leq
		\mathfrak{C}_n
		+C\del[4]{\log^2(d)\sbr[3]{\overline{\eta}_n^{1-\frac{2}{m}}+\del[2]{\frac{\log(dn)}{n}}^{1-\frac{1}{(m/2)\wedge 2}}}}^{1/2}	,
	\end{align*}	
	where~$\mathfrak{C}_n$ is the upper bound on~$\rho_{n,T}$ in~\eqref{eq:trimGauss} of Theorem~\ref{thm:HDGaussTrim} and~$C$ is a constant depending only on~$b_1,b_2,\lambda_1',\lambda_2'$ and~$m$. 
	
	In particular,~$\tilde{\rho}_{n,T}\to 0$ in probability if~$\sqrt{n\log(d)}\overline{\eta}_n^{1-\frac{1}{m}}\to 0$ and~$\log(d)/n^{\frac{m-2}{5m-2}}\to 0$.
\end{theorem}
Theorem~\ref{thm:HDBootstrapTrim} allows for exactly the same growth rates of~$\overline{\eta}_n$ and~$d$ as Theorem~\ref{thm:HDBootstrap} for winsorized means. Note that we have chosen to keep the estimator~$\tilde{\Sigma}_n$, which is based on winsorized means. Of course, one could also use a trimmed mean based estimator obeying the same performance guarantees (in fact, \emph{any} estimator obeying the same performance guarantees as~$\tilde{\Sigma}_n$ in~\eqref{eq:covestimGram} suffices, cf.~also the discussion following~\eqref{eq:Gaussiancomparison}).

\section{Numerical results}\label{sec:sims}
The P-P (probability-probability) plots in Figures~\ref{fig:1} and~\ref{fig:2} below illustrate the theoretical results of the paper. To allow a comparison to existing Gaussian approximations to the distribution of~$S_n=n^{-1/2}\sum_{i=1}^n(X_i-\mu)$, which all are for uncontaminated~$X_i$, we set~$\overline{\eta}_n=0$ throughout. The~$X_{i,j}$ are i.i.d., following a~$t(\nu)$-distribution with~$\nu$ degrees of freedom such that~$\mu=0$ and~$\Sigma=\text{diag}(\nu/(\nu-2),\hdots,\nu/(\nu-2))$. We consider~$\nu=3.01$ (Figure~\ref{fig:1}) and~$\nu=4.01$ (Figure~\ref{fig:2}), corresponding to the~$X_{i,j}$ having all moments strictly less than~$3.01$ and~$4.01$, respectively. Throughout,~$n=200$ and we study~$d\in\cbr[0]{500,5{,}000}$. When implementing~$S_{n,W}$ and~$S_{n,T}$,~$\lambda_1$ and~$\lambda_1'$ are irrelevant as~$\overline{\eta}_n=0$. We set~$\lambda_2=0.1$ throughout. We implement~$\tilde{\Sigma}_n$ with~$\lambda_2'=0.07$ since then the number of observations winsorized in implementing~$\tilde{\Sigma}_n$ is typically the same as in implementing~$S_{n,W}$ and~$S_{n,T}$, i.e.,~$\eps_nn=\lambda_2\log(dn)\approx \lambda_2'\log(d^2n)=\eps_n'n$, which is found to work well. All results are based on~$10{,}000$ replications and we write~$||x||_\infty=\max_{1\leq j\leq d}|x_j|$ for~$x\in\R^d$.

We first discuss how the numerical results illustrate the Gaussian approximations from Theorems~\ref{thm:HDGauss} and~\ref{thm:HDGaussTrim} by comparing $\P(||T_{n}||_\infty\leq t)$,~$t\in\R$ for~$T_n\in\cbr[0]{S_{n,W},S_{n,T},S_n}$ to~$\P(||Z||_\infty\leq t)$ for~$Z\sim\mathsf{N}_d(0,\Sigma)$ (amounting to~$H=[-t,t]^d$ for~$t\in\R$ in these theorems). These results can be found in the leftmost two panels of Figures~\ref{fig:1} and~\ref{fig:2}. For the~$\nu$,~$n$ and~$d$ studied, the trimmed mean of~\cite{resende2024robust} requires trimming more observations than the sample size~$n$ and is hence not implementable (explaining the absent green line in the figures). More generally, even when implementable, that procedure trims an excessive number of observations resulting in inferior practical performance, cf.~the right panel in Figure~\ref{fig:intro} in Section~\ref{sec:Intro}.\footnote{The trimmed mean of \cite{resende2024robust} trims more observations than our implementation of the trimmed mean, which already has lighter tails than its Gaussian approximation, cf.~the left panels of Figures~\ref{fig:1} and~\ref{fig:2}.} For~$\nu=3.01$, Figure~\ref{fig:1} shows that in accordance with Theorem~2.1 of~\cite{kock2024remark}, the Gaussian approximation to the distribution of~$S_n$ can be very imprecise. Crucially, for small~$\alpha\in(0,1)$ it fails exactly at those quantiles~$$c_{1-\alpha}=\inf \cbr[1]{t\in\R: \P(||Z||_\infty\leq t)\geq 1-\alpha },\qquad Z\sim\mathsf{N}_d(0,\Sigma),$$ of the approximating Gaussian, which are natural candidates for critical values of the common test $\varphi_n=\mathds{1}\del[1]{||S_n||_\infty >c_{1-\alpha}}$ [grant that~$\Sigma$ is known]. In fact, one has that $\P(||S_n||_\infty >c_{0.95})\approx 0.39$ for~$d=500$ and~$\P(||S_n||_\infty >c_{0.95})\approx 0.93$ for~$d=5{,}000$ instead of the desired~$0.05=\P(||Z||_\infty>c_{0.95})$.\footnote{The figures plot the cumulative distribution functions (cdf) of~$||S_{n,W}||_\infty$,~$||S_{n,T}||_\infty$, and~$||S_n||_\infty$. Thus, to get the rejection frequencies mentioned in the previous sentence, one calculates one \emph{minus} the cdf.} In conclusion,~$||S_n||_\infty$ has a much \emph{heavier} tail than the approximating~$||Z||_\infty$ and the test~$\varphi_n$ is severely oversized.

On the other hand, as predicted by Theorem~\ref{thm:HDGauss},~$\P(||S_{n,W}||_\infty >c_{0.95})\approx 0.05$ for~$d=500$ and $\P(||S_{n,W}||_\infty >c_{0.95})\approx 0.09$ for~$d=5{,}000$ are much closer to~$0.05=\P(||Z||_\infty>c_{0.95})$. The Gaussian approximation to the trimmed means~$||S_{n,T}||_\infty$, albeit better than that to~$||S_n||_\infty$, is worse than that to the winsorized means~$||S_{n,W}||_\infty$. In particular, the tail of~$||S_{n,T}||_\infty$ is lighter than that of~$||Z||_\infty$.

Next, the middle two panels of Figures~\ref{fig:1} and~\ref{fig:2} illustrate the bootstrap approximation in Theorems~\ref{thm:HDBootstrap} and~\ref{thm:HDBootstrapTrim} by comparing~$\P(||T_{n}||_\infty\leq t)$,~$t\in\R$ for~$T_n\in\cbr[0]{S_{n,W},S_{n,T}}$ to~$\P(||\tilde{Z}||_\infty\leq t\mid X_1,\hdots,X_n)$ [recall that~$\tilde{Z}\sim\mathsf{N}_d(0,\tilde{\Sigma}_n)$ given~$X_1,\hdots,X_n$ with~$\tilde{\Sigma}_n$ from \eqref{eq:tildeSigma}]. These approximations are important as~$\Sigma$ is typically unknown. Defining~$$c_{1-\alpha}^B=\inf \cbr[1]{t\in\R: \P(||\tilde{Z}||_\infty\leq t\mid X_1,\hdots,X_n)\geq 1-\alpha },$$ it is encouraging that~$\P(||S_{n,W}||_\infty >c_{0.95}^B)\approx 0.06$ 
for~$d=500$ and $\P(||S_{n,W}||_\infty >c_{0.95}^B)\approx 0.05$ for~$d=5{,}000$. For~$n=200$ the bootstrap approximation to the distribution of the trimmed mean has a better worst-case performance than for the winsorized mean, but it is less precise in the tails where it is most needed (for testing).

The right two panels of Figures~\ref{fig:1} and~\ref{fig:2} compare~$\P(||S_n||_\infty \leq t)$ to~$\P(||Z^\dagger||_\infty\leq t\mid X_1,\hdots,X_n)$ for~$Z^\dagger\sim\mathsf{N}_d(0,\Sigma_n^\dagger)$ conditionally on~$X_1,\hdots,X_n$ with~$\Sigma_n^\dagger$ being the empirical covariance matrix. Despite the Gaussian approximation to the distribution of~$||S_n||_\infty$ being poor (cf.~the left panels), the bootstrap approximation is passable. This is a result of~$\Sigma_n^\dagger$ ``overestimating''~$\Sigma$, resulting in~$||Z^\dagger||_\infty$ having much heavier tails than~$||Z||_\infty$ and thus the distribution of~$||Z^\dagger||_\infty$ being closer to that of~$||S_n||_\infty$.

Figure~\ref{fig:2} increases~$\nu$ to~$4.01$ from~$3.01$. As expected, this generally improves the performance of all approximations. In conclusion, the Gaussian and bootstrap approximations to the winsorized mean work well in all settings considered. We note, however, that whenever~$n\eps_n$ or~$n\eps_n'$ cross an integer value, the number of winsorized observations changes, which can lead to a ``discontinuity'' in the performance. We also emphasize again that we have investigated scenarios \emph{without} contamination. The normal and bootstrap approximations for the sample mean break down completely under contamination.

\begin{figure}
	\begin{center}
		\footnotesize $\nu=3.01$,~$n=200$ and~$d=500$	
	\end{center}
	
	\vspace{-0.8cm}
	\hspace{-0.8cm}
	\includegraphics[width=5cm]{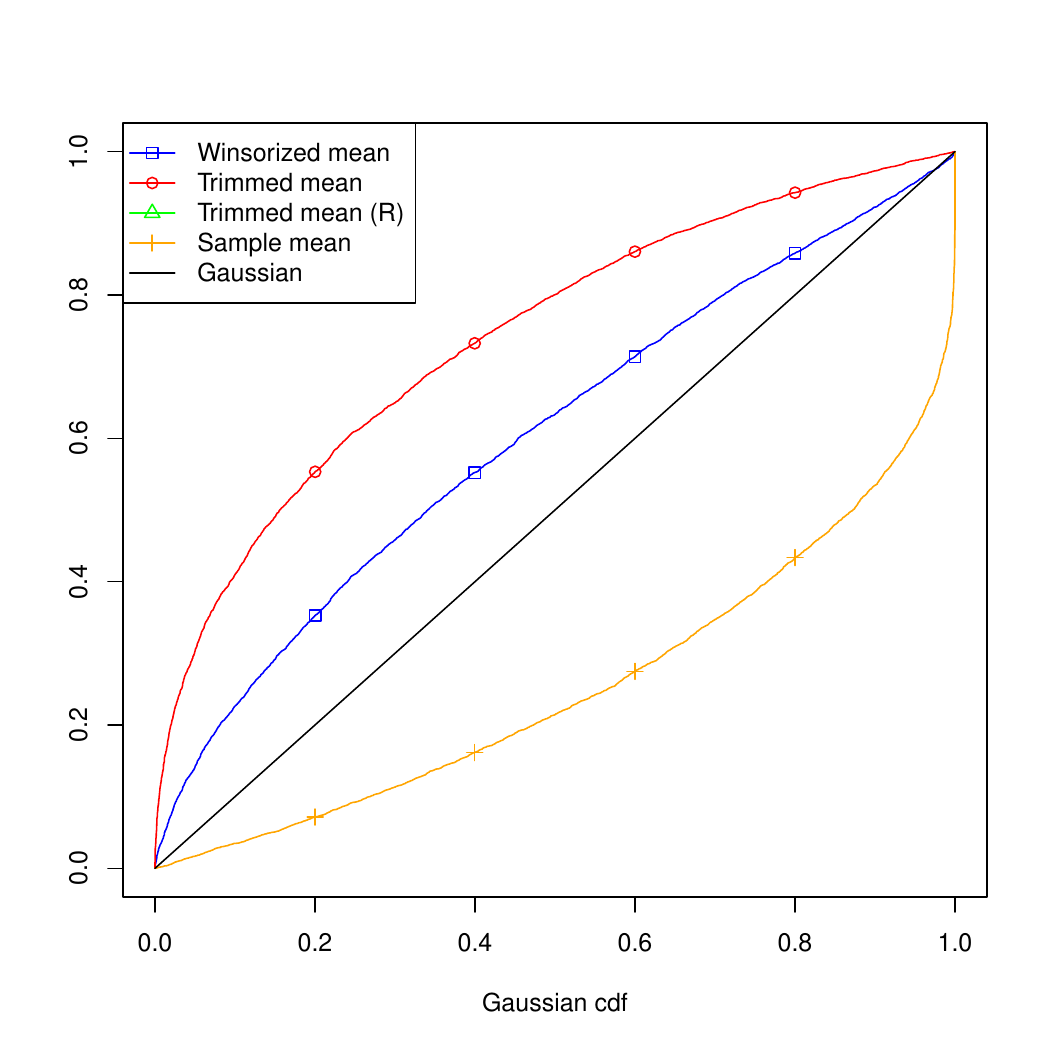}
	\hspace{-0.7cm}
	\includegraphics[width=5cm]{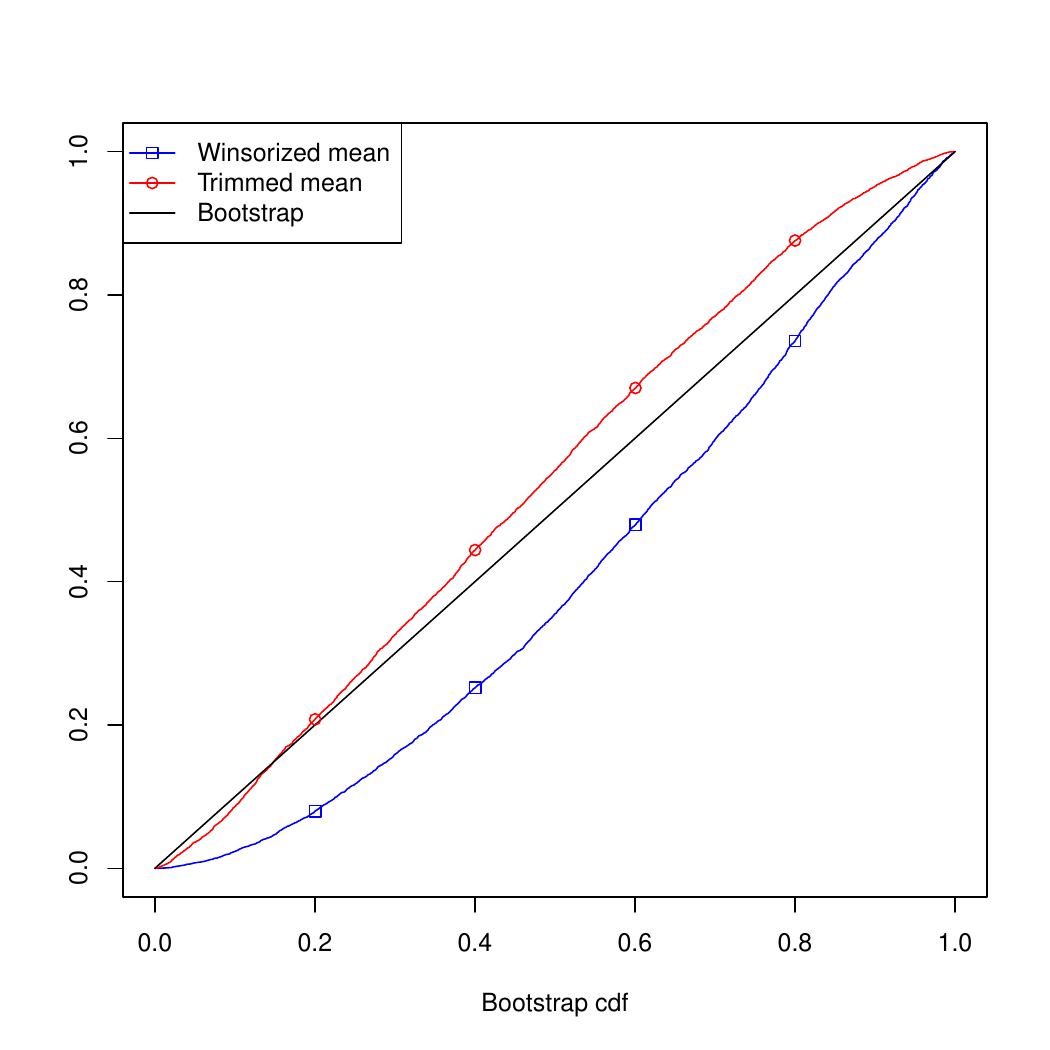}
	\hspace{-0.7cm}
	\includegraphics[width=5cm]{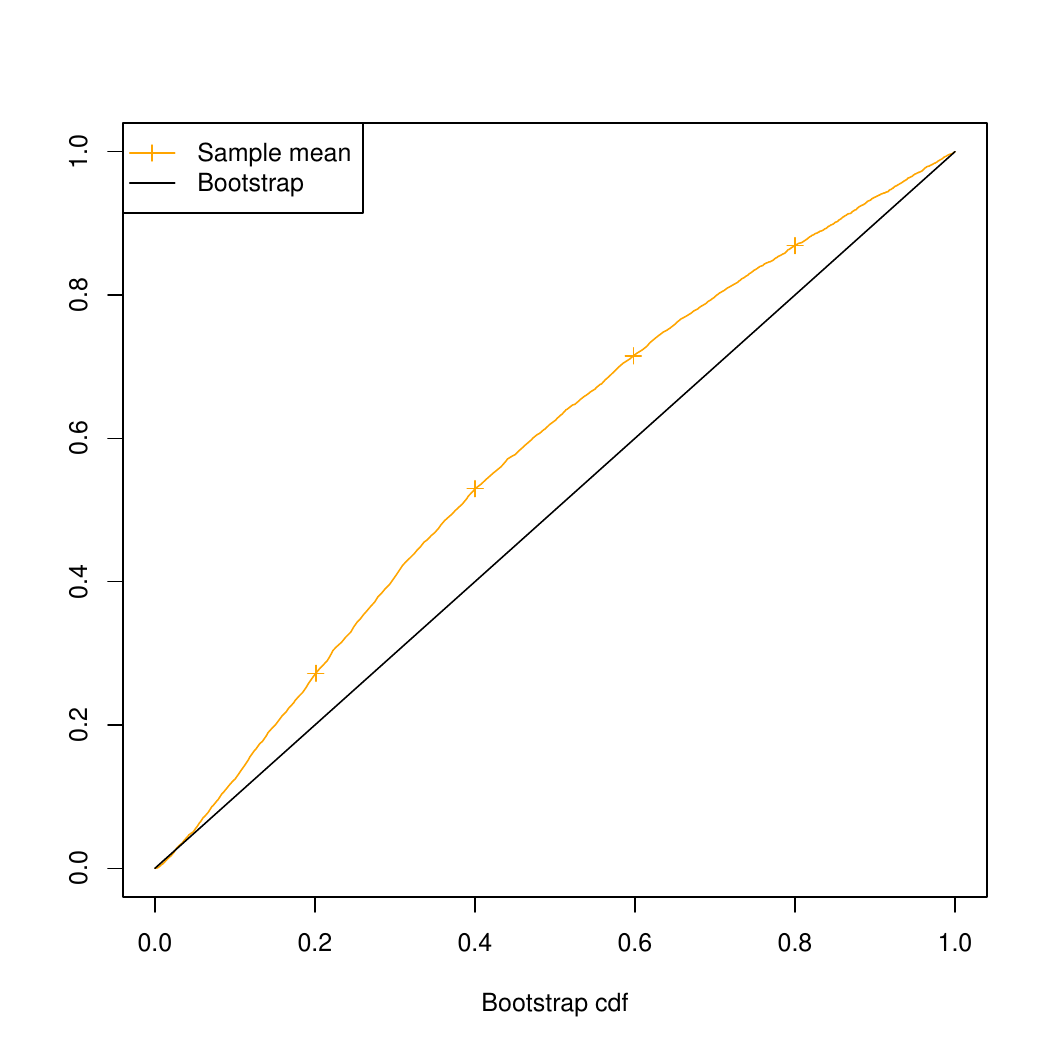}
	
	\begin{center}
		\footnotesize $\nu=3.01$,~$n=200$ and~$d=5{,}000$	
	\end{center}
	\vspace{-0.8cm}
	\hspace{-0.8cm}
	\includegraphics[width=5cm]{Gauss_n=200d=5000df=3.01.pdf}
	\hspace{-0.7cm}
	\includegraphics[width=5cm]{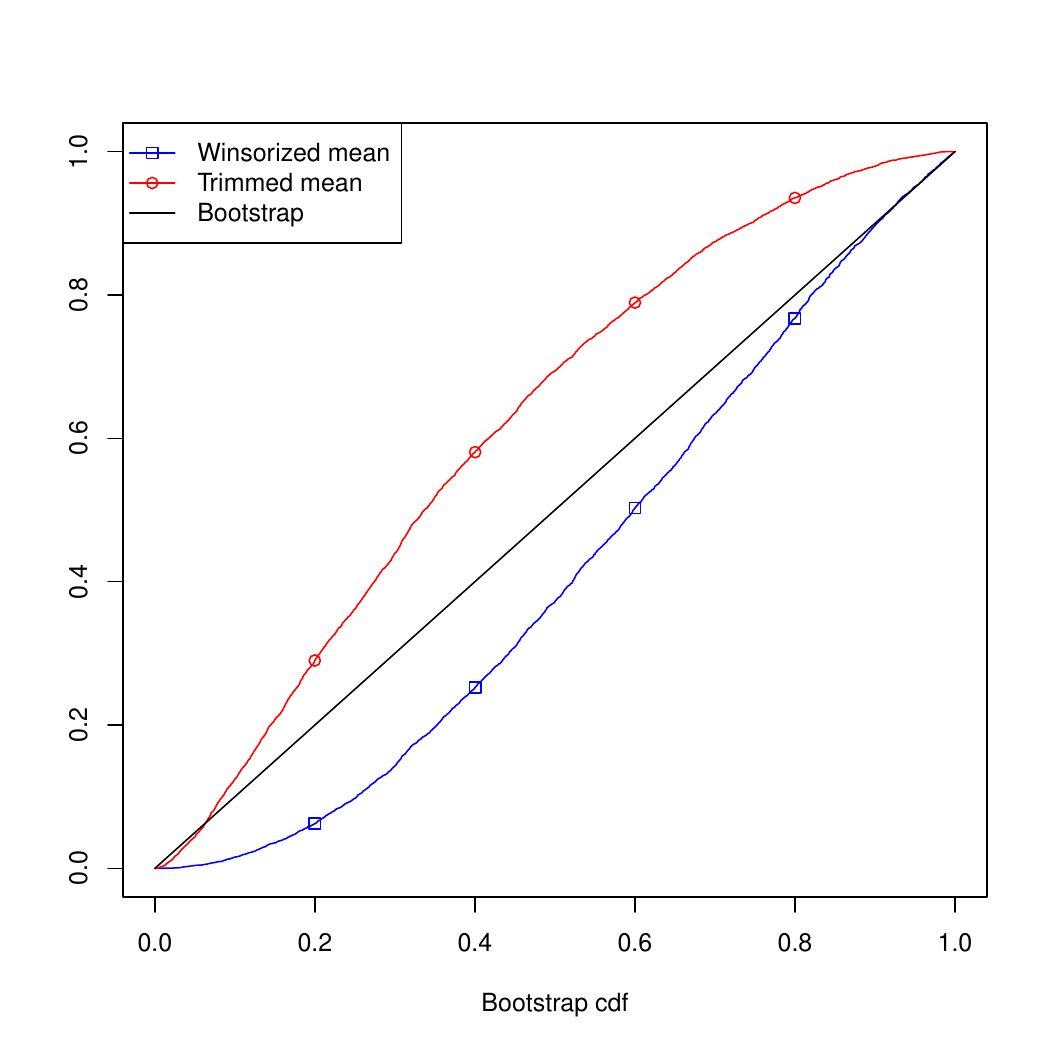}
	\hspace{-0.7cm}
	\includegraphics[width=5cm]{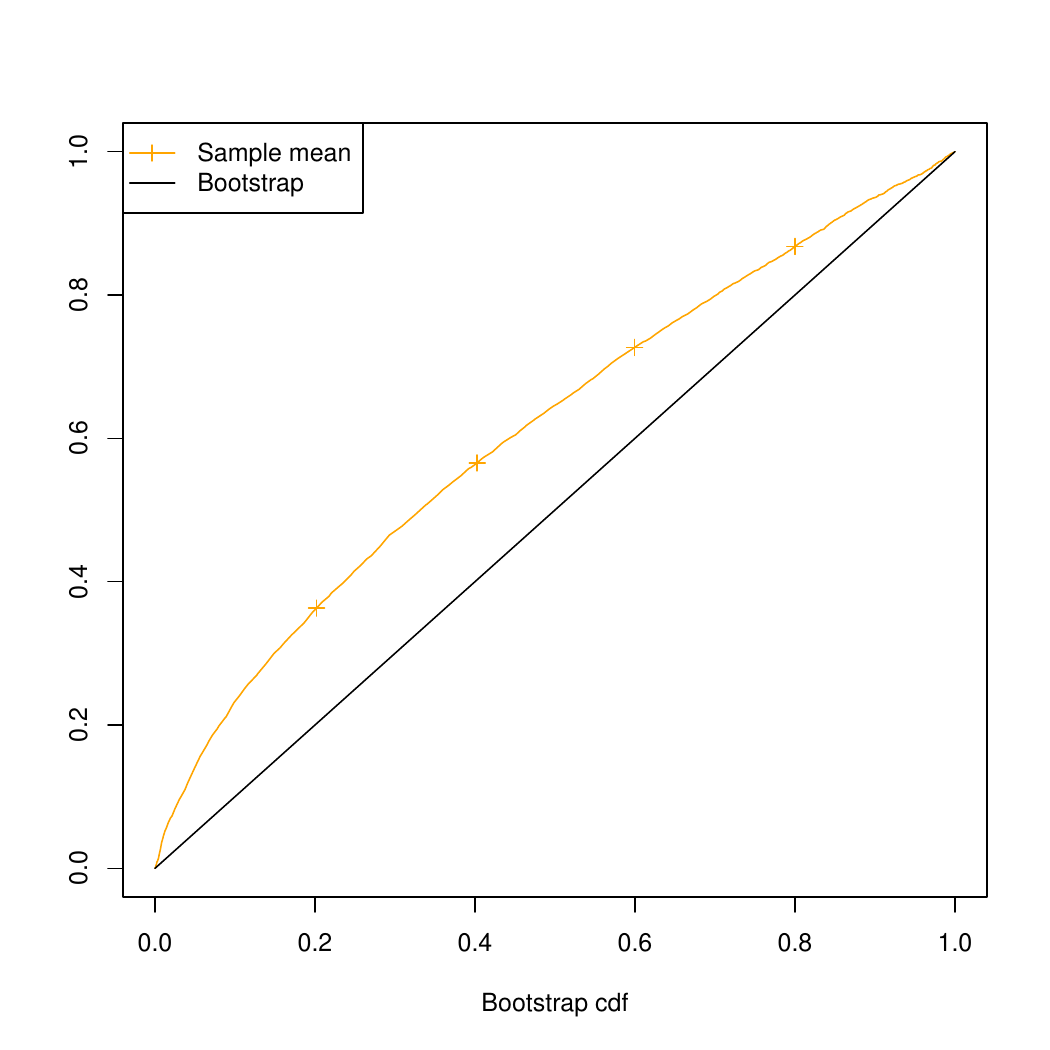}
	
	\caption{\footnotesize P-P (probability-probability) plots: $X_{i,j}\sim$ i.i.d.~$t(3.01)$,~$\overline{\eta}_n=0$. For each pair~$(n,d)$, the left plot illustrates the Gaussian approximations from Theorems~\ref{thm:HDGauss} and~\ref{thm:HDGaussTrim} by comparing $\P(||T_{n}||_\infty\leq t)$,~$t\in\R$ for~$T_n\in\cbr[0]{{\color{blue}S_{n,W}},{\color{red}S_{n,T}},{\color{orange} S_n}}$ to~$\P(||Z||_\infty\leq t)$ for~$Z\sim\mathsf{N}_d(0,\Sigma)$ and~$\Sigma=\text{diag}(2.980,\hdots,2.980)$ the covariance matrix of~$X_1$. This setting corresponds to studying hyperrectangles of the form~$H=[-t,t]^d$ in these theorems. The absent green line indicates that the trimmed mean of~\cite{resende2024robust} requires trimming more than~$n$ observations and is hence not possible to implement for the studied values of~$\nu$,~$n$, and~$d$. For each pair~$(n,d)$, the middle plot illustrates the bootstrap approximations from Theorems~\ref{thm:HDBootstrap} and~\ref{thm:HDBootstrapTrim} by comparing $\P(||T_{n}||_\infty\leq t)$,~$t\in\R$ for~$T_n\in\cbr[0]{S_{n,W},S_{n,T}}$ to~$\P(||\tilde{Z}||_\infty\leq t\mid X_1,\hdots,X_n)$ for~$\tilde{Z}\sim \mathsf{N}_d(0,\tilde{\Sigma}_n)$ conditionally on~$X_1,\hdots,X_n$. Finally, for each pair~$(n,d)$, the right plot illustrates the bootstrap approximations to~$\P(||S_n||_\infty\leq t)$ by~$\P(||Z^\dagger||_\infty\leq t\mid X_1,\hdots,X_n)$ for~$Z^\dagger\sim\mathsf{N}_d(0,\Sigma_n^\dagger)$ conditionally on~$X_1,\hdots,X_n$ and~$\Sigma_n^\dagger$ being the empirical covariance matrix.}
	\label{fig:1}
\end{figure}

\begin{figure}
	\begin{center}
		\footnotesize $\nu=4.01$,~$n=200$ and~$d=500$	
	\end{center}
	\vspace{-0.8cm}
	\hspace{-0.8cm}
	\includegraphics[width=5cm]{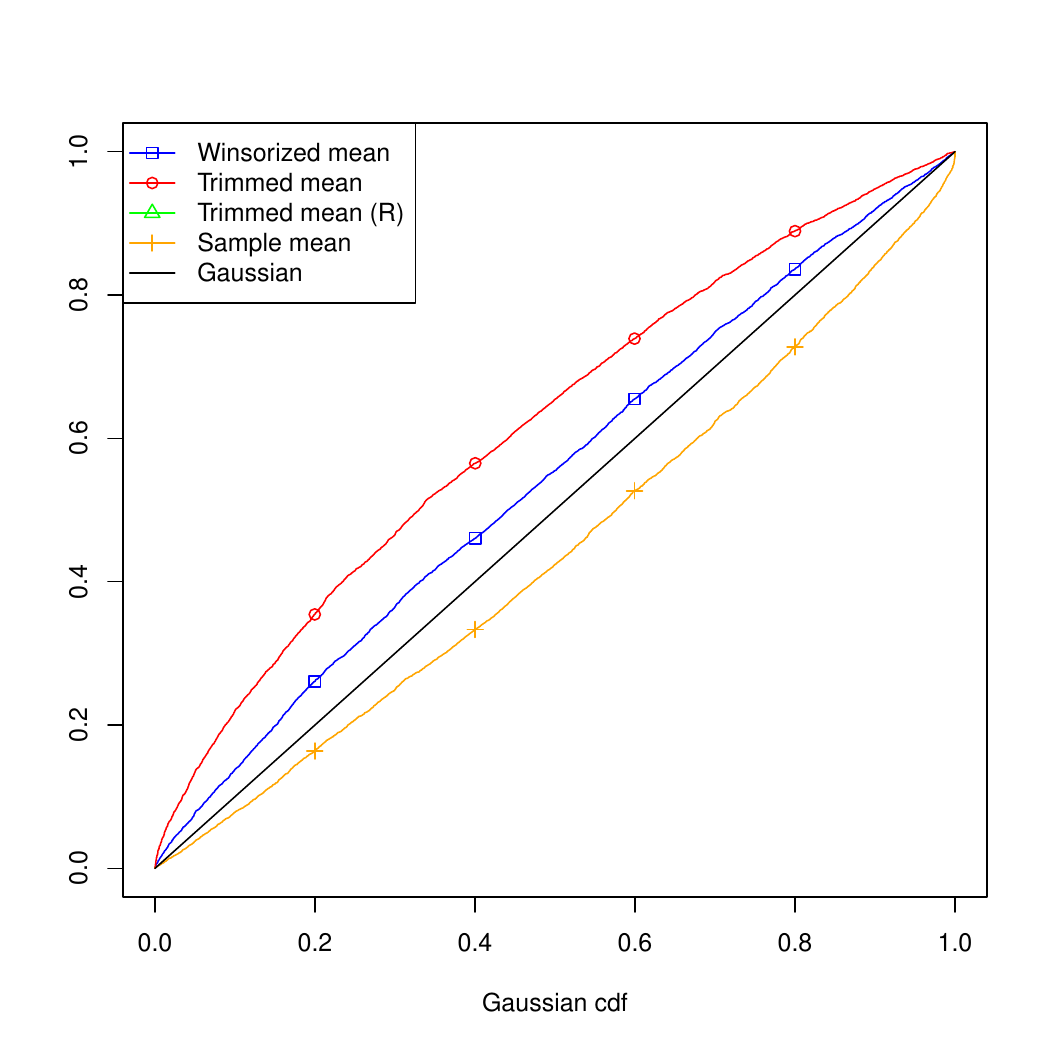}
	\hspace{-0.7cm}
	\includegraphics[width=5cm]{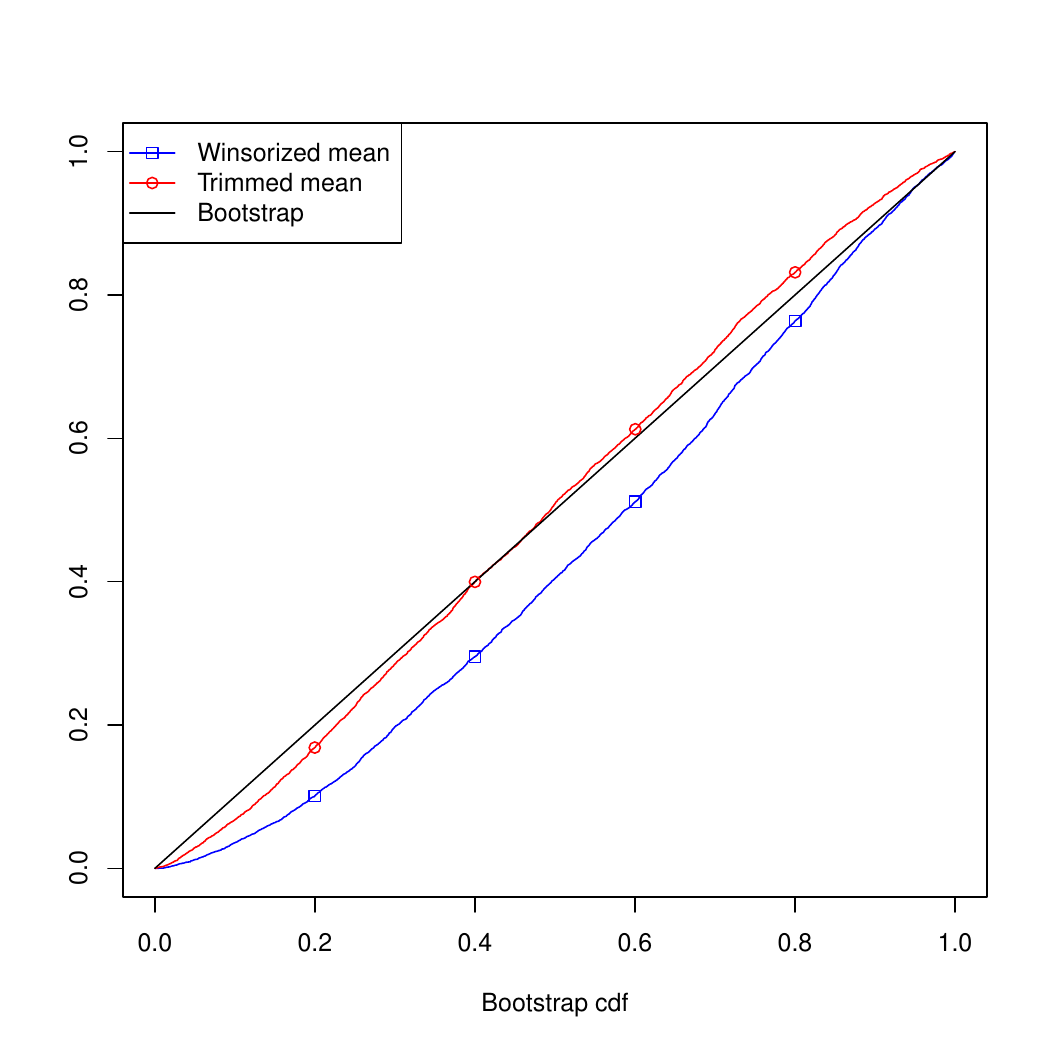}
	\hspace{-0.7cm}
	\includegraphics[width=5cm]{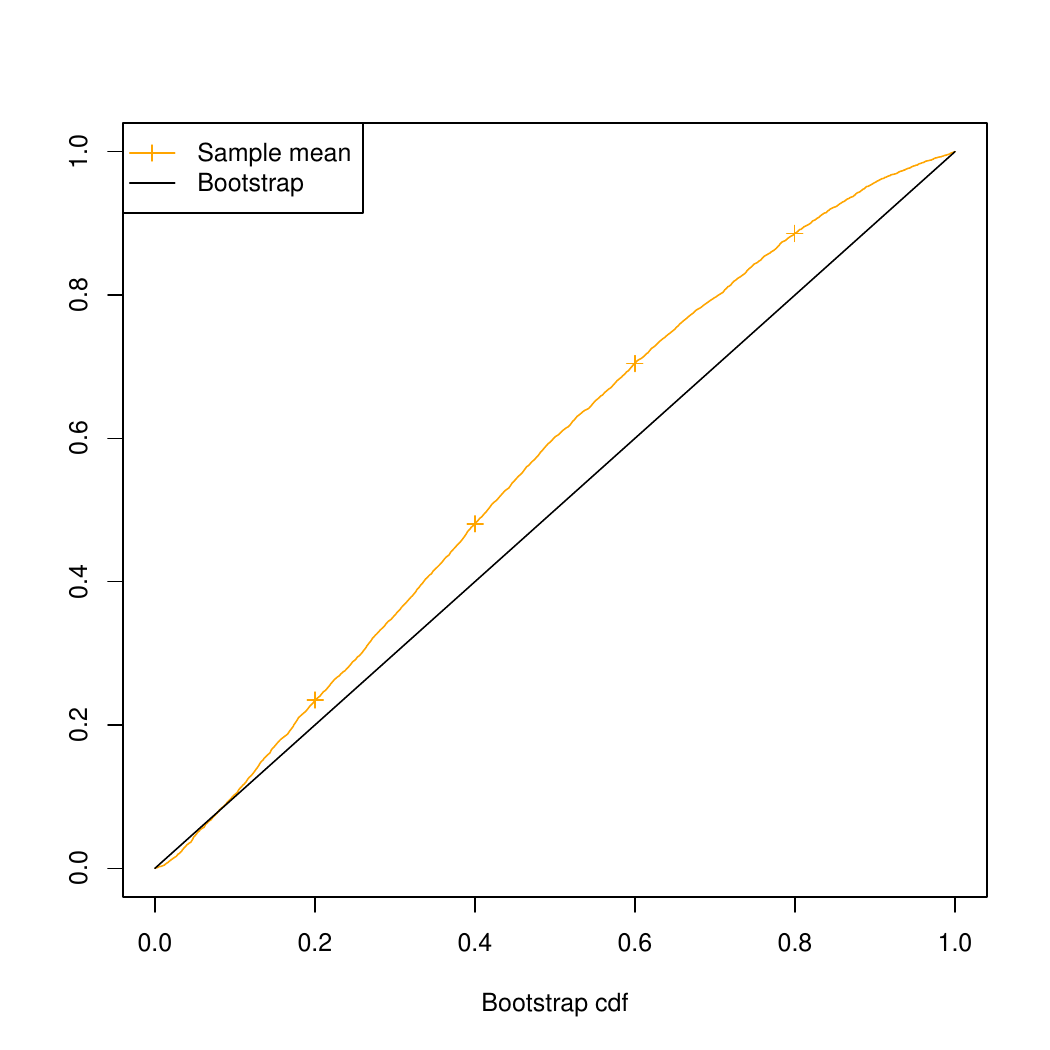}
	
	\begin{center}
		\footnotesize $\nu=4.01$,~$n=200$ and~$d=5{,}000$	
	\end{center}
	\vspace{-0.8cm}
	\hspace{-0.8cm}
	\includegraphics[width=5cm]{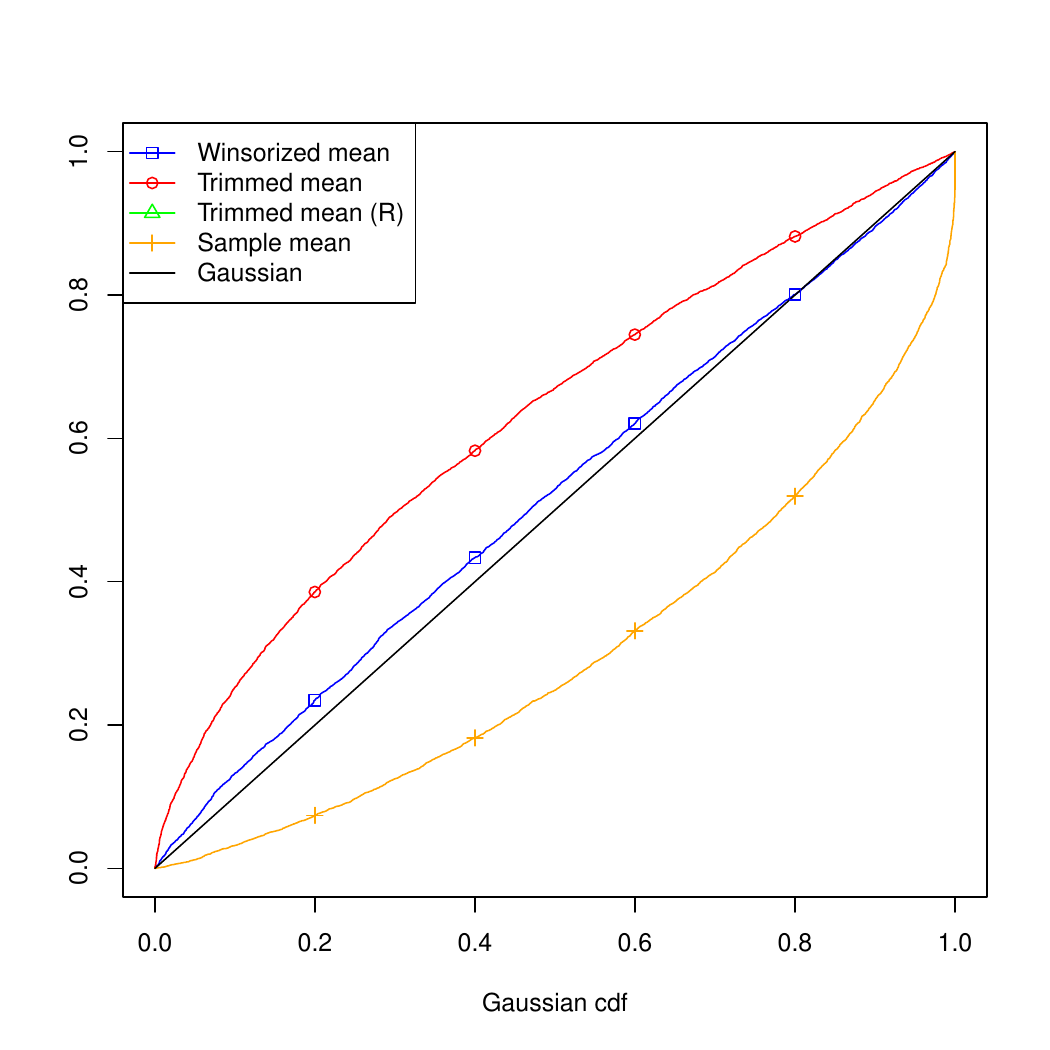}
	\hspace{-0.7cm}
	\includegraphics[width=5cm]{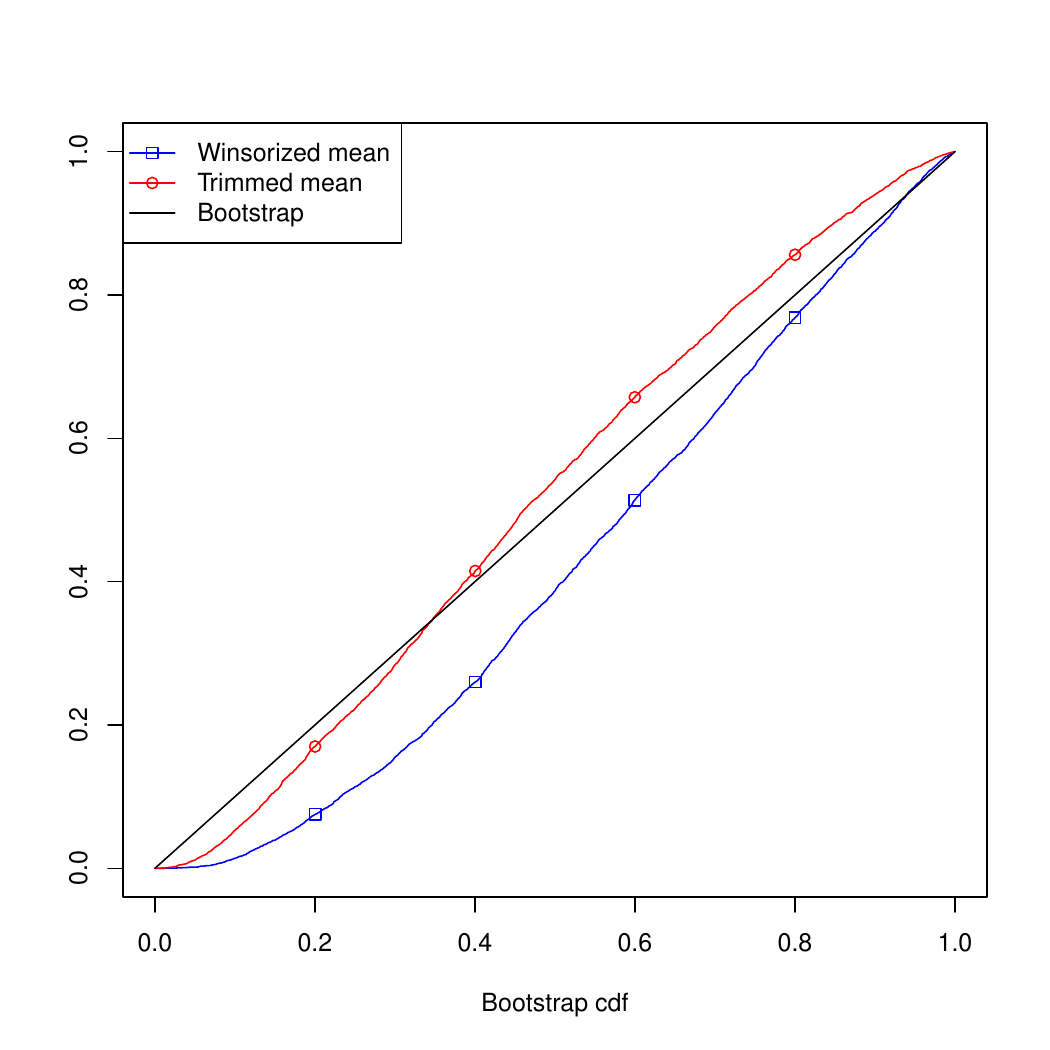}
	\hspace{-0.7cm}
	\includegraphics[width=5cm]{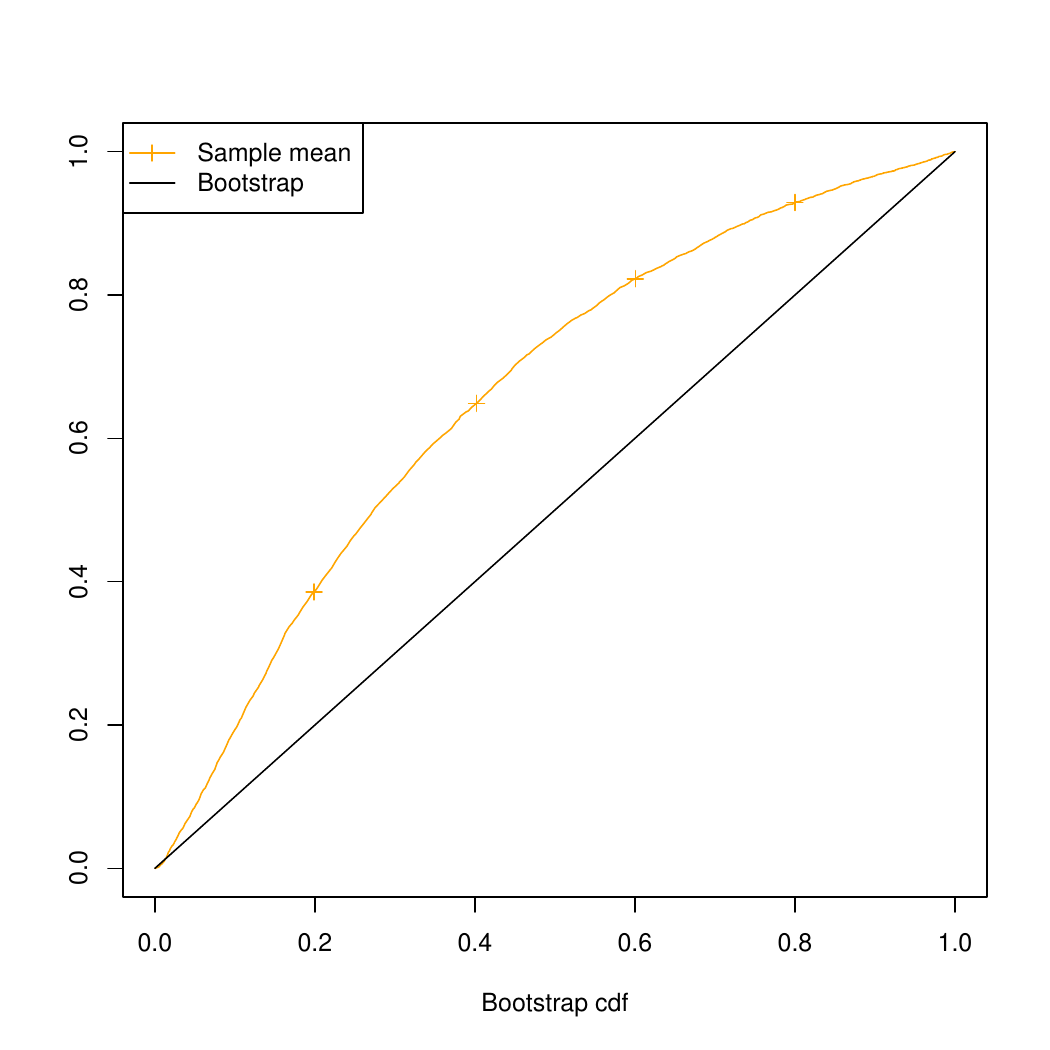}
	\caption{\footnotesize P-P (probability-probability) plots: $X_{i,j}\sim$ i.i.d.~$t(4.01)$,~$\overline{\eta}_n=0$. For each pair~$(n,d)$, the left plot illustrates the Gaussian approximations from Theorems~\ref{thm:HDGauss} and~\ref{thm:HDGaussTrim} by comparing $\P(||T_{n}||_\infty\leq t)$,~$t\in\R$ for~$T_n\in\cbr[0]{{\color{blue}S_{n,W}},{\color{red}S_{n,T}},{\color{orange} S_n}}$ to~$\P(||Z||_\infty\leq t)$ for~$Z\sim\mathsf{N}_d(0,\Sigma)$ and~$\Sigma=\text{diag}(1.995,\hdots,1.995)$ the covariance matrix of~$X_1$. This setting corresponds to studying hyperrectangles of the form~$H=[-t,t]^d$ in these theorems. The absent green line indicates that the trimmed mean of~\cite{resende2024robust} requires trimming more than~$n$ observations and is hence not possible to implement for the studied values of~$\nu$,~$n$, and~$d$. For each pair~$(n,d)$, the middle plot illustrates the bootstrap approximations from Theorems~\ref{thm:HDBootstrap} and~\ref{thm:HDBootstrapTrim} by comparing $\P(||T_{n}||_\infty\leq t)$,~$t\in\R$ for~$T_n\in\cbr[0]{S_{n,W},S_{n,T}}$ to~$\P(||\tilde{Z}||_\infty\leq t\mid X_1,\hdots,X_n)$ for~$\tilde{Z}\sim \mathsf{N}_d(0,\tilde{\Sigma}_n)$ conditionally on~$X_1,\hdots,X_n$. Finally, for each pair~$(n,d)$, the right plot illustrates the bootstrap approximations to~$\P(||S_n||_\infty\leq t)$ by~$\P(||Z^\dagger||_\infty\leq t\mid X_1,\hdots,X_n)$ for~$Z^\dagger\sim\mathsf{N}_d(0,\Sigma_n^\dagger)$ conditionally on~$X_1,\hdots,X_n$ and~$\Sigma_n^\dagger$ being the empirical covariance matrix.}
	\label{fig:2}
\end{figure}

\clearpage

\begin{appendix}
	\numberwithin{equation}{section}

	\section{A useful decomposition}\label{sec:appnotation}

	We begin by outlining a decomposition of~$S_{n,W}$ used in the proof of Theorem~\ref{thm:HDGauss} along with some further notation used in the appendix. For~$p\in(0,1)$ and a random variable~$X$, denote by~$Q_p(X)$ the $p$-quantile of the distribution of~$X$, that is
	\begin{equation*}
		Q_p(X)=\inf\cbr[1]{z\in \R:\P(X\leq z)\geq p}.
	\end{equation*}
	For~$A_+=1-\lambda_1^{-1}\mathds{1}\cbr[0]{\overline{\eta}_n>0}$ and~$A_-=1+\lambda_1^{-1}\mathds{1}\cbr[0]{\overline{\eta}_n>0}$, define (recall~$\varepsilon_n$ from~\eqref{eq:epsfam})
	\begin{equation}\label{eq:c1}
		c_{1,n}:=-A_+W_0\del[1]{-e^{-(\frac{\log(dn)}{\eps_n n}+A_+)/A_+}}\in (0,A_+) 
	\end{equation}
	as well as
	\begin{equation}\label{eq:c2}
		c_{2,n}:=-A_-W_{-1}\del[1]{-e^{-(\frac{\log(dn)}{\eps_n n}+A_-)/A_-}}\in(A_-,\infty);
	\end{equation}
	here~$W_0$ and~$W_{-1}$ denote the principal and lower branch, respectively, of Lambert's~$W$ function (cf., e.g.,~\cite{Corless1996}). Because~$A_+\leq A_-$, note that~$0<c_{1,n}<c_{2,n}<\infty$. Next, for~$j=1,\hdots,d$, let
	\begin{equation}\label{eq:alphas}
		\underline{\alpha}_{j}:=Q_{c_{1,n}\eps_n}(X_{1,j})\qquad\text{and}\qquad 	\overline{\alpha}_j:= Q_{c_{2,n}\eps_n}(X_{1,j})
	\end{equation}
	as well as
	\begin{equation}\label{eq:betas}
		\underline{\beta}_j:=Q_{1-c_{2,n}\eps_n
		}(X_{1,j})\qquad \text{and}
		\qquad\overline{\beta}_j:= Q_{1-c_{1,n}\eps_n}(X_{1,j}), 
	\end{equation}
	In the \emph{appendix only} we shall sometimes impose
	\begin{align}\label{eq:epscond}
		\mathfrak{T}_n^{(1)}:=2\eps_n +\frac{\log(dn)}{n}+\sqrt{\del[2]{\frac{\log(dn)}{n}}^2+4\frac{\log(dn)}{n}\eps_n}<\frac{1}{2},
	\end{align}
	which, however, is \emph{not} imposed for any of our Gaussian and bootstrap approximations in this paper, as these will be argued to hold trivially when~\eqref{eq:epscond} is not satisfied. By Lemma~\ref{lem:cControl} one has that~\eqref{eq:epscond} ensures that
	\begin{align*}
		0
		<
		\eps_n\min(c_{1,n},c_{2,n})
		\leq 
		\eps_n(c_{1,n}+c_{2,n})<1/2,
	\end{align*}
	implying that~$\underline{\alpha}_{j}\leq\underline{\beta}_j$ and~$\overline{\alpha}_j\leq \overline{\beta}_j$, such that the decompositions in~\eqref{eq:decomp} and~\eqref{eq:decomp2} below are well-defined. 
	
	By definition,~$\hat{\alpha}_j=\tilde X_{\lceil \eps_n n \rceil,j}^*$ and $\hat{\beta}_j=\tilde X_{\lfloor(1-\eps_n )n\rfloor+1,j}^*$. Lemma~\ref{lem:quantiles} and the union bound, together with obvious monotonicity properties of~$(a,b) \mapsto \phi_{a,b}$, show that with probability at least~$1-\frac{4}{n}$ we have simultaneously for~$j=1,\hdots,d$ (the inequalities $\phi_{\underline{\alpha}_j, \underline{\beta}_j} \leq \phi_{\hat{\alpha}_j, \hat{\beta}_j} \leq \phi_{\overline{\alpha}_j, \overline{\beta}_j}$ and hence)
	\begin{equation}\label{eq:lbub}
		\frac{1}{n}\sum_{i=1}^n\sbr[1]{\phi_{\underline{\alpha}_j,\underline{\beta}_{j}}(\tilde{X}_{i,j})-\mu_j}
		\leq
		\frac{1}{n}\sum_{i=1}^n\sbr[1]{\phi_{\hat\alpha_j,\hat\beta_j}(\tilde{X}_{i,j})-\mu_j}
		\leq
		\frac{1}{n}\sum_{i=1}^n\sbr[1]{\phi_{\overline{\alpha}_j,\overline{\beta}_j}(\tilde{X}_{i,j})-\mu_j}.
	\end{equation}
	The far right-hand side of~\eqref{eq:lbub} can be decomposed as 
	\begin{align}
		\frac{1}{n}\sum_{i=1}^n\sbr[1]{\phi_{\overline{\alpha}_j,\overline{\beta}_j}(\tilde{X}_{i,j})-\mu_j}
		&=
		\underbrace{\frac{1}{n}\sum_{i=1}^n\sbr[1]{\phi_{\overline{\alpha}_j,\overline{\beta}_j}(\tilde{X}_{i,j})-\phi_{\overline{\alpha}_j,\overline{\beta}_j}(X_{i,j})}}_{\overline I_{n,j,1}}\notag\\
		&+
		\underbrace{\frac{1}{n}\sum_{i=1}^n\sbr[1]{\phi_{\overline{\alpha}_j,\overline{\beta}_j}(X_{i,j})-\E\phi_{\overline{\alpha}_j,\overline{\beta}_j}(X_{i,j})}}_{\overline I_{n,j,2}}\notag\\
		&+\underbrace{\frac{1}{n}\sum_{i=1}^n\sbr[1]{\E\phi_{\overline{\alpha}_j,\overline{\beta}_j}(X_{i,j})-\mu_j}}_{\overline I_{n,j,3}}\label{eq:decomp}.
	\end{align}
	Similarly, the left-hand side of~\eqref{eq:lbub} can be decomposed as
	\begin{equation}\label{eq:decomp2}
		\frac{1}{n}\sum_{i=1}^n\sbr[1]{\phi_{\underline{\alpha}_j,\underline{\beta}_j}(\tilde{X}_{i,j})-\mu_j}
		=
		\underline I_{n,j,1}+\underline I_{n,j,2}+\underline I_{n,j,3},
	\end{equation}
	with~$\underline I_{n,j,k}$, $j=1,\hdots, d,\ k=1,2,3$, defined analogously to the~$\overline I_{n,j,k}$. Define
	\begin{equation}\label{eq:I_13_max}
		\overline I_{n,k}=\max_{j=1,\hdots, d}\envert[1]{\overline I_{n,j,k}}\qquad\text{and}\qquad \underline I_{n,k}=\max_{j=1,\hdots, d}\envert[1]{\underline I_{n,j,k}},\qquad k=1,3,
	\end{equation}
	as well as
	\begin{equation*}
		\overline I_{n,2}=\del[1]{\overline I_{n,1,2},\hdots,\overline I_{n,d,2}}'\qquad\text{and}\qquad \underline I_{n,2}=\del[1]{\underline I_{n,1,2},\hdots,\underline I_{n,d,2}}'.	
	\end{equation*}
	Throughout, $\overline{\R}=\R\cup\cbr[0]{-\infty,\infty}$. Thus, writing~$Y_{n,j}=n^{-1/2}\sum_{i=1}^n\sbr[1]{\phi_{\hat\alpha_j,\hat\beta_j}(\tilde{X}_{i,j})-\mu_j}$ and~$Y_n=\del[1]{Y_{n,1},\hdots,Y_{n,d}}'$, one has with probability at least~$1-\frac{4}{n}$ that
	\begin{equation}\label{eq:LBUPHD}
		\sqrt{n}(\underline{I}_{n,2}-\underline{I}_{n,1}-\underline{I}_{n,3})
		\leq 
		Y_n
		\leq 
		\sum_{k=1}^3\sqrt{n}\overline I_{n,k},
	\end{equation}
	where \emph{here and in the sequel} (i) all inequalities between vectors are understood elementwise, and (ii) with some abuse of notation we define, for a vector~$x \in \overline{\R}^m$, say, and a real number~$a$, the sum~$x + a$ coordinate-wise as~$(x_1 + a, \hdots, x_m + a)'$ , with the usual convention that~$\infty + a = \infty$ and~$-\infty + a = -\infty$. In the next section we study the left- and right-hand sides of the previous display.
	
	\section{Preparatory lemmas}\label{sec:preplem}
	Let~$\varepsilon_n$ be as in~\eqref{eq:epsfam} and recall the quantities defined in~\eqref{eq:alphas} and~\eqref{eq:betas}. Let~$\overline{\Sigma}_{\eps_n}$ and~$\underline{\Sigma}_{\eps_n}$ be the matrices with~$(j,k)$th entry
	\begin{equation*}
		\overline{\Sigma}_{\eps_n,j,k}
		=
		\E\sbr[2]{\del[1]{\phi_{\overline{\alpha}_j,\overline{\beta}_j}(X_{1,j})-\E \phi_{\overline{\alpha}_j,\overline{\beta}_j}(X_{1,j})}\del[1]{\phi_{\overline{\alpha}_k,\overline{\beta}_k}(X_{1,k})-\E \phi_{\overline{\alpha}_k,\overline{\beta}_k}(X_{1,k}}} 
	\end{equation*}
	and
	\begin{equation*}
		\underline{\Sigma}_{\eps_n,j,k}
		=
		\E\sbr[2]{\del[1]{\phi_{\underline{\alpha}_j,\underline{\beta}_j}(X_{1,j})-\E \phi_{\underline{\alpha}_j,\underline{\beta}_j}(X_{1,j})}\del[1]{\phi_{\underline{\alpha}_k,\underline{\beta}_k}(X_{1,k})-\E \phi_{\underline{\alpha}_k,\underline{\beta}_k}(X_{1,k}}}, 
	\end{equation*}
	respectively. The following lemma bounds the distance of these covariance matrices to the covariance matrix~$\Sigma$ of the vector~$X_1$. 
	\begin{lemma}\label{lem:trunccov}
		Let Assumption~\ref{ass:setting} be satisfied with~$m>2$. If~$\eps_n$ is chosen as in~\eqref{eq:epsfam} satisfying~\eqref{eq:epscond}, then, for a positive constant~$C$ depending only on~$\lambda_2$ and~$m$, it holds that
		\begin{equation}\label{eq:truncvar}
			\max_{j,k\in\cbr[0]{1,\hdots,d}}\del[2]{\envert[1]{\overline{\Sigma}_{\eps_n,j,k}-\Sigma_{j,k}}\vee\envert[1]{\underline{\Sigma}_{\eps_n,j,k}-\Sigma_{j,k}}}
			\leq
			C\sigma_m^2\eps_n^{1-\frac{2}{m}}=:\mathfrak{T}_n^{(2)}.
		\end{equation}
	\end{lemma}
	
	\begin{proof}
		We only establish~\eqref{eq:truncvar} for~$\max_{j,k\in\cbr[0]{1,\hdots,d}}\envert[1]{\overline{\Sigma}_{\eps_n,j,k}-\Sigma_{j,k}}$ as the proof is identical for~$\max_{j,k\in\cbr[0]{1,\hdots,d}}\envert[1]{\underline{\Sigma}_{\eps_n,j,k}-\Sigma_{j,k}}$. 
		
		Fix~$1\leq j\leq d$ and note that
		\begin{align*}
			&\phi_{\overline{\alpha}_j,\overline{\beta}_j}(X_{1,j})-\E \phi_{\overline{\alpha}_j,\overline{\beta}_j}(X_{1,j}) \\
			&=
			X_{1,j}+(\overline{\alpha}_j-X_{1,j})\mathds{1}\del[1]{X_{1,j}<\overline{\alpha}_j}+	(\overline{\beta}_j-X_{1,j})\mathds{1}\del[1]{X_{1,j}>\overline{\beta}_j}\\
			&-\E \del[1]{ X_{1,j}+(\overline{\alpha}_j-X_{1,j})\mathds{1}\del[1]{X_{1,j}<\overline{\alpha}_j}+	(\overline{\beta}_j-X_{1,j})\mathds{1}\del[1]{X_{1,j}>\overline{\beta}_j}}\\
			&=\sbr[1]{X_{1,j}-\E X_{1,j}}+\sbr[1]{(\overline{\alpha}_j-X_{1,j})\mathds{1}\del[1]{X_{1,j}<\overline{\alpha}_j}-\E (\overline{\alpha}_j-X_{1,j})\mathds{1}\del[1]{X_{1,j}<\overline{\alpha}_j}}\\
			&+\sbr[1]{(\overline{\beta}_j-X_{1,j})\mathds{1}\del[1]{X_{1,j}>\overline{\beta}_j}-\E (\overline{\beta}_j-X_{1,j})\mathds{1}\del[1]{X_{1,j}>\overline{\beta}_j}}\\
			&	=:T_{1,j}+T_{2,j}+T_{3,j}.
		\end{align*}
		Thus, for~$j,k\in\cbr[0]{1,\hdots,d}$, it follows from the Cauchy-Schwarz inequality that
		\begin{align}
			\envert[1]{\overline{\Sigma}_{\eps_n,j,k}-\Sigma_{j,k}}
			&=
			\envert[1]{\E(T_{1,j}+T_{2,j}+T_{3,j})(T_{1,k}+T_{2,k}+T_{3,k})-\E (T_{1,j}T_{1,k})}\notag\\
			&\leq
			|\E T_{1,j}T_{2,k}|+|\E T_{1,j}T_{3,k}|+|\E T_{2,j}T_{1,k}|+|\E T_{2,j}T_{2,k}|\notag\\
			&+
			|\E T_{2,j}T_{3,k}|+|\E T_{3,j}T_{1,k}|+|\E T_{3,j}T_{2,k}|+|\E T_{3,j}T_{3,k}|\notag\\
			&\leq
			|\E T_{1,j}T_{2,k}|+|\E T_{1,j}T_{3,k}|+|\E T_{2,j}T_{1,k}|+|\E T_{3,j}T_{1,k}|\notag\\
			&+
			\del[1]{\E T_{2,j}^2\E T_{2,k}^2}^{1/2}+\del[1]{\E T_{2,j}^2\E T_{3,k}^2}^{1/2}+\del[1]{\E T_{3,j}^2\E T_{2,k}^2|}^{1/2}+\del[1]{\E T_{3,j}^2\E T_{3,k}^2}^{1/2}.\label{eq:decompcov}
		\end{align}
		We proceed by bounding the far right-hand side of~\eqref{eq:decompcov}. To this end, note first that for all~$j=1,\hdots,d$ it holds that 
		\begin{align*}
			\E T_{2,j}^2
			&=
			\E\sbr[1]{(\overline{\alpha}_j-X_{1,j})\mathds{1}\del[1]{X_{1,j}<\overline{\alpha}_j}-\E (\overline{\alpha}_j-X_{1,j})\mathds{1}\del[1]{X_{1,j}<\overline{\alpha}_j}}^2\\
			&\leq
			\E(\overline{\alpha}_j-X_{1,j})^2\mathds{1}\del[1]{X_{1,j}<\overline{\alpha}_j}.
		\end{align*}
		Recall that~$\overline{\alpha}_j=Q_{c_{2,n}\eps_n}(X_{1,j})$ such that~$\P(X_{1,j}<\overline{\alpha}_j)\leq c_{2,n}\eps_n$. Thus, it follows by Lemma~\ref{lem:quantile_mean}, followed by H{\"o}lder's inequality, that there exists a positive constant~$C$ depending only on~$\lambda_2$ and~$m$ such that
		\begin{align}
			\E T_{2,j}^2
			&\leq
			\E(\mu_j+\sigma_m/(1-c_{2,n}\eps_n)^{1/m}-X_{1,j})^2\mathds{1}\del[1]{X_{1,j}<\overline{\alpha}_j}\notag\\
			&\leq 
			\del[2]{\E\envert[1]{\mu_j-X_{1,j}+\sigma_m/(1-c_{2,n}\eps_n)^{1/m}}^m}^{2/m}\P(X_{1,j}<\overline{\alpha}_j)^{1-\frac{2}{m}}\notag\\
			&\leq
			\del[2]{2^m\sigma_m^m+2^m\sigma_m^m/(1-c_{2,n}\eps_n)}^{2/m}(c_{2,n}\eps_n)^{1-\frac{2}{m}}\notag\\
			&\leq C\sigma_m^2\eps_n^{1-\frac{2}{m}},\label{eq:boundT2sq}
		\end{align}
		the last inequality using that~$c_{2,n}\eps_n<1/2$ by Lemma~\ref{lem:cControl} and~\eqref{eq:epscond} and the same lemma also bounding~$c_{2,n}$ from above in terms of~$\lambda_2$ only.
		
		By the same arguments, using~$\P(X_{1,j}>\overline{\beta}_j)=1-\P(X_{1,j}\leq \overline{\beta}_j)\leq c_{1,n}\eps_n$ and adjusting~$C$, one also has that
		\begin{equation}\label{eq:boundT3sq}
			\E T_{3,j}^2
			\leq
			C\sigma_m^2\eps_n^{1-\frac{2}{m}}.
		\end{equation}
		Furthermore,
		\begin{align*}
			\E T_{1,j}T_{2,k}
			&=
			\E\del[2]{[X_{1,j}-\mu_j]\sbr[1]{(\overline{\alpha}_k-X_{1,k})\mathds{1}\del[1]{X_{1,k}<\overline{\alpha}_k}-\E (\overline{\alpha}_k-X_{1,k})\mathds{1}\del[1]{X_{1,k}<\overline{\alpha}_k}}}\\
			&=
			\E\del[1]{[X_{1,j}-\mu_j](\overline{\alpha}_k-X_{1,k})\mathds{1}\del[1]{X_{1,k}<\overline{\alpha}_k}}.	
		\end{align*}
		Recalling that~$\overline{\alpha}_j=Q_{c_{2,n}\eps_n}(X_{1,j})$, it follows by Lemma~\ref{lem:quantile_mean} and H{\"o}lder's inequality, that
		\begin{align}
			&\envert[1]{\E T_{1,j}T_{2,k}} \notag\\
			&\leq
			\E\del[1]{|X_{1,j}-\mu_j||\mu_k+\sigma_m/(1-c_{2,n}\eps_n)^{1/m}-X_{1,k}|\mathds{1}\del[1]{X_{1,k}<\overline{\alpha}_k}}\notag\\
			&\leq
			\E\del[1]{|X_{1,j}-\mu_j||X_{1,k}-\mu_k|+|X_{1,j}-\mu_j|\sigma_m/(1-c_{2,n}\eps_n)^{1/m}}\mathds{1}\del[1]{X_{1,k}<\overline{\alpha}_k}\notag\\
			&\leq
			\del[2]{\E\del[1]{|X_{1,j}-\mu_j||X_{1,k}-\mu_k|+|X_{1,j}-\mu_j|\sigma_m/(1-c_{2,n}\eps_n)^{1/m}}^{\frac{m}{2}}}^{\frac{2}{m}}(c_{2,n}\eps_n)^{1-\frac{2}{m}}\notag\\
			&\leq
			\del[2]{2^{m/2}(\sigma_m^m+\sigma_{m/2}^{m/2}\sigma_m^{m/2}/(1-c_{2,n}\eps_n)^{1/2}}^{2/m}(c_{2,n}\eps_n)^{1-\frac{2}{m}}\notag\\
			&\leq
			C\sigma_m^2\eps_n^{1-\frac{2}{m}},\label{eq:boundT1T2}
		\end{align}
		the last inequality using that~$c_{2,n}\eps_n<1/2$ by Lemma~\ref{lem:cControl} and~\eqref{eq:epscond} and the same lemma also bounding~$c_{2,n}$ from above in terms of~$\lambda_2$ only (adjust~$C$ if needed).
		
		By the same arguments,
		\begin{equation}
			\envert[1]{\E T_{1,j}T_{3,k}}
			\leq
			C\sigma_m^2\eps_n^{1-\frac{2}{m}}.\label{eq:boundT1T3}	
		\end{equation}
		Therefore, observing that none of the upper bounds in~\eqref{eq:boundT2sq}--\eqref{eq:boundT1T3} depend on~$j$ and~$k$ it follows that (adjust~$C$ if needed)
		\begin{equation*}
			\max_{j,k\in\cbr[0]{1,\hdots,d}}\envert[1]{\overline{\Sigma}_{\eps_n,j,k}-\Sigma_{j,k}}
			\leq 
			C\sigma_m^2\eps_n^{1-\frac{2}{m}}.
		\end{equation*}
	\end{proof}

	The following lemma applies a high-dimensional Gaussian approximation to~$\underline{I}_{n,2}$ and~$\overline{I}_{n,2}$, which are both sums of winsorized, and hence bounded, random variables. 
	\begin{lemma}\label{lem:HDGaussBD}
		Let Assumption~\ref{ass:setting} be satisfied with~$m>2$. Let~$\eps_n$ be chosen as in~\eqref{eq:epsfam} satisfying~\eqref{eq:epscond} and suppose that there exists a strictly positive constant~$\mathfrak{b}_1$ such that for all~$j=1,\hdots,d$
		\begin{equation*}
			\E\sbr[1]{\phi_{\underline{\alpha}_j,\underline{\beta}_j}(X_{1,j})-\E\phi_{\underline{\alpha}_j,\underline{\beta}_j}(X_{1,j})}^2>\mathfrak{b}_1^2\quad\text{and}\quad\E\sbr[1]{\phi_{\overline{\alpha}_j,\overline{\beta}_j}(X_{1,j})-\E\phi_{\overline{\alpha}_j,\overline{\beta}_j}(X_{1,j})}^2>\mathfrak{b}_1^2.
		\end{equation*}
		Then, for~$\overline{Z}\sim \mathsf{N}_d(0,\overline{\Sigma}_{\eps_n})$ and~$\underline{Z}\sim \mathsf{N}_d(0,\underline{\Sigma}_{\eps_n})$,
		\begin{equation*}
			\sup_{H\in \mc{H}}\envert[3]{\P\del[3]{\frac{1}{n^{1/2}}\sum_{i=1}^n\sbr[1]{\phi_{\overline{\alpha}_j,\overline{\beta}_j}(X_{i,j})-\E\phi_{\overline{\alpha}_j,\overline{\beta}_j}(X_{i,j})}\in H}-\P\del[1]{\overline{Z}\in H}}
		\end{equation*}
		and
		\begin{equation*}
			\sup_{H\in\mc{H}}\envert[3]{\P\del[3]{\frac{1}{n^{1/2}}\sum_{i=1}^n\sbr[1]{\phi_{\underline{\alpha}_j,\underline{\beta}_j}(X_{i,j})-\E\phi_{\underline{\alpha}_j,\underline{\beta}_j}(X_{i,j})}\in H}-\P\del[1]{\underline{Z}\in H}}
		\end{equation*}
		are bounded from above by
		\begin{equation*}
			C\sigma_m^{1/2}\del[3]{\frac{\log^{5-\frac{2}{m}}(dn)}{n^{1-\frac{2}{m}}}}^{\frac{1}{4}},	
		\end{equation*}
		where~$C$ is a constant depending only on~$\mathfrak{b}_1$,~$b_2$,~$\lambda_1$,~$\lambda_2$, and~$m$.
	\end{lemma}
	\begin{proof}
		We establish the first bound by verifying the conditions of Theorem 2.1 in~\cite{chernozhuokov2022improved}, cf.~in particular the display following that theorem. The second bound is proven analogously. First, note that for~$j=1,\hdots,d$
		\begin{align*}
			\phi_{\overline{\alpha}_j,\overline{\beta}_j}(X_{i,j})-\E\phi_{\overline{\alpha}_j,\overline{\beta}_j}(X_{i,j})
			\leq
			\overline{\beta}_j-\overline{\alpha}_j=Q_{1-c_{1,n}\eps_n}(X_{1,j})-Q_{c_{2,n}\eps_n}(X_{1,j}),	
		\end{align*}
		Furthermore, by Lemma~\ref{lem:quantile_mean} and~$\eps_n\geq \lambda_2\frac{\log(dn)}{n}$, the far right-hand side of the previous display is upper bounded by
		\begin{align*}
			\frac{\sigma_m}{(c_{1,n}\eps_n)^{1/m}}+\frac{\sigma_m}{(c_{2,n}\eps_n)^{1/m}}
			\leq
			C_1\sigma_m\del[3]{\frac{n}{\log(dn)}}^{1/m},
		\end{align*}
		where, by Lemma~\ref{lem:cControl}, $C_1$ is a constant depending only on~$\lambda_1$,~$\lambda_2$ and~$m$. Similarly, 
		\begin{align*}
			\phi_{\overline{\alpha}_j,\overline{\beta}_j}(X_{i,j})-\E\phi_{\overline{\alpha}_j,\overline{\beta}_j}(X_{i,j})
			&\geq
			\overline{\alpha}_j-\overline{\beta}_j
			=
			Q_{c_{2,n}\eps_n}(X_{1,j})-Q_{1-c_{1,n}\eps_n}(X_{1,j})\\
			&\geq
			-C_1\sigma_m\del[3]{\frac{n}{\log(dn)}}^{1/m},
		\end{align*}
		and we define
		\begin{equation*}
			B_n=2C_1\sigma_m\del[3]{\frac{n}{\log(dn)}}^{1/m}.
		\end{equation*}
		By construction,~$\envert[1]{\phi_{\overline{\alpha}_j,\overline{\beta}_j}(X_{i,j})-\E\phi_{\overline{\alpha}_j,\overline{\beta}_j}(X_{i,j})}\leq \frac{1}{2}B_n$ for all~$i=1,\hdots, n$ and~$j=1,\hdots,d$. Furthermore, by assumption,~$\E\sbr[1]{\phi_{\overline{\alpha}_j,\overline{\beta}_j}(X_{1,j})-\E\phi_{\overline{\alpha}_j,\overline{\beta}_j}(X_{1,j})}^2>\mathfrak{b}_1^2$, as well as
		\begin{align*}
			\E\sbr[1]{\phi_{\overline{\alpha}_j,\overline{\beta}_j}(X_{1,j})-\E\phi_{\overline{\alpha}_j,\overline{\beta}_j}(X_{1,j})}^4 &\leq B_n^2 \E\sbr[1]{\phi_{\overline{\alpha}_j,\overline{\beta}_j}(X_{1,j})-\E\phi_{\overline{\alpha}_j,\overline{\beta}_j}(X_{1,j})}^2 \\
			&\leq B_n^2 \E\sbr[1]{X_{1,j}-\E(X_{1,j})}^2 \leq B_n^2 \sigma_2^2 \leq B_n^2 b_2^2;
		\end{align*}
		a proof of the second inequality can be found in, e.g., Corollary 3 in~\cite{cstudd}. Therefore, by Theorem 2.1 and the display following it in~\cite{chernozhuokov2022improved}, there exists a constant~$C_2$ depending only on~$\mathfrak{b}_1$ and~$b_2$ such that
		\begin{align*}
			& \sup_{H\in\mc{H}}\envert[3]{\P\del[3]{\frac{1}{n^{1/2}}\sum_{i=1}^n\sbr[1]{\phi_{\overline{\alpha}_j,\overline{\beta}_j}(X_{i,j})-\E\phi_{\overline{\alpha}_j,\overline{\beta}_j}(X_{i,j})}\in H}-\P\del[1]{\overline{Z}\in H}}\\
			\leq & C_2\del[3]{\frac{B_n^2\log^5(dn)}{n}}^{\frac{1}{4}}
			=
			C_2\cdot(2C_1)^{1/2}\sigma_m^{1/2}\del[3]{\frac{\log^{5-\frac{2}{m}}(dn)}{n^{1-\frac{2}{m}}}}^{\frac{1}{4}},
		\end{align*}
		from which the conclusion of the lemma follows.
	\end{proof}

	\begin{lemma}\label{lem:I1I3BD}
		Let Assumption~\ref{ass:setting} be satisfied with~$m>2$. If~$\eps_n$ is chosen as in~\eqref{eq:epsfam} and satisfies~\eqref{eq:epscond},
		then for~$\overline I_{n,3}$ and~$\underline I_{n,3}$ as defined in \eqref{eq:I_13_max}, it holds that
		\begin{align*}
			\underline{I}_{n,3} \vee \overline{I}_{n,3}
			\leq 
			\sigma_m
			C\del[3]{\overline{\eta}_n^{1-\frac{1}{m}}+\sbr[2]{\frac{\log(dn)}{n}}^{1-\frac{1}{m}}},
		\end{align*}
		for a positive constant~$C$ depending only on~$\lambda_1$,~$\lambda_2$, and~$m$.
		
	\end{lemma}
	\begin{proof}
		This is a trivial consequence of Lemma~\ref{lem:meancontrol} and subadditivity of~$z\mapsto z^{1-\frac{1}{m}}$.	
	\end{proof}
	
	The following lemma collects some simple limits for later reference.
	\begin{lemma}\label{lem:covcompare}
		Let~$m>2$ and assume that~$\sqrt{n\log(d)}\overline{\eta}_n^{1-\frac{1}{m}}\to 0$ as well as~$\log(d)/n^{\frac{m-2}{5m-2}}\to 0$. Then
		\begin{equation}\label{eq:conv1}
			\log^2(d)\del[3]{\overline{\eta}_n^{1-\frac{2}{m}}+\sbr[2]{\frac{\log(dn)}{n}}^{1-\frac{2}{m}}}\to 0,
		\end{equation}	
		\begin{equation}\label{eq:conv2}
			\sbr[3]{\frac{\log^{5-\frac{2}{m}}(dn)}{n^{1-\frac{2}{m}}}}^{\frac{1}{4}}
			+
			\del[3]{\overline{\eta}_n^{1-\frac{1}{m}}+\sbr[2]{\frac{\log(dn)}{n}}^{1-\frac{1}{m}}}\sqrt{n\log(d)}\to 0,
		\end{equation}
		\begin{equation}\label{eq:conv3}
			\log^2(d)\sbr[3]{\overline{\eta}_n^{1-\frac{2}{m}}+\del[2]{\frac{\log(dn)}{n}}^{1-\frac{1}{(m/2)\wedge 2}}}\to 0,
		\end{equation}
		and
		\begin{equation}\label{eq:conv4}
			\sqrt{\log(d)\log(dn)}\del[3]{\overline{\eta}_n^{1-\frac{2}{m}}+\sbr[2]{\frac{\log(dn)}{n}}^{1-\frac{1}{(m/2)\wedge 2}}}
			\to 
			0.
		\end{equation}
	\end{lemma}
	\begin{proof}
		We begin by proving~\eqref{eq:conv1}, to which end we first establish that $\log^2(d)\overline{\eta}_n^{1-\frac{2}{m}}\to 0$. Since~$\log(d)/n^{\frac{m-2}{5m-2}}\to 0$, one has that~$\log^2(d)\overline{\eta}_n^{1-\frac{2}{m}}
		=
		o\del[1]{n^{\frac{2m-4}{5m-2}}\overline{\eta}_n^{1-\frac{2}{m}}}$ and
		\begin{equation*}
			n^{\frac{2m-4}{5m-2}}\overline{\eta}_n^{1-\frac{2}{m}}\to 0
			\iff 
			n^{\frac{2m^2-4m}{(5m-2)(m-2)}}\overline{\eta}_n\to 0.
		\end{equation*}
		The latter convergence is satisfied since~$\sqrt{n}\overline{\eta}^{1-\frac{1}{m}}\to 0$ by assumption, which is equivalent to~$n^{\frac{m}{2(m-1)}}\overline{\eta}_n\to 0$, and~$n^{\frac{2m^2-4m}{(5m-2)(m-2)}}\overline{\eta}_n\leq n^{\frac{m}{2(m-1)}}\overline{\eta}_n\to 0$ for~$m>2$.
		
		Next,
		\begin{equation*}
			\log^2(d)\sbr[2]{\frac{\log(d)}{n}}^{1-\frac{2}{m}}
			=
			\frac{\log(d)^{3-\frac{2}{m}}}{n^{1-\frac{2}{m}}}
			\to 0
			\iff 
			\frac{\log(d)}{n^{\frac{m-2}{3m-2}}}\to 0
			,
		\end{equation*}
		the latter convergence following from~$\log(d)/n^{\frac{m-2}{5m-2}}\to 0$ by assumption.
		
		Finally, 
		\begin{equation*}
			\log^2(d)\sbr[2]{\frac{\log(n)}{n}}^{1-\frac{2}{m}}\to 0
		\end{equation*}
		by a standard subsequence argument, considering separately the cases of subsequences along which i)~$n\leq d$ for which the convergence follows from the penultimate display and ii)~$d< n$ for which the convergence follows from~$m>2$. This establishes~\eqref{eq:conv1}.
		
		To prove~\eqref{eq:conv2}, note that
		\begin{equation*}
			\frac{\log^{5-\frac{2}{m}}(dn)}{n^{1-\frac{2}{m}}}\to 0 \iff \frac{\log(d)}{n^{\frac{m-2}{5m-2}}}\to 0,
		\end{equation*}
		the latter convergence being true by assumption. Furthermore,~$\sqrt{n\log(d)}\overline{\eta}_n^{1-\frac{1}{m}}\to 0$ by assumption and it remains to prove that
		\begin{equation*}
			\sqrt{n\log(d)}\del[2]{\frac{\log(dn)}{n}}^{1-\frac{1}{m}}\to 0.
		\end{equation*}
		To this end, note that
		\begin{equation*}
			\sqrt{n\log(d)}\del[2]{\frac{\log(d)}{n}}^{1-\frac{1}{m}}
			=
			\frac{\sbr[1]{\log(d)}^{\frac{3}{2}-\frac{1}{m}}}{n^{\frac{1}{2}-\frac{1}{m}}}\to 0
			\iff
			\frac{\log(d)}{n^{\frac{m-2}{3m-2}}}\to 0,	
		\end{equation*}
		where the latter convergence was already verified above. The convergence in the penultimate display now follows by considering separately subsequences along which~$n\leq d$ and~$d<n$, respectively (as in the concluding argument of the proof of~\eqref{eq:conv1}). 
		
		Next, to establish~\eqref{eq:conv3}, note that we have already shown that~$\log^2(d)\overline{\eta}_n^{1-\frac{2}{m}}\to 0$. Furthermore, in case~$m\leq 4$,
		\begin{equation*}
			\log^2(d)\del[2]{\frac{\log(d)}{n}}^{1-\frac{1}{(m/2)\wedge 2}}
			=
			\log^2(d)\del[2]{\frac{\log(d)}{n}}^{1-\frac{2}{m}}
			=
			\frac{\sbr[1]{\log(d)}^{3-\frac{2}{m}}}{n^{1-\frac{2}{m}}}
			\to 0
			\iff
			\frac{\log(d)}{n^{\frac{m-2}{3m-2}}}\to 0,
		\end{equation*}
		the latter convergence being true because it is assumed that even~$\log(d)/n^{\frac{m-2}{5m-2}}\to 0$. Conclude, once more, by a subsequence argument. The case of~$m>4$ is handled analogously.
		
		Finally, we establish~\eqref{eq:conv4} by showing that any subsequence possesses a further subsequence along which~\eqref{eq:conv4} is true. Thus, fix a subsequence. In case there exists a subsequence thereof along which~$n\leq d$, it suffices to show that along this subsequence
		\begin{equation*}
			\log(d)\del[3]{\overline{\eta}_n^{1-\frac{2}{m}}+\sbr[2]{\frac{\log(d)}{n}}^{1-\frac{1}{(m/2)\wedge 2}}}
			\to 
			0,	
		\end{equation*}
		which is implied by~\eqref{eq:conv3}. 
		
		In the remaining case (where there does not exist a further subsequence along which~$n\leq d$), we have~$d<n$ for~$n$ sufficiently large, so that it suffices to show that
		\begin{equation*}
			\log(n)\del[3]{\overline{\eta}_n^{1-\frac{2}{m}}+\sbr[2]{\frac{\log(n)}{n}}^{1-\frac{1}{(m/2)\wedge 2}}}
			\to 
			0.
		\end{equation*}
		Since~$m>2$ this, in turn, is true if
		\begin{equation*}
			\log(n)\overline{\eta}_n^{1-\frac{2}{m}}
			=
			\sqrt{n}\overline{\eta}_n^{1-\frac{1}{m}}\cdot \frac{\log(n)}{\sqrt{n}}\overline{\eta}_n^{-\frac{1}{m}}
			\to 0,	
		\end{equation*}
		which, since~$\sqrt{n}\overline{\eta}_n^{1-\frac{1}{m}}\to 0$ by assumption, is true in case~$\frac{\log(n)}{\sqrt{n}}\overline{\eta}_n^{-\frac{1}{m}}$ is bounded. In case~$\frac{\log(n)}{\sqrt{n}}\overline{\eta}_n^{-\frac{1}{m}}$ is unbounded,~$\frac{\log(n)}{\sqrt{n}}\overline{\eta}_n^{-\frac{1}{m}}\geq 1$ along a further subsequence. The latter implies~$\overline{\eta}_n\leq \frac{\log^m(n)}{n^{\frac{m}{2}}}$, such that
		\begin{equation*}
			\log(n)\overline{\eta}_n^{1-\frac{2}{m}}
			\leq
			\frac{[\log(n)]^{m-1}}{n^{\frac{m-2}{2}}}
			\to 0,
		\end{equation*}
		because~$m>2$.
	\end{proof}
	
	The following lemma shows that to establish Theorem~\ref{thm:HDGauss}, it suffices to prove it for the class of one-sided rectangles~$H=\bigtimes_{j=1}^d [-\infty,t_j]$ with~$t_j\in \overline{\R}=\R\cup\cbr[0]{-\infty,\infty}$ for~$j=1,\hdots,d$.
	\begin{lemma}\label{lem:onesidedtorectangle}
		Suppose Theorem~\ref{thm:HDGauss} holds for~$\rho_{n, W}$ throughout replaced by
		\begin{equation}
			\rho_{n, W}^{os} := \sup_{t\in\overline{\R}^d}\envert[1]{\P\del[1]{S_{n,W}\leq t}-\P\del[1]{Z\leq t}},
		\end{equation}
		and with~$\lambda_2$ in the definition of~$\eps_n$ in~\eqref{eq:epsfam} replaced by~$\lambda_{2,n}\in[\lambda_2/2,\lambda_2]$. Then, Theorem~\ref{thm:HDGauss} holds (with a different constant~$C$ depending only on~$b_1,b_2,\lambda_1,\lambda_2$, and~$m$).
	\end{lemma}
	
	\begin{proof}
		Fix~$a$ and~$b$ in~$\overline{\R}^d$ such that~$-\infty \leq a_j\leq b_j\leq \infty$ and let~$H=\bigtimes_{j=1}^d [a_j,b_j]\in\mc{H}$. Clearly,
		\begin{equation}\label{eq:twodim1}
			\P\del[1]{S_{n,W}\in H}
			=
			\P\del[1]{(S_{n,W}',-S_{n,W}')'\leq (b',-a')'}.
		\end{equation}
		Furthermore, by~\eqref{eqn:winsmeandef}
		\begin{equation*}
			-S_{n,W,j}
			=
			n^{-1/2}\sum_{i=1}^n\del[1]{\phi_{-\hat\beta_j,-\hat\alpha_j}(-\tilde{X}_{i,j})-(-\mu_j)}
			,\qquad j=1,\hdots,d.
		\end{equation*}
		Thus, the vectors~$(X_i',-X_i')'$ and~$(\tilde{X}_i',-\tilde{X}_i')'$ in~$\R^{2d}$ satisfy Assumption~\ref{ass:setting} with~$d$ there replaced by~$2d$ (but with the same~$m$, $b_1$,~$b_2$ and~$\overline{\eta}_n$). Note that the covariance matrix of~$(X_1',-X_1')'$ is
		\begin{equation*}
			\Xi= 
			\begin{pmatrix}
				\Sigma & -\Sigma\\
				-\Sigma & \Sigma
			\end{pmatrix},
		\end{equation*}
		and let~$Z_2\sim~\mathsf{N}_{2d}(0,\Xi)=:\nu'$. Note also that for fixed~$\lambda_1\in(1,\infty)$ and~$\lambda_2\in(0,\infty)$,~$\eps_n$ defined in~\eqref{eq:epsfam} can be written as
		\begin{equation*}
			\eps_n
			=			
			\lambda_1\cdot \overline{\eta}_n +\lambda_{2,n}\cdot \frac{\log(2dn)}{n},\qquad\text{where}\qquad\lambda_{2,n}=\lambda_2\cdot\frac{\log(dn)}{\log(2dn)}\in[\lambda_2/2,\lambda_2],
		\end{equation*}
		the inclusion following from~$d\geq 2$ and~$n>3$. Hence, the version of Theorem~\ref{thm:HDGauss} that replaces~$\rho_{n, W}$ by~$\rho_{n, W}^{os}$ and~$\lambda_2$ by~$\lambda_{2,n}\in[\lambda_2/2,\lambda_2]$ (and which is assumed to hold true) applies to~$(S_{n,W}',-S_{n,W}')'$ and yields, for~$C^*$ a constant depending only on~$b_1,b_2,\lambda_1,\lambda_2$, and~$m$, the upper bound 
		\begin{align*}
			C^*\del[4]{\sbr[3]{\frac{\log^{5-\frac{2}{m}}(dn)}{n^{1-\frac{2}{m}}}}^{\frac{1}{4}}
				+
				\sbr[3]{\overline{\eta}_n^{1-\frac{1}{m}}+\sbr[2]{\frac{\log(dn)}{n}}^{1-\frac{1}{m}}}\sqrt{n\log(d)}} \\
			+
			C^*\del[4]{\log^2(d)\sbr[3]{\overline{\eta}_n^{1-\frac{2}{m}}+\sbr[2]{\frac{\log(dn)}{n}}^{1-\frac{2}{m}}}}^{1/2}
		\end{align*}
		on
		\begin{equation}\label{eq:twodim2}
			\envert[1]{\P\del[1]{(S_{n,W},-S_{n,W})\leq (b',-a')'}-\P\del[1]{Z_2\leq(b',-a')'}} ,
		\end{equation}
		Furthermore, for~$\nu=\mathsf{N}_d(0,\Sigma)$,~$\nu'$ is the image measure of~$\nu$ under the mapping~$\R^d\ni z\mapsto (z,-z)\in\R^{2d}$, that is~$\nu'=\nu\circ \del[1]{z\mapsto (z,-z)}^{-1}$. Thus,
		\begin{align}\label{eq:twodim3}
			\P\del[1]{Z_2\leq(b',-a')'}
			&=
			\nu'\del[1]{[-\infty,b]\times [-\infty,-a]}\notag\\
			&=
			\nu\del[1]{z\in\R^d: (z,-z)\in [-\infty,b]\times [-\infty,-a]}\notag\\
			&=
			\P\del[1]{(Z,-Z)\leq (b,-a)}\notag\\
			&=
			\P\del[1]{Z\in H}.
		\end{align}
		Combining~\eqref{eq:twodim1} and~\eqref{eq:twodim3} with the upper bound in~\eqref{eq:twodim2} obtained above (which does not depend on~$a$ or~$b$) delivers the claim. 
	\end{proof}
	
	\section{Proof of Theorem~\ref{thm:HDGauss}}
	By Lemma~\ref{lem:onesidedtorectangle} it suffices to prove~\eqref{eq:Gaussapprox} for hyperrectangles~$H$ of the form~$H=\bigtimes_{j=1}^d[-\infty,t_j]$ with~$ t_j\in \overline{\R}$ and~$\eps_n$ as defined in~\eqref{eq:epsfam}, but with~$\lambda_2$ replaced by~$\lambda_{2,n}\in[\lambda_2/2,\lambda_2]$. Throughout the proof, we consider the case where~$\lambda_{2,n}=\lambda_2$ for notational convenience. Inspection of the argument shows that it remains valid for~$\lambda_{2,n}\in[\lambda_2/2,\lambda_2]$. Furthermore, recalling the definitions of~$\mathfrak{T}_n^{(1)}$ and~$\mathfrak{T}_n^{2}$ in~\eqref{eq:epscond} and \eqref{eq:truncvar}, respectively, and using that~$\log(dn)/n\leq \lambda_2^{-1}\eps_n$, one has with~$C$ being the constant from~\eqref{eq:truncvar} that
	\begin{equation}\label{eq:auxgood}
		\mathfrak{T}_n^{(1)}+\mathfrak{T}_n^{(2)}
		\leq
		2\eps_n+\lambda_2^{-1}\eps_n+\sqrt{\lambda_2^{-2}\eps_n^2+4\lambda_2^{-1}\eps_n^2}+Cb_2^2\eps_n^{1-\frac{2}{m}}
		\leq
		C_1\eps_n^{1-\frac{2}{m}},
	\end{equation}
	for~$C_1$ a constant depending only on~$b_2$,~$\lambda_2$~and~$m$, and where we used that~$\eps_n\leq \eps_n^{1-\frac{2}{m}}$ as~$\eps_n\in(0,0.5)$. In the following we first suppose that
	\begin{equation}\label{eq:CaseGood}
		C_1\eps_n^{1-\frac{2}{m}}
		<
		\frac{(b_1^2\wedge 1)}{2},
	\end{equation}
	which, by the penultimate display, implies that~$\mathfrak{T}_n^{(1)}<1/2$ such that~$\eqref{eq:epscond}$ is satisfied. By Lemma~\ref{lem:noisecontrol}, there exists a constant~$C_2$ depending only on~$\lambda_1$,~$\lambda_2$, and~$m$ such that
	\begin{equation}\label{eq:I1}
		\underline{I}_{n,1}\vee \overline{I}_{n,1}
		\leq
		C_2 \sigma_m \overline{\eta}_n^{1-\frac{1}{m}}
		:=I_{n,1}
	\end{equation}
	with~$I_{n,1}$ being non-random. Thus, from~\eqref{eq:LBUPHD} it holds with probability at least~$1-\frac{4}{n}$ that
	\begin{equation*}
		\underline{Y}_{n}:=\sqrt{n}\underline{I}_{n,2}-\sqrt{n}(I_{n,1}+\underline{I}_{n,3})
		\leq 
		Y_n
		\leq 
		\sqrt{n}\overline{I}_{n,2}+\sqrt{n}(I_{n,1}+\overline{I}_{n,3})=:\overline{Y}_n.
	\end{equation*}
	Next, it holds for all~$t=(t_1,\hdots,t_d)$ with~$t_j\in \overline{\R}$ that
	\begin{equation*}
		\P\del[1]{Y_n\leq t}
		\geq
		\P\del[1]{\overline{Y}_n\leq t, Y_n\leq \overline{Y}_n}
		\geq
		\P\del[1]{\overline{Y}_n\leq t}+\P\del[1]{Y_n\leq \overline{Y}_n}-1
		\geq
		\P\del[1]{\overline{Y}_n\leq t}-\frac{4}{n}.
	\end{equation*}
	Similarly,
	\begin{equation*}
		\P\del[1]{Y_n\leq t}
		\leq
		\P\del[1]{\underline{Y}_n\leq t}+\P\del[1]{\underline{Y}_n>Y_n}
		\leq	
		\P\del[1]{\underline{Y}_n\leq t}+\frac{4}{n}.
	\end{equation*}
	Thus,
	\begin{align}\label{eq:1overn}
		\sup_{t\in\overline{\R}^d}\envert[2]{\P\del[1]{Y_n\leq t}-\P\del[1]{Z\leq t}}	
		&\leq
		\sup_{t\in\overline{\R}^d}\envert[2]{\P\del[1]{\underline{Y}_n\leq t}-\P\del[1]{Z\leq t}}\notag\\
		&+\sup_{t\in\overline{\R}^d}\envert[2]{\P\del[1]{\overline{Y}_n\leq t}-\P\del[1]{Z\leq t}}
		+\frac{4}{n}.
	\end{align}
	%
	We proceed by bounding the second summand on the right-hand side of the previous display (the argument for the first summand is analogous and hence skipped). To this end, note that~\eqref{eq:auxgood} and~\eqref{eq:CaseGood} imply that~$\mathfrak{T}_n^{(2)}\leq b_1^2/2$ which then together with Lemma~\ref{lem:trunccov} implies that~$\min_{j=1,\hdots, d}\overline{\Sigma}_{\eps_n,j,j}\geq b_1^2/2$. Next, by the definition of~$\overline{Y}_n$
	\begin{equation*}
		\P\del[1]{\overline{Y}_n\leq t}
		=
		\P\del[1]{\sqrt{n} \overline{I}_{n,2}\leq t-\sqrt{n}\del[0]{I_{n,1}+\overline{I}_{n,3}}},\qquad t\in \overline{\R}^d,
	\end{equation*}
	whence, noting also that~$I_{n,1}$ and~$\overline{I}_{n,3}$ are non-random, Lemma~\ref{lem:HDGaussBD} applied with~$\mathfrak{b}_1^2=b_1^2/2$ implies the existence of a constant~$C_3$ depending only on~$b_1,b_2,\lambda_1,\lambda_2$, and~$m$ such that
	\begin{equation*}
		\sup_{t\in\overline{\R}^d}\envert[3]{\P\del[1]{\overline{Y}_n\leq t}-\P\del[3]{\overline{Z}\leq t-\sqrt{n}\del[0]{I_{n,1}+\overline{I}_{n,3}}}}
		\leq 
		C_3\del[3]{\frac{\log^{5-\frac{2}{m}}(dn)}{n^{1-\frac{2}{m}}}}^{\frac{1}{4}},	
	\end{equation*}
	where~$\overline{Z}\sim \mathsf{N}_d(0,\overline{\Sigma}_{\eps_n})$.
	Furthermore, by the Gaussian anti-concentration inequality as stated in Theorem 1 of~\cite{chernozhukov2017detailed}, (cf.~also Lemma A.1 in~\cite{chernozhukov2017central}),
	\begin{align*}
		0
		&\geq
		\P\del[1]{\overline{Z}\leq t-\sqrt{n}\del[0]{I_{n,1}+\overline{I}_{n,3}}}-\P\del[1]{\overline{Z}\leq t}\\
		&\geq
		-\frac{\sqrt{n}\del[0]{I_{n,1}+\overline{I}_{n,3}}}{\sqrt{\min_{j=1,\hdots, d}\overline{\Sigma}_{\eps_n,j,j}}}\del[1]{\sqrt{2\log(d)}+2}\\
		&\geq
		-C_4\del[3]{\overline{\eta}_n^{1-\frac{1}{m}}+\sbr[2]{\frac{\log(dn)}{n}}^{1-\frac{1}{m}}}\sqrt{n\log(d)},
	\end{align*}
	the final inequality following from Lemma~\ref{lem:I1I3BD},~\eqref{eq:I1},~$\min_{j=1,\hdots, d}\overline{\Sigma}_{\eps_n,j,j}\geq b_1^2/2$, and~$C_4$ being a constant depending on~$b_1$, $b_2$,~$\lambda_1,\lambda_2$, and~$m$ only. Thus, combining the previous two displays, there exists a constant~$C_5$ depending on~$b_1,b_2,\lambda_1,\lambda_2$, and~$m$ only such that when~\eqref{eq:CaseGood} is satisfied one has
	\begin{align*}
		\sup_{t\in\overline{\R}^d}\envert[2]{\P\del[1]{\overline{Y}_n\leq t}-\P\del[1]{\overline{Z}\leq t}}
		\leq
		C_5\del[4]{\sbr[3]{\frac{\log^{5-\frac{2}{m}}(dn)}{n^{1-\frac{2}{m}}}}^{\frac{1}{4}}
			+
			\sbr[3]{\overline{\eta}_n^{1-\frac{1}{m}}+\sbr[2]{\frac{\log(dn)}{n}}^{1-\frac{1}{m}}}\sqrt{n\log(d)}}.
	\end{align*}
	Finally, by the Gaussian-to-Gaussian comparison inequality as stated in Proposition 2.1 of~\cite{chernozhuokov2022improved} and Lemma~\ref{lem:trunccov}, there exists a constant~$C_6$, depending on~$b_1$,~$b_2$,~$\lambda_1$,~$\lambda_2$, and~$m$ only ($C_6$ changes value in the second inequality below), such that
	\begin{align*}
		\sup_{t\in\overline{\R}^d}\envert[2]{\P\del[1]{\overline{Z}\leq t}-\P\del[1]{Z\leq t}}	
		&\leq
		C_6\del[2]{\log^2(d)\eps_n^{1-\frac{2}{m}}}^{1/2}\\
		&\leq
		C_6\del[4]{\log^2(d)\sbr[3]{\overline{\eta}_n^{1-\frac{2}{m}}+\sbr[2]{\frac{\log(dn)}{n}}^{1-\frac{2}{m}}}}^{1/2}.
	\end{align*}
	The previous two displays yield that there exists a constant~$C_7$ depending only on~$b_1$,~$b_2$,~$\lambda_1,\lambda_2$, and~$m$ such that when~\eqref{eq:CaseGood} is satisfied it holds that
	\begin{align*}
		\sup_{t\in\overline{\R}^d}\envert[2]{\P\del[1]{\overline{Y}_n\leq t}-\P\del[1]{Z\leq t}}
		\leq
		C_7\del[4]{\log^2(d)\sbr[3]{\overline{\eta}_n^{1-\frac{2}{m}}+\sbr[2]{\frac{\log(dn)}{n}}^{1-\frac{2}{m}}}}^{1/2}\\
		+\ C_7\del[4]{\sbr[3]{\frac{\log^{5-\frac{2}{m}}(dn)}{n^{1-\frac{2}{m}}}}^{\frac{1}{4}}
			+
			\sbr[3]{\overline{\eta}_n^{1-\frac{1}{m}}+\sbr[2]{\frac{\log(dn)}{n}}^{1-\frac{1}{m}}}\sqrt{n\log(d)}}.
	\end{align*}
	This implies the claimed bound on~$\rho_{n,W}$ in~\eqref{eq:Gaussapprox} by using that~$4/n$ in~\eqref{eq:1overn} is dominated by, e.g., $\sbr[1]{\frac{\log(dn)}{n}}^{1-\frac{1}{m}}\sqrt{n\log(d)}$ (if necessary adjust~$C_7$). 
	
	If, on the other hand,~\eqref{eq:CaseGood} is not satisfied, then 
	\begin{align*}
		C_1\eps_n^{1-\frac{2}{m}}
		\geq
		\frac{(b_1^2\wedge 1)}{2}	\quad\iff \quad 
		\eps_n
		\geq
		\sbr[3]{\frac{(b_1^2\wedge 1)}{2C_1}}^{\frac{m}{m-2}}=:K(b_1,b_2,\lambda_2,m),	
	\end{align*}
	where~$K=K(b_1,b_2,\lambda_2,,m)$ does not depend on~$n$. Thus, by the definition of~$\eps_n$, 
	\begin{align*}
		\overline{\eta}_n\geq\frac{K}{2\lambda_1}\qquad\text{or}\qquad \frac{\log(dn)}{n}\geq\frac{K}{2\lambda_2}, 
	\end{align*}
	In either case, since~$\rho_{n,W}\leq 1$ the bound in~\eqref{eq:Gaussapprox} remains valid by adjusting the constant~$C$ there, if necessary.
	
	Finally,~$\rho_{n,W}\to 0$ by~\eqref{eq:conv1} and~\eqref{eq:conv2} of Lemma~\ref{lem:covcompare}.
	
	\section{Proofs for Section~\ref{sec:Bootstrap}}
	We start with the proof of Theorem~\ref{thm:covestimGRAM}. Before we give the proof, we first define~$c_{1,n}'$ and~$c_{2,n}'$ that appear in the proof: Recall the definition of~$\eps_n'$ in~\eqref{eq:epsprime}. For~$\lambda_1'\in(1,\infty)$, define~$A_+'=1-\lambda_1'^{-1}\mathds{1}\cbr[0]{\overline{\eta}_n>0}$,~$A_-'=1+\lambda_1'^{-1}\mathds{1}\cbr[0]{\overline{\eta}_n>0}$, and 
	\begin{equation}\label{eq:c1'}
		c_{1,n}':=-A_+'W_0\del[1]{-e^{-(\frac{\log(d^2n)}{\eps_n' n}+A_+')/A_+'}}\in(0,A_+')	\end{equation}
	as well as
	\begin{equation}\label{eq:c2'}
		c_{2,n}':=-A_-'W_{-1}\del[1]{-e^{-(\frac{\log(d^2n)}{\eps_n' n}+A_-')/A_-'}}\in(A_-',\infty);
	\end{equation}
	where we recall that~$W_0$ and~$W_{-1}$ denote the principal and lower branch, respectively, of Lambert's~$W$ function (cf., e.g.,~\cite{Corless1996}).
	\begin{proof}[Proof of Theorem~\ref{thm:covestimGRAM}]
		Fix~$1\leq j,k\leq d$ and note that
		\begin{align*}
			\tilde{\Sigma}_{n,j,k}
			&=
			\frac{1}{n}\sum_{i=1}^n\sbr[1]{\phi_{\hat{a}_j,\hat{b}_j}(\tilde{X}_{i,j})-\mu_j-(\tilde{\mu}_{n,j}-\mu_j)}\sbr[1]{\phi_{\hat{a}_k,\hat{b}_k}(\tilde{X}_{i,k})-\mu_k-(\tilde{\mu}_{n,k}-\mu_k)}\\
			&=
			\frac{1}{n}\sum_{i=1}^n\sbr[1]{\phi_{\hat{a}_j,\hat{b}_j}(\tilde{X}_{i,j})-\mu_j}\sbr[1]{\phi_{\hat{a}_k,\hat{b}_k}(\tilde{X}_{i,k})-\mu_k}-(\tilde{\mu}_{n,j}-\mu_j)(\tilde{\mu}_{n,k}-\mu_k)\\
			&=
			\frac{1}{n}\sum_{i=1}^n\phi_{\check{a}_j,\check{b}_j}(\tilde{Y}_{i,j})\phi_{\check{a}_k,\check{b}_k}(\tilde{Y}_{i,k})-(\tilde{\mu}_{n,j}-\mu_j)(\tilde{\mu}_{n,k}-\mu_k),
		\end{align*}
		where~$Y_{i,j}=X_{i,j}-\mu_j$,~$\tilde{Y}_{i,j}=\tilde{X}_{i,j}-\mu_j$,~$\check{a}_j=\tilde{Y}_{\lceil \eps_n'n\rceil,j}^*=\tilde{X}_{\lceil \eps_n'n\rceil,j}^*-\mu_j=\hat{a}_j-\mu_j$, and~$\check{b}_j=\tilde{Y}_{\lfloor (1-\eps_n')n\rfloor+1,j}^*=\tilde{X}_{\lfloor (1-\eps_n')n\rfloor+1,j}^*-\mu_j=\hat{b}_j-\mu_j$.
		
		By Theorem 2.1 in~\cite{Wins1} (the assumptions of which are implied by Assumption~\ref{ass:setting} of the present paper) applied with~$\delta=6/(d^2n)$,~$\lambda_i=\lambda_i'$ for~$i=1,2$, there exists a constant~$C_1$ depending only on~$b_2,\lambda_1',\lambda_2'$, and~$m$, such that with probability at least~$1-12/(d^2n)$ 
		\begin{align*}
			|\tilde{\mu}_{n,j}-\mu_j||\tilde{\mu}_{n,k}-\mu_k|
			&\leq 
			C_1\del[3]{\overline{\eta}_n^{1-\frac{1}{m}}+\sbr[2]{\frac{\log(dn)}{n}}^{1-\frac{1}{m\wedge 2}}}^2\\
			&\leq
			C_1\del[3]{\overline{\eta}_n^{1-\frac{2}{m}}+\sbr[2]{\frac{\log(dn)}{n}}^{1-\frac{1}{(m/2)\wedge 2}}},
		\end{align*}
		the second estimate using that~$\overline{\eta}_n$ and~$\frac{\log(dn)}{n}$ are both strictly smaller than 1 by~\eqref{eq:epscond'} (and potentially adjusting~$C_1$).
		
		Next, write~$U_{i,j,k}=Y_{i,j}Y_{i,k}$,~$\tilde{U}_{i,j,k}=\tilde{Y}_{i,j}\tilde{Y}_{i,k}$,~$\check{\alpha}_{j,k}=\tilde{U}_{\lceil \eps_n'n\rceil,j,k}^*$, and~$\check{\beta}_{j,k}=\tilde{U}_{\lfloor (1-\eps_n')n\rfloor+1,j,k}^*$. Note that~$\E U_{1,j,k}=\Sigma_{j,k}$ and~$\E|U_{1,j,k}|^{m/2}\leq \sigma_m^m\leq b_2^m$. Therefore, by Theorem 2.1 in~\cite{Wins1} applied with~$\delta=6/(d^2n)$,~$m$ \emph{there} being~$m/2$,~$\lambda_i=\lambda_i'$ for~$i=1,2$, there exists a constant~$C_2$ depending only on~$b_2,\lambda_1',\lambda_2'$, and~$m$ such that with probability at least~$1-6/(d^2n)$ 
		\begin{align*}
			\envert[2]{\frac{1}{n}\sum_{i=1}^n\phi_{\check{\alpha}_{j,k},\check{\beta}_{j,k}}(\tilde{U}_{i,j,k})-\Sigma_{j,k}}
			\leq C_2\del[3]{\overline{\eta}_n^{1-\frac{2}{m}}+\sbr[2]{\frac{\log(dn)}{n}}^{1-\frac{1}{(m/2)\wedge 2}}}.
		\end{align*}
		Thus, defining
		\begin{align*}
			\mathfrak{S}_{n,j,k}
			=:
			\frac{1}{n}\sum_{i=1}^n\phi_{\check{a}_j,\check{b}_j}(\tilde{Y}_{i,j})\phi_{\check{a}_k,\check{b}_k}(\tilde{Y}_{i,k})-\frac{1}{n}\sum_{i=1}^n\phi_{\check{\alpha}_{j,k},\check{\beta}_{j,k}}(\tilde{U}_{i,j,k}),
		\end{align*}
		in order to prove the theorem, it suffices (by the union bound over~$d^2$) to show that there exists a constant~$C_3$ depending only on~$b_2,\lambda_1',\lambda_2'$, and~$m$ such that with probability at least~$1-6/(d^2n)$ 
		\begin{align}\label{eq:Toshow}
			\envert[1]{\mathfrak{S}_{n,j,k}}
			\leq C_3\del[3]{\overline{\eta}_n^{1-\frac{2}{m}}+\sbr[2]{\frac{\log(dn)}{n}}^{1-\frac{1}{(m/2)\wedge 2}}}	
		\end{align}
		To establish~\eqref{eq:Toshow}, note that with
		\begin{align*}
			A_{n,j,k}:=\cbr[2]{i\in\cbr[1]{1,\hdots,n}: \check{\alpha}_{j,k}\leq \tilde{U}_{i,j,k}\leq \check{\beta}_{j,k} \text{ and }\check{a}_j\leq \tilde{Y}_{i,j}\leq\check{b}_j\text{ and }\check{a}_k\leq \tilde{Y}_{i,k}\leq\check{b}_k},
		\end{align*}
		one has that
		\begin{align*}
			\phi_{\check{\alpha}_{j,k},\check{\beta}_{j,k}}\del[0]{\tilde{U}_{i,j,k}}-\phi_{\check{a}_j,\check{b}_j}(\tilde{Y}_{i,j})\phi_{\check{a}_k,\check{b}_k}(\tilde{Y}_{i,k})
			=
			\tilde{U}_{i,j,k}-\tilde{Y}_{i,j}\tilde{Y}_{i,k}
			=0\qquad\text{for }i\in A_{n,j,k}. 
		\end{align*}
		Hence,
		\begin{align*}
			\mathfrak{S}_{n,j,k}
			&=
			\frac{1}{n}\sum_{i\in A_{n,j,k}^c}\sbr[1]{\phi_{\check{a}_j,\check{b}_j}(\tilde{Y}_{i,j})\phi_{\check{a}_k,\check{b}_k}(\tilde{Y}_{i,k})-\phi_{\check{\alpha}_{j,k},\check{\beta}_{j,k}}\del[0]{\tilde{U}_{i,j,k}}},
		\end{align*}	
		and for every~$i\in A_{n,j,k}^c$
		\begin{align}\label{eq:tailterms}
			\envert[2]{\phi_{\check{\alpha}_{j,k},\check{\beta}_{j,k}}\del[0]{\tilde{U}_{i,j,k}}-\phi_{\check{a}_j,\check{b}_j}(\tilde{Y}_{i,j})\phi_{\check{a}_k,\check{b}_k}(\tilde{Y}_{i,k})}
			\leq \del[1]{|\check{\alpha}_{j,k}| \vee |\check{\beta}_{j,k}|}	+\del[1]{|\check{a}_j| \vee |\check{b}_j|}\del[1]{|\check{a}_k| \vee |\check{b}_k|}.
		\end{align}
		To bound the right-hand side of~\eqref{eq:tailterms}, note that by Lemma~\ref{lem:quantiles2}
		\begin{align*}
			Q_{c_{1,n}'\eps_n'}(U_{1,j,k})
			\leq
			\check{\alpha}_{j,k}
			\leq 
			\check{\beta}_{j,k}
			\leq
			Q_{1-c_{1,n}'\eps_n'}(U_{1,j,k})
		\end{align*}
		with probability at least~$1-\frac{2}{d^2n}$, and where the second inequality used that~$\eps_n'\in(0,1/2)$ by~\eqref{eq:epscond'}. Therefore, with at least this probability,
		\begin{align*}
			|\check{\alpha}_{j,k}| \vee |\check{\beta}_{j,k}|
			\leq
			|Q_{c_{1,n}'\eps_n'}(U_{1,j,k})|\vee |Q_{1-c_{1,n}'\eps_n'}(U_{1,j,k})|	
		\end{align*}
		By the same argument, it holds with probability at least~$1-\frac{4}{d^2n}$ that
		\begin{align*}
			\del[1]{|\check{a}_j| \vee |\check{b}_j|}\del[1]{|\check{a}_k| \vee |\check{b}_k|}
			\leq &
			\del[1]{|Q_{c_{1,n}'\eps_n'}(Y_{1,j})|\vee  |Q_{1-c_{1,n}'\eps_n'}(Y_{1,j})|}\\
			\cdot & \del[1]{|Q_{c_{1,n}'\eps_n'}(Y_{1,k})|\vee  |Q_{1-c_{1,n}'\eps_n'}(Y_{1,k})|}.
		\end{align*}
		Thus, with probability at least~$1-\frac{6}{d^{2}n}$ the right-hand side of~\eqref{eq:tailterms} is bounded from above by
		\begin{align*}
			|Q_{c_{1,n}'\eps_n'}(U_{1,j,k})|\vee |Q_{1-c_{1,n}'\eps_n'}(U_{1,j,k})|
			+&\del[1]{|Q_{c_{1,n}'\eps_n'}(Y_{1,j})|\vee  |Q_{1-c_{1,n}'\eps_n'}(Y_{1,j})|}\\
			\cdot & \del[1]{|Q_{c_{1,n}'\eps_n'}(Y_{1,k})|\vee  |Q_{1-c_{1,n}'\eps_n'}(Y_{1,k})|}.
		\end{align*}
		Next, recall that~$\E U_{1,j,k}=\Sigma_{j,k}$,~$\E\envert[0]{U_{1,j,k}-\E U_{1,j,k}}^{m/2}
		\leq 2^{m/2}\E|U_{1,j,k}|^{m/2}\leq 
		2^{m/2}\sigma_m^m$, $\E Y_{1,j}=0$, and~$\E |Y_{1,j}|^m\leq \sigma_m^m$ for all~$1\leq j, k\leq d$. Hence, Lemma~\ref{lem:quantile_mean} implies that the previous display is bounded from above by 
		\begin{align*}
			\del[3]{|\Sigma_{j,k}|+\frac{2\sigma_m^2}{(c_{1,n}'\eps_n')^{2/m}}}+\frac{\sigma_m^2}{(c_{1,n}'\eps_n')^{2/m}}
			\leq
			\frac{4\sigma_m^2}{(c_{1,n}'\eps_n')^{2/m}},
		\end{align*}
		the inequality following from~$|\Sigma_{j,k}|\leq \sigma_2^2\leq \sigma_m^2$. Thus, since~$|A_{n,j,k}^c|\leq 6\eps_n' n$, with probability at least~$1-\frac{6}{d^2n}$ it holds that
		\begin{align*}
			\envert[1]{\mathfrak{S}_{n,j,k}}
			\leq
			6\eps_n'\cdot \frac{4\sigma_m^2}{(c_{1,n}'\eps_n')^{2/m}}
			\leq 
			\frac{24\sigma_m^2}{(c_{1,n}')^{2/m}}\cdot{\eps_n'}^{1-\frac{2}{m}}
			\leq 
			C_3\del[2]{\overline{\eta}_n^{1-\frac{2}{m}}+\sbr[2]{\frac{\log(d^2n)}{n}}^{1-\frac{2}{m}}},
		\end{align*}
		where we inserted~$\eps_n'$ from~\eqref{eq:epsprime}, used subadditivity of~$z\mapsto z^{1-\frac{2}{m}}$, the lower bound on~$c_{1,n}'$ from Lemma~\ref{lem:cControl'}, and~$C_3$ is the desired constant depending only on~$b_2,\lambda_1',\lambda_2'$, and~$m$. Hence,~\eqref{eq:Toshow} follows because~$1-\frac{2}{m}\geq 1-\frac{1}{(m/2)\wedge 2}$ and~$\frac{\log(d^2n)}{n}<1$ by~\eqref{eq:epscond'}.

	\end{proof}
	
	\begin{proof}[Proof of Theorem~\ref{thm:HDBootstrap}]
		Note first that writing~$\mathfrak{T}_n^{(3)}$ for the left-hand side of~\eqref{eq:epscond'}, it follows from~$\frac{\log(d^2n)}{n}\leq {\lambda_2'}^{-1}{\eps_n'}$ that
		\begin{align*}
			\mathfrak{T}_n^{(3)}
			\leq
			2\eps_n'+{\lambda_2'}^{-1}{\eps_n'}+\sqrt{{\lambda_2'}^{-2}{\eps_n'}^2+4{\lambda_2'}^{-1}{\eps_n'}^2}
			\leq
			C_1\eps_n',	
		\end{align*}
		for a constant~$C_1$ depending only on~$\lambda_2'$. 
		
		Suppose first that~$C_1\eps_n'<1$ such that~\eqref{eq:epscond'} is satisfied by the previous display. By the triangle inequality~$\tilde{\rho}_{n,W}$ is bounded from above by the sum of
		\begin{align*}
			\rho_{n,W}=\sup_{H\in\mc{H}}\envert[2]{\P\del[1]{S_{n,W}\in H}-\P\del[1]{Z\in H}},	
		\end{align*}
		where we recall that~$Z\sim\mathsf{N}_d(0,\Sigma)$, and
		\begin{align*}
			B_n:=\sup_{H\in\mc{H}}\envert[2]{\P\del[1]{Z\in H}-\P\del[1]{\tilde{Z}\in H\mid\tilde{X}_1,\hdots,\tilde{X}_n}}.
		\end{align*}
		First, $\rho_{n,W}\leq \mathfrak{A}_n$, where~$\mathfrak{A}_n$ is the upper bound on~$\rho_{n,W}$ in Theorem~\ref{thm:HDGauss}. Next, by the Gaussian-to-Gaussian comparison inequality as stated in Proposition 2.1 of~\cite{chernozhuokov2022improved} (cf.~also Proposition 2 in~\cite{chernozhukov2023high}),
		\begin{align*}
			B_n
			\leq 
			C\del[2]{\log^2(d)\max_{1\leq j,k\leq d}\envert[1]{\tilde{\Sigma}_{n,j,k}-\Sigma_{j,k}}}^{1/2}
			\leq
			C\del[4]{\log^2(d)\sbr[3]{\overline{\eta}_n^{1-\frac{2}{m}}+\del[2]{\frac{\log(dn)}{n}}^{1-\frac{1}{(m/2)\wedge 2}}}}^{1/2}	,
		\end{align*}
		the last inequality holding with probability at least~$1-\frac{24}{n}$ by Theorem~\ref{thm:covestimGRAM} and~$C$ being a constant depending on~$b_1,b_2,\lambda_1',\lambda_2'$, and~$m$ only (recall that~\eqref{eq:epscond'} is satisfied). 
		
		If, on the other hand,~$C_1\eps_n'\geq 1$ then by the definition of~$\eps_n'$ in~\eqref{eq:epsprime}%
		\begin{align*}
			\overline{\eta}_n\geq \frac{1}{2C_1\lambda_1'}\qquad\text{or}\qquad \frac{\log(dn)}{n}
			=
			\frac{\log([dn]^2)/2}{n}
			\geq
			\frac{\log(d^2n)}{2n}
			\geq
			\frac{1}{4C_1\lambda_2'}.
		\end{align*} 
		In either case, the upper bound on~$\tilde{\rho}_{n,W}$ remains valid by suitably adjusting~$C$ (if needed).
		
		That~$\tilde{\rho}_{n,W}\to 0$ in probability follows from~\eqref{eq:conv1}--\eqref{eq:conv3}	of Lemma~\ref{lem:covcompare}.
	\end{proof}
	
	\section{Proofs for Section~\ref{sec:Normalized}}
	
	The following ``two-sided'' Gaussian anti-concentration inequality is a simple consequence of the ``one-sided'' one stated in Theorem 1 of~\cite{chernozhukov2017detailed}, (cf.~also Lemma A.1 in~\cite{chernozhukov2017central}), but as we could not pinpoint it in the literature we state it here for completeness. For all~$u,v$ in~$\overline{\R}^d$, we define the set~$[u,v]=\cbr[1]{x\in\R^d: u_j\leq x_j\leq v_j\text{ for all }j=1,\hdots,d}$, which may be empty. Recall the notational conventions introduced after~\eqref{eq:LBUPHD}.
	\begin{lemma}\label{lem:anticonc}
		Let~$Z$ in~$\R^d$ with~$d\geq 1$ be such that~$Z\sim~\mathsf{N}_d(0,\Sigma)$ with~$\Sigma_{j,j} \geq \underline{\sigma}^2$ for all~$j=1,\hdots, d$ and some~$\underline{\sigma}^2>0$. Then, for all real numbers~$\delta_1$ and~$\delta_2$ and all~$a$ and~$b$ in~$\overline{\R}^d$, it holds that
		\begin{equation}\label{eq:twogauss1}
			\P\del[1]{Z\in[a+\delta_1,b+\delta_2]}
			\leq
			\P\del[1]{Z \in[a,b]}+\frac{\overline{\delta}}{\underline{\sigma}}\del[1]{\sqrt{2\log(d)}+4},
		\end{equation}
		where~$\overline{\delta}=|\delta_1|\vee |\delta_2|$, and
		\begin{equation}\label{eq:twogauss2}
			\P\del[1]{Z\in[a+\delta_1,b+\delta_2]}
			\geq
			\P\del[1]{Z \in[a,b]}-\frac{\overline{\delta}}{\underline{\sigma}}\del[1]{\sqrt{2\log(d)}+4}.
		\end{equation}
	\end{lemma}
	\begin{proof}
		Consider first~\eqref{eq:twogauss1}. Clearly,
		\begin{equation*}
			\P\del[1]{Z\in[a+\delta_1,b+\delta_2]}
			\leq 
			\P\del[1]{Z\in[a-\overline{\delta},b+\overline{\delta}]}
			=
			\P\del[1]{(Z',-Z')'\leq (b'+\overline{\delta},-a'+\overline{\delta})'}.
		\end{equation*}
		By Theorem 1 of~\cite{chernozhukov2017detailed}, which trivially remains valid for~$y$ there taking values in~$\overline{\R}^d$, the far right-hand side of the previous display is bounded from above by the sum of~$\P\del[0]{Z \in[a,b]}$ and
		\begin{align*}
			\frac{\overline{\delta}}{\underline{\sigma}}\del[1]{\sqrt{2\log(2d)}+2}	
			=
			\frac{\overline{\delta}}{\underline{\sigma}}\del[1]{\sqrt{2(\log(2)+\log(d))}+2}
			\leq
			\frac{\overline{\delta}}{\underline{\sigma}}\del[1]{\sqrt{2\log(d)}+2+\sqrt{2(\log(2)}},	
		\end{align*}
		which yields the desired result. To prove \eqref{eq:twogauss2}, note that
		\begin{align*}
			\P\del[1]{Z\in[a+\delta_1,b+\delta_2]}
			\geq 
			\P\del[1]{Z\in[a+\overline{\delta},b-\overline{\delta}]}
			=
			\P\del[1]{(Z',-Z')\leq (b'-\overline{\delta},-a'-\overline{\delta})'},
		\end{align*}
		which, by Theorem 1 of~\cite{chernozhukov2017detailed}, is bounded from below by~$\P\del[0]{Z \in[a,b]}$ minus the left-hand side of the penultimate display, implying~\eqref{eq:twogauss2}.
	\end{proof}
	
	We will use the following notation in the proof of Theorem~\ref{thm:HDGauss_studentized}.  For~$x\in\R^d$, let~$\enVert[0]{x}_\infty=\max_{j=1,\hdots,d}|x_j|$. For any $A\subseteq \R^d$ and~$\zeta>0$, let $$A^{\zeta,\infty}=\cbr[1]{x\in\R^d:\inf_{y\in A}\|x-y\|_\infty\leq \zeta}.$$ Furthermore, $$A^{-\zeta,\infty}=\cbr[1]{x\in\R^d:\mc{B}_\infty(x,\zeta) \subseteq A} \quad \text{ where } \quad \mc{B}_{\infty}(x,\zeta)=\cbr[1]{y\in\R^d:\|y-x\|_{\infty} \leq \zeta}.$$
	
	\begin{proof}[Proof of Theorem~\ref{thm:HDGauss_studentized}]
		
		We first establish~\eqref{eq:HDGaussstudentized}. Writing~$\mathfrak{T}_n^{(3)}$ for the left-hand side of~\eqref{eq:epscond'}, it follows from~$\frac{\log(d^2n)}{n}\leq {\lambda_2'}^{-1}{\eps_n'}$ that
		\begin{align*}
			\mathfrak{T}_n^{(3)}
			\leq
			2\eps_n'+{\lambda_2'}^{-1}{\eps_n'}+\sqrt{{\lambda_2'}^{-2}{\eps_n'}^2+4{\lambda_2'}^{-1}{\eps_n'}^2}
			\leq
			K\eps_n',	
		\end{align*}
		for a constant~$K$ depending only on~$\lambda_2'$. Assume that
		\begin{align}\label{eq:CaseNormalize}
			K\eps_n'<1 \qquad\text{and}\qquad C_1\sbr[3]{\overline{\eta}_n^{1-\frac{2}{m}}+\del[2]{\frac{\log(dn)}{n}}^{1-\frac{1}{(m/2)\wedge 2}}}\leq \frac{b_1^2}{2},
		\end{align}
		where~$C_1$ is the constant from Theorem~\ref{thm:covestimGRAM} depending only on~$b_2,\lambda_1',\lambda_2'$, and~$m$. Note that~$K\eps_n'<1$ implies, by the penultimate display, that~\eqref{eq:epscond'} is satisfied which is used when invoking Theorem~\ref{thm:covestimGRAM} below. If either of the conditions in~\eqref{eq:CaseNormalize} is violated, we conclude~\eqref{eq:HDGaussstudentized} with similar arguments as those employed in the end of the proofs of Theorems~\ref{thm:HDGauss} and \ref{thm:HDBootstrap}, respectively.

		Recall the definition of~$S_{n,W,S}$ in~\eqref{eq:S_nWS}, let~$T_n\in\R^d$ have entries
		\begin{equation*}
			T_{n,j}=\frac{1}{\sqrt{n}\sigma_{2,j}}\sum_{i=1}^n\sbr[1]{\phi_{\hat\alpha_j,\hat\beta_j}(\tilde{X}_{i,j})-\mu_j},\qquad j=1,\hdots,d,	
		\end{equation*}
		and, observe that (grant the quotients are well-defined)
		\begin{equation}
			A_n
			:=
			\enVert[1]{S_{n,W,S}-T_n}_\infty
			\leq 
			\max_{j=1,\hdots,d}\envert[3]{\frac{1}{\tilde{\sigma}_{n,j}}-\frac{1}{\sigma_{2,j}}}
			\cdot\max_{j=1,\hdots,d}\envert[3]{\frac{1}{\sqrt{n}}\sum_{i=1}^n\sbr[1]{\phi_{\hat\alpha_j,\hat\beta_j}(\tilde{X}_{i,j})-\mu_j}}\label{eq:An1}.
		\end{equation}
		By Theorem~\ref{thm:covestimGRAM},~$\min_{j=1,\hdots,d}\sigma_{2,j}\geq b_1$, and the mean-value theorem, there exists a constant~$C_2$ depending only on $b_1,b_2,\lambda_1',\lambda_2'$, and~$m$ such that on a set of probability at least~$1-24/n$ it holds that~$\tilde{\sigma}_{n,j} \geq b_1/\sqrt{2} > 0$ for every~$j = 1, \hdots, d$ (we used~\eqref{eq:CaseNormalize}) and
		\begin{equation}
			\max_{j=1,\hdots,d}\envert[3]{\frac{1}{\tilde{\sigma}_{n,j}}-\frac{1}{\sigma_{2,j}}}
			=
			\max_{j=1,\hdots,d}\frac{|\tilde{\sigma}_{n,j}-\sigma_{2,j}|}{\tilde{\sigma}_{n,j}\sigma_{2,j}}
			\leq 
			C_2\del[3]{\overline{\eta}_n^{1-\frac{2}{m}}+\del[2]{\frac{\log(dn)}{n}}^{1-\frac{1}{(m/2)\wedge 2}}}\label{eq:An2}.
		\end{equation}
		Furthermore, as~$\sigma_{2,j}\leq \sigma_{m,j}\leq b_2$, the union bound and~$\P\del[0]{|\bm{z}|> t}\leq 2\exp(-t^2/2)$ for~$t\geq 0$ and~$\bm{z} \sim \mathsf{N}_1(0,1)$ yields for~$Z \sim \mathsf{N}_d(0, \Sigma)$ that
		\begin{equation*}
			\P\del{\max_{j=1,\hdots,d}|Z_j|> b_2\sqrt{2\log(2dn)}}
			\leq
			\frac{1}{n}.
		\end{equation*}
		Thus, by Theorem~\ref{thm:HDGauss}, and writing
		\begin{align*}
			r_n:=&\sbr[3]{\frac{\log^{5-\frac{2}{m}}(dn)}{n^{1-\frac{2}{m}}}}^{\frac{1}{4}}
			+
			\sbr[3]{\overline{\eta}_n^{1-\frac{1}{m}}+\sbr[2]{\frac{\log(dn)}{n}}^{1-\frac{1}{m}}}\sqrt{n\log(d)} \\
			&+
			\del[4]{\log^2(d)\sbr[3]{\overline{\eta}_n^{1-\frac{2}{m}}+\sbr[2]{\frac{\log(dn)}{n}}^{1-\frac{2}{m}}}}^{1/2},
		\end{align*}
		it follows that there exists a constant~$K_1$ depending only on~$b_1,b_2,\lambda_1,\lambda_2,$ and~$m$ such that
		\begin{align}
			\P\del[3]{\max_{j=1,\hdots,d}\envert[3]{\frac{1}{\sqrt{n}}\sum_{i=1}^n\sbr[1]{\phi_{\hat\alpha_j,\hat\beta_j}(\tilde{X}_{i,j})-\mu_j}}> b_2\sqrt{2\log(2dn)}}
			\leq 
			K_1r_n+\frac{1}{n}
			\leq
			K_1r_n\label{eq:An3}, 
		\end{align} 
		where the value of~$K_1$ is suitably adjusted to justify the second inequality. Hence, by~\eqref{eq:An1}--\eqref{eq:An3} there exists a constant~$C_3$ depending only on~$b_1, b_2,\lambda_1',\lambda_2'$, and~$m$ such that
		\begin{equation}
			\P\del[3]{A_n\leq C_3\sqrt{\log(dn)}\sbr[2]{\overline{\eta}_n^{1-\frac{2}{m}}+\del[2]{\frac{\log(dn)}{n}}^{1-\frac{1}{(m/2)\wedge 2}}}}
			\geq 
			1-\frac{24}{n}-K_1r_n
			\geq
			1-Kr_n\label{eq:An4},
		\end{equation}
		where~$K$ depends only on~$b_1,b_2,\lambda_1,\lambda_2$, and~$m$. Thus, writing
		\begin{equation*}
			\overline{A}_n:=C_3\sqrt{\log(dn)}\sbr[3]{\overline{\eta}_n^{1-\frac{2}{m}}+\del[2]{\frac{\log(dn)}{n}}^{1-\frac{1}{(m/2)\wedge 2}}},
		\end{equation*}
		it holds for all~$H\in\mc{H}$ that
		\begin{equation}\label{eq:setdecomp1}
			\cbr[1]{S_{n,W,S}\in H}
			\subseteq 
			\cbr[1]{T_n\in H^{\overline{A}_n,\infty}}
			\bigcup \cbr[1]{A_n>\overline{A}_n}.
		\end{equation} 
		Writing~$Y_i=D^{-1}X_i$,~$\tilde{Y}_i=D^{-1}\tilde{X}_i$ as well as~$\check{\alpha}_j=\tilde Y_{\lceil \eps_n n \rceil,j}^*$ and $\check{\beta}_j=\tilde Y_{\lfloor(1-\eps_n )n\rfloor+1,j}^*$, note that
		\begin{equation*}
			T_{n,j}=\frac{1}{\sqrt{n}\sigma_{2,j}}\sum_{i=1}^n\sbr[1]{\phi_{\hat\alpha_j,\hat\beta_j}(\tilde{X}_{i,j})-\mu_j}
			=
			\frac{1}{\sqrt{n}}\sum_{i=1}^n\sbr[1]{\phi_{\check\alpha_j,\check\beta_j}(\tilde{Y}_{i,j})-\E Y_{1,j}},\quad \text{for }j=1,\hdots,d,
		\end{equation*}
		such that by Theorem~\ref{thm:HDGauss} and the covariance matrix of~$Y_1$ being~$\Sigma_0=D^{-1}\Sigma D^{-1}$, there exists a constant~$C$ depending only on~$b_1,b_2,\lambda_1,\lambda_2$, and~$m$ such that (using that~$H^{\overline{A}_n,\infty}\in\mc{H}$)
		\begin{equation*}
			\envert[2]{\P\del[1]{T_n\in H^{\overline{A}_n,\infty}}-\P\del[1]{Z'\in H^{\overline{A}_n,\infty}}}
			\leq
			Cr_n.
		\end{equation*}
		Furthermore, by Lemma~\ref{lem:anticonc} (applied with~$\delta_1=-\overline{A}_n$ and~$\delta_2=\overline{A}_n$) and~$Z'\sim\mathsf{N}_d(0,\Sigma_0)$ where~$\Sigma_{0,j,j}=1$ for all~$j=1,\hdots,d$, 
		\begin{equation*}
			0
			\leq 
			\P\del[1]{Z'\in H^{\overline{A}_n,\infty}}
			-
			\P\del[1]{Z'\in H}
			\leq 
			\overline{A}_n\del[1]{\sqrt{2\log(d)}+4}
			\leq
			7\overline{A}_n\sqrt{\log(d)},
		\end{equation*}
		the last inequality following from~$d\geq 2$. Using this in~\eqref{eq:setdecomp1} along with~\eqref{eq:An4} yields that
		\begin{equation*}
			\P\del[1]{S_{n,W,S}\in H}-\P\del[1]{Z'\in H}
			\leq
			(C+K)r_n+7\overline{A}_n\sqrt{\log(d)}. 	
		\end{equation*}
		Finally, since 
		\begin{equation*}
			\cbr[1]{T_n\in H^{-\overline{A}_n,\infty}}
			\subseteq 
			\cbr[1]{S_{n,W,S}\in H}	
			\bigcup \cbr[1]{A_n>\overline{A}_n},
		\end{equation*}
		the same arguments as those following~\eqref{eq:setdecomp1} lead to
		\begin{equation*}
			-(C+K)r_n-7\overline{A}_n\sqrt{\log(d)}
			\leq
			\P\del[1]{S_{n,W,S}\in H}-\P\del[1]{Z'\in H},
		\end{equation*}
		which yields~\eqref{eq:HDGaussstudentized} under the assumption of~\eqref{eq:CaseNormalize}, upon adjusting constants as the upper and lower bounds just obtained do not depend on~$H$. 
		
		Finally, since~$r_n\to 0$ by~\eqref{eq:conv1} and~\eqref{eq:conv2} as well as~$\overline{A}_n\sqrt{\log(d)}\to 0$ by~\eqref{eq:conv4} of Lemma~\ref{lem:covcompare} it follows that~$\rho_{n,W,S}\to 0$. 
	\end{proof}

	\begin{proof}[Proof of Theorem~\ref{thm:HDBootstrap_studentized}]
		The proof is similar to that of Theorem~\ref{thm:HDBootstrap}, but we include it here for completeness. By the triangle inequality~$\tilde{\rho}_{n,W,S}$ is bounded from above by the sum of
		\begin{align*}
			\rho_{n,W,S}=\sup_{H\in\mc{H}}\envert[2]{\P\del[1]{S_{n,W,S}\in H}-\P\del[1]{Z'\in H}},	
		\end{align*}
		where we recall that~$Z'\sim \mathsf{N}_d(0,\Sigma_0)$, and
		\begin{align*}
			B_n:=\sup_{H\in\mc{H}}\envert[2]{\P\del[1]{Z'\in H}-\P\del[1]{\hat{Z}'\in H\mid\tilde{X}_1,\hdots,\tilde{X}_n}}.
		\end{align*}
		First, $\rho_{n,W,S}\leq \mathfrak{B}_n$, where~$\mathfrak{B}_n$ is the upper bound on~$\rho_{n,W,S}$ in Theorem~\ref{thm:HDGauss_studentized}. Next, by the Gaussian-to-Gaussian comparison inequality as stated in Proposition 2.1 of~\cite{chernozhuokov2022improved} (cf.~also Proposition 2 in~\cite{chernozhukov2023high}),
		\begin{align*}
			B_n
			\leq 
			C\del[2]{\log^2(d)\max_{1\leq j,k\leq d}\envert[1]{\tilde{\Sigma}_{n,0,j,k}-\Sigma_{0,j,k}}}^{1/2},
		\end{align*}
		for an absolute constant~$C$ (since~$\Sigma_{0,j,j}=1$ for all~$j=1,\hdots,d$).
		
		Arguing as in the opening paragraph of the proof of Theorem~\ref{thm:HDGauss_studentized}, it suffices to consider the case 
		of~\eqref{eq:CaseNormalize} being satisfied.
		
		Let~$C_1$ be the constant from Theorem~\ref{thm:covestimGRAM} depending only on~$b_2,\lambda_1',\lambda_2'$, and~$m$. By that theorem (note that~\eqref{eq:epscond'} is ensured to be true by~\eqref{eq:CaseNormalize} so that Theorem~\ref{thm:covestimGRAM} can be employed) there exists a set~$E_n$ of probability at least~$1-\frac{24}{n}$ on which
		\begin{equation*}
			\max_{1\leq j,k\leq d}\envert[1]{\tilde{\Sigma}_{n,j,k}-\Sigma_{j,k}}
			\leq	
			C_1\sbr[3]{\overline{\eta}_n^{1-\frac{2}{m}}+\del[2]{\frac{\log(dn)}{n}}^{1-\frac{1}{(m/2)\wedge 2}}}\leq \frac{b_1^2}{2}.
		\end{equation*}
		Therefore, on~$E_n$,~$\min_{j=1,\hdots,d}\tilde{\sigma}_{n,j}^2\geq b_1^2/2$ because by assumption~$\min_{j=1,\hdots,d}\sigma_{2,j}^2\geq b_1^2$. Furthermore, note that for all~$1\leq j,k\leq d$
		\begin{equation*}
			\tilde{\Sigma}_{n,0,j,k}-\Sigma_{0,j,k}
			=
			\frac{\tilde{\Sigma}_{n,j,k}}{\tilde{\sigma}_{n,j}\tilde{\sigma}_{n,k}}-\frac{\Sigma_{j,k}}{\sigma_{2,j}\sigma_{2,k}}
			=
			\frac{(\tilde{\Sigma}_{n,j,k}-\Sigma_{j,k})\sigma_{2,j}\sigma_{2,k}+\Sigma_{j,k}(\sigma_{2,j}\sigma_{2,k}-\tilde{\sigma}_{n,j}\tilde{\sigma}_{n,k})}{\tilde{\sigma}_{n,j}\tilde{\sigma}_{n,k}\sigma_{2,j}\sigma_{2,k}},	
		\end{equation*}
		which is well-defined on~$E_n$. Thus, by the mean-value theorem, there exists a constant~$C_2$ depending only on~$b_1,b_2,\lambda_1',\lambda_2'$, and~$m$ such that on~$E_n$
		\begin{align*}
			\max_{1\leq j,k\leq d}\envert[1]{\tilde{\Sigma}_{n,0,j,k}-\Sigma_{0, j,k}}
			\leq
			C_2\del[3]{\overline{\eta}_n^{1-\frac{2}{m}}+\del[2]{\frac{\log(dn)}{n}}^{1-\frac{1}{(m/2)\wedge 2}}}.
		\end{align*}
		Hence, with probability at least~$1-\frac{24}{n}$
		\begin{align*}
			B_n
			\leq
			C\del[4]{\log^2(d)\sbr[3]{\overline{\eta}_n^{1-\frac{2}{m}}+\del[2]{\frac{\log(dn)}{n}}^{1-\frac{1}{(m/2)\wedge 2}}}}^{1/2}	
		\end{align*}
		for a constant~$C$ depending only on~$b_1,b_2,\lambda_1',\lambda_2'$, and~$m$ which implies the bound on~$\tilde{\rho}_{n,W,S}$ in~\eqref{eq:StudentizeBS} in case of~\eqref{eq:CaseNormalize}.
		
		Finally,~$\tilde{\rho}_{n,W,S}\to 0$ by~\eqref{eq:conv1}--\eqref{eq:conv4} of Lemma~\ref{lem:covcompare}.
	\end{proof}

	\section{Proofs for Section~\ref{sec:trimmedmean}}
	
	\begin{proof}[Proof of Theorem~\ref{thm:HDGaussTrim}]
		Since~$\log(dn)/n\leq \lambda_2^{-1}\eps_n$, one has that~$\mathfrak{T}_n^{(1)}$ in~\eqref{eq:epscond} is bounded as
		\begin{align*}
			\mathfrak{T}_n^{(1)}
			\leq
			2\eps_n+\lambda_2^{-1}\eps_n+\sqrt{\lambda_2^{-2}\eps_n^2+4\lambda_2^{-1}\eps_n^2}\leq
			C_1\eps_n,
		\end{align*}
		for a constant~$C_1$ depending only on~$\lambda_2$. In the sequel we impose that~$C_1\eps_n<1/2$ such that by the previous display~\eqref{eq:epscond} is satisfied allowing us to invoke Lemma \ref{lem:quantiles} below. In case~$C_1\eps_n\geq 1/2$,~\eqref{eq:trimGauss} trivially holds by the same arguments as in the end of the proof of Theorem~\ref{thm:HDGauss}. 
		
		Fix~$j\in\cbr[0]{1,\hdots,d}$. Since~$\phi_{\hat\alpha_j,\hat\beta_j}(\tilde{X}_{i,j}^*)=\tilde{X}_{i,j}^*$ for $i\in I_n=\cbr[1]{\lceil\eps_nn\rceil,\hdots,\lfloor(1-\eps_n)n\rfloor+1}$
		\begin{align*}
			A_{n,j}:&=
			\envert[1]{S_{n,T,j}-S_{n,W,j}}\\
			&=
			\sqrt{n}\envert[4]{\frac{1}{|I_n|}\sum_{i=\lceil \eps_nn\rceil}^{\lfloor(1-\eps_n)n\rfloor+1}\sbr[1]{\tilde{X}_{i,j}^*-\mu_j}-\frac{1}{n}\sum_{i=1}^n\sbr[1]{\phi_{\hat{\alpha}_j,\hat{\beta}_j}(\tilde{X}_{i,j}^*)-\mu_j}}\\
			&=
			\sqrt{n}\envert[4]{\frac{2(\lceil\eps_nn\rceil-1)}{n|I_n|}\sum_{i=\lceil \eps_nn\rceil}^{\lfloor(1-\eps_n)n\rfloor+1}\sbr[1]{\tilde{X}_{i,j}^*-\mu_j}-\frac{\lceil \eps_nn\rceil-1}{n}\sbr[1]{\hat{\alpha}_j-\mu_j}-\frac{\lceil \eps_nn\rceil-1}{n}\sbr[1]{\hat{\beta}_j-\mu_j}}\\
			&\leq 
			\sqrt{n}\del[4]{\frac{2(\lceil\eps_nn\rceil-1)}{n|I_n|}\sum_{i=\lceil \eps_nn\rceil}^{\lfloor(1-\eps_n)n\rfloor+1}\envert[1]{\tilde{X}_{i,j}^*-\mu_j}
				+
				\frac{\lceil \eps_nn\rceil-1}{n}\envert[1]{\hat{\alpha}_j-\mu_j}+\frac{\lceil \eps_nn\rceil-1}{n}\envert[1]{\hat{\beta}_j-\mu_j}}.	
		\end{align*}
		Next, since
		\begin{align*}
			\frac{2(\lceil\eps_nn\rceil-1)}{n|I_n|}\sum_{i=\lceil \eps_nn\rceil}^{\lfloor(1-\eps_n)n\rfloor+1}\envert[1]{\tilde{X}_{i,j}^*-\mu_j}
			&\leq
			\frac{2(\lceil\eps_nn\rceil-1)}{n}\max_{i=\lceil \eps_nn\rceil,\hdots,\lfloor(1-\eps_n)n\rfloor+1}\envert[1]{\tilde{X}_{i,j}^*-\mu_j}\\
			&=
			\frac{2(\lceil\eps_nn\rceil-1)}{n}\sbr[1]{|\hat{\alpha}_j-\mu_j|\vee |\hat{\beta}_j-\mu_j|},	
		\end{align*}
		it suffices to bound~$|\hat{\alpha}_j-\mu_j|\vee|\hat{\beta}_j-\mu_j|$ from above. To this end, by Lemma~\ref{lem:quantiles} and~$\hat{\alpha}_j-\mu_j\leq \hat{\beta}_j-\mu_j$ (\eqref{eq:epscond} implies~$\eps_n<1/4$ such that~$\hat{\alpha}_j=\tilde{X}_{\lceil{\eps_nn}\rceil,j}^*\leq \tilde{X}_{\lfloor(1-\eps_n)n\rfloor+1,j}^*= \hat{\beta}_j$) it holds with probability at least~$1-\frac{2}{dn}$ that 
		\begin{equation*}
			|\hat{\alpha}_j-\mu_j|\vee|\hat{\beta}_j-\mu_j|
			\leq
			|Q_{c_{1,n}\eps_n}(X_{1,j})-\mu_j|\vee |Q_{1-c_{1,n}\eps_n}(X_{1,j})-\mu_j|
			\leq
			\frac{\sigma_m}{(c_{1,n}\eps_n)^{1/m}},
		\end{equation*}
		where the last estimate is by Lemma~\ref{lem:quantile_mean}. Therefore, by Lemma~\ref{lem:cControl} there exists a constant~$C$ depending only on~$b_2,\lambda_1,\lambda_2$, and~$m$ such that with probability at least~$1-\frac{2}{dn}$,
		\begin{equation*}
			A_{n,j}
			\leq
			\frac{4\sqrt{n}(\lceil{\eps_nn}\rceil-1)\sigma_m}{n(c_{1,n}\eps_n)^{1/m}}
			\leq
			C\sqrt{n}\eps_n^{1-\frac{1}{m}}
			\leq
			C\sqrt{n}\del[3]{\overline{\eta}_n^{1-\frac{1}{m}}+\sbr[2]{\frac{\log(dn)}{n}}^{1-\frac{1}{m}}}=:\overline{A}_n,
		\end{equation*}
		with~$C$ potentially changing values in the last inequality and subadditivity of~$z\mapsto z^{1-\frac{1}{m}}$ was used along with the definition of~$\eps_n$ in~\eqref{eq:epsfam}. Therefore, since the right-hand side does not depend on~$j$, it follows by the union bound over~$j=1,\hdots,d$ that with probability at least~$1-\frac{2}{n}$
		\begin{equation*}
			A_n:=\max_{j=1,\hdots,d}A_{n,j}\leq \overline{A}_n.
		\end{equation*}
		Next, observe that for all~$H\in\mc{H}$ (recalling also the notation introduced prior to the proof of Theorem~\ref{thm:HDGauss_studentized})
		\begin{equation*}
			\cbr[1]{S_{n,T}\in H}
			\subseteq 
			\cbr[1]{S_{n,W}\in H^{\overline{A}_n,\infty}}
			\bigcup \cbr[1]{A_n>\overline{A}_n}.
		\end{equation*} 
		By Theorem~\ref{thm:HDGauss},~$H^{\overline{A}_n,\infty}\in\mc{H}$, and Lemma~\ref{lem:anticonc},
		\begin{equation*}
			\P\del[1]{S_{n,W}\in H^{\overline{A}_n,\infty}}
			-\P\del[1]{Z\in H} 
			\leq 
			\mathfrak{A}_n+\frac{\overline{A}_n}{b_1}\del[1]{\sqrt{2\log(d)}+4}
			\leq
			\mathfrak{A}_n+C\overline{A}_n\sqrt{\log(d)},
		\end{equation*}
		for a constant~$C$ depending only on~$b_1$. Thus, since ~$\P\del[1]{A_n> \overline{A}_n}\leq \frac{2}{n}$, we conclude that
		\begin{equation}\label{eq:Truncupper}
			\P\del[1]{S_{n,T}\in H}
			-
			\P\del[1]{Z\in H}
			\leq
			\mathfrak{A}_n+C\overline{A}_n\sqrt{\log(d)}+\frac{2}{n}.
		\end{equation}
		Since also
		\begin{equation*}
			\cbr[1]{S_{n,W}\in H^{-\overline{A}_n,\infty}}
			\subseteq 
			\cbr[1]{S_{n,T}\in H}
			\bigcup \cbr[1]{A_n>\overline{A}_n},
		\end{equation*} 
		an identical argument shows that
		\begin{equation}\label{eq:Trunclower}
			-\mathfrak{A}_n-C\overline{A}_n\sqrt{\log(d)}-\frac{2}{n}
			\leq
			\P\del[1]{S_{n,T}\in H}
			-
			\P\del[1]{Z\in H},
		\end{equation}
		which, together with the penultimate display,~\eqref{eq:Truncupper}--\eqref{eq:Trunclower} not depending on~$H\in\mc{H}$, and dominating~$2/n$ by~$\overline{A}_n\sqrt{\log(d)}$ implies the bound in~\eqref{eq:trimGauss}. That~$\rho_{n,T}\to 0$ follows from~\eqref{eq:conv1} and~\eqref{eq:conv2} of Lemma~\ref{lem:covcompare}.
	\end{proof}

	\begin{proof}[Proof of Theorem~\ref{thm:HDBootstrapTrim}]
		The proof is almost identical to that of~Theorem~\ref{thm:HDBootstrap}, but is included for completeness. Note first that writing~$\mathfrak{T}_n^{(3)}$ for the left-hand side of~\eqref{eq:epscond'}, it follows from~$\frac{\log(d^2n)}{n}\leq {\lambda_2'}^{-1}{\eps_n'}$ that
		\begin{align*}
			\mathfrak{T}_n^{(3)}
			\leq
			2\eps_n'+{\lambda_2'}^{-1}{\eps_n'}+\sqrt{{\lambda_2'}^{-2}{\eps_n'}^2+4{\lambda_2'}^{-1}{\eps_n'}^2}
			\leq
			C_1\eps_n',	
		\end{align*}
		for a constant~$C_1$ depending only on~$\lambda_2'$. 
		
		Suppose first that~$C_1\eps_n'<1$ such that~\eqref{eq:epscond'} is satisfied by the previous display. By the triangle inequality~$\tilde{\rho}_{n,T}$ is bounded from above by the sum of
		\begin{align*}
			\rho_{n,T}=\sup_{H\in\mc{H}}\envert[2]{\P\del[1]{S_{n,T}\in H}-\P\del[1]{Z\in H}},	
		\end{align*}
		where we recall that~$Z\sim\mathsf{N}_d(0,\Sigma)$, and
		\begin{align*}
			B_n:=\sup_{H\in\mc{H}}\envert[2]{\P\del[1]{Z\in H}-\P\del[1]{\tilde{Z}\in H\mid\tilde{X}_1,\hdots,\tilde{X}_n}}.
		\end{align*}
		First, $\rho_{n,T}\leq \mathfrak{C}_n$, where~$\mathfrak{C}_n$ is the upper bound on~$\rho_{n,T}$ in Theorem~\ref{thm:HDGaussTrim}. Next, by the Gaussian-to-Gaussian comparison inequality as stated in Proposition 2.1 of~\cite{chernozhuokov2022improved} (cf.~also Proposition 2 in~\cite{chernozhukov2023high}),
		\begin{align*}
			B_n
			\leq 
			C\del[2]{\log^2(d)\max_{1\leq j,k\leq d}\envert[1]{\tilde{\Sigma}_{n,j,k}-\Sigma_{j,k}}}^{1/2}
			\leq
			C\del[4]{\log^2(d)\sbr[3]{\overline{\eta}_n^{1-\frac{2}{m}}+\del[2]{\frac{\log(dn)}{n}}^{1-\frac{1}{(m/2)\wedge 2}}}}^{1/2}	,
		\end{align*}
		the last inequality holding with probability at least~$1-\frac{24}{n}$ by Theorem~\ref{thm:covestimGRAM} and~$C$ being a constant depending on~$b_1,b_2,\lambda_1',\lambda_2'$ and~$m$ only.
		
		If, on the other hand,~$C_1\eps_n'\geq 1$ then by the definition of~$\eps_n'$ in~\eqref{eq:epsprime}%
		\begin{align*}
			\overline{\eta}_n\geq \frac{1}{2C_1\lambda_1'}\qquad\text{or}\qquad \frac{\log(dn)}{n}
			=
			\frac{\log([dn]^2)/2}{n}
			\geq
			\frac{\log(d^2n)}{2n}
			\geq
			\frac{1}{4C_1\lambda_2'}.
		\end{align*} 
		In either case, the upper bound on~$\tilde{\rho}_{n,T}$ remains valid by suitably adjusting~$C$ (if needed).
		
		That~$\tilde{\rho}_{n,T}\to 0$ in probability follows from~\eqref{eq:conv1}--\eqref{eq:conv3}	 of Lemma~\ref{lem:covcompare}.
	\end{proof}

	\section{Auxiliary lemmas}
	This section gathers some auxiliary lemmas, the proofs of which largely follow from related results in~\cite{Wins1}.
	
	The following standard result is Lemma C.1 from~\cite{Wins1}. It bounds the difference between the mean and quantile of a distribution of a random variable~$Z$ (which is not necessarily continuous).
	\begin{lemma}\label{lem:quantile_mean}
		Let~$Z$ satisfy~$\sigma_m^m:=\E|Z-\E Z|^m\in[0,\infty)$ for some~$m\in[1,\infty)$. Then, for all~$p\in(0,1)$,
		\begin{equation}\label{eq:QM}
			\E Z-\frac{\sigma_m}{p^{1/m}}\leq Q_p(Z)\leq \E Z+\frac{\sigma_m}{(1-p)^{1/m}}.
		\end{equation}
	\end{lemma}
	
	\begin{lemma}\label{lem:cControl}
		For~$c_{1,n}$ as defined in~\eqref{eq:c1},~$c_{2,n}$ as defined in~\eqref{eq:c2}, and $\eps_n$ as defined in~\eqref{eq:epsfam} and satisfying~\eqref{eq:epscond}, it holds that
		\begin{align*}
			0<(1-\lambda_1^{-1})\exp \del[2]{{-\frac{1}{\lambda_2(1-\lambda_1^{-1})}-1}}\leq\ &c_{1,n}<1,\\
			1<\ &c_{2,n} \leq 2+\lambda_2^{-1}+\sqrt{\lambda_2^{-2}+4\lambda_2^{-1}},
		\end{align*}
		and
		\begin{align}\label{eq:welldefined}
			0
			<
			\eps_n\min(c_{1,n},c_{2,n})
			&\leq 
			\eps_n(c_{1,n}+c_{2,n})\notag\\
			&\leq
			2\eps_n +\frac{\log(dn)}{n}+\sqrt{\del[2]{\frac{\log(dn)}{n}}^2+4\frac{\log(dn)}{n}\eps_n}.
		\end{align}	
	\end{lemma}
	\begin{proof}
		Apply Lemma~B.3 in~\cite{Wins1} with~$\delta=\frac{6}{dn}$ (recall that $d n > 6$ is assumed throughout; cf.~the sentence right before~Theorem~\ref{thm:HDGauss}),~$\epsilon = \varepsilon_n$, and~$\eta=\overline{\eta}_n$. To this end, note that~$c_1$, and~$c_2$ \emph{there} equal~$c_{1,n}$ and~$c_{2,n}$, respectively, for~$\delta=\frac{6}{dn}$, that~$\eps_n \geq \lambda_2 \log(dn)/n$ by definition, and that~$\varepsilon_n < 1/4$ follows from~\eqref{eq:epscond}.
	\end{proof}

	The following lemma shows that for~$\eps_n$ as defined in~\eqref{eq:epsfam}, the lower and upper~$\eps_n$ order statistics of the contaminated data are close to related population quantiles of the uncontaminated data. 
	\begin{lemma}\label{lem:quantiles}
		Fix~$j\in\cbr[0]{1,\hdots,d}$,~$n\in\N$, let Assumption~\ref{ass:setting} be satisfied, and~$\eps_n$ be as defined in~\eqref{eq:epsfam} satisfying~\eqref{eq:epscond}. Then, for~$c_{1,n}$ as in~\eqref{eq:c1} and~$c_{2,n}$ as in~\eqref{eq:c2}, we have~$$0<\eps_n\min(c_{1,n},c_{2,n})\leq \eps_n(c_{1,n}+c_{2,n})<1/2,$$ and each of \eqref{eq:LBLQ}--\eqref{eq:UBUQ} below holds with probability at least~$1-\frac{1}{dn}$:
		\begin{eqnarray}
			\tilde X_{\lceil \eps_n n \rceil, j}^* 
			&\geq&
			Q_{c_{1,n}\eps_n}(X_{1,j})\label{eq:LBLQ};\\
			\tilde X_{\lceil(1-\eps_n)n\rceil , j}^*
			&\geq&
			Q_{1-c_{2,n}\eps_n}(X_{1,j})\label{eq:LBUQ};\\
			\tilde{X}_{\lfloor \eps_n n \rfloor+1, j}^*
			&\leq&
			Q_{c_{2,n}\eps_n}(X_{1,j})\label{eq:UBLQ};\\
			\tilde{X}_{\lfloor (1-\eps_n)n\rfloor+1, j}^*
			&\leq&
			Q_{1-c_{1,n}\eps_n}(X_{1,j}).\label{eq:UBUQ}
		\end{eqnarray}	
	\end{lemma}	
	
	\begin{proof}
		The first statement in the present lemma follows from~\eqref{eq:epscond} along with~\eqref{eq:welldefined} of~Lemma~\ref{lem:cControl}. For the remaining statements, apply Lemma~B.5 in~\cite{Wins1} with~$\delta=\frac{6}{dn}$ (recall that $d n > 6$ is implied by~$d\geq 2$ and~$n>3$, which is assumed throughout; cf.~the sentence right before~Theorem~\ref{thm:HDGauss}),~$\epsilon = \varepsilon_n$,  and~$\eta=\overline{\eta}_n$. To this end, note that
		\begin{enumerate}
			\item $c_1$ and~$c_2$ \emph{there} equal~$c_{1,n}$ and~$c_{2,n}$, respectively, for~$\delta=\frac{6}{dn}$,
			\item that~$\eps_n \geq \lambda_1 \overline{\eta}_n$ by definition and~$\varepsilon_n c_{2,n} < 1/2$ (which has already been verified),
			\item $\varepsilon_n < 1/4$ follows from~\eqref{eq:epscond},
		\end{enumerate}
		and that our Assumption~\ref{ass:setting} implies the remaining assumptions \emph{there}. 
	\end{proof}

	\begin{lemma}\label{lem:noisecontrol}
		There exists a positive constant~$C$ depending only on~$\lambda_1,\lambda_2$, and~$m$, such that if~$\eps_n$ is as in~\eqref{eq:epsfam} satisfying~\eqref{eq:epscond} and Assumption~\ref{ass:setting} is satisfied, then for all~$j=1,\hdots,d$,
		\begin{equation}\label{eq:distcorrupt1}
			\envert[3]{\frac{1}{n}\sum_{i=1}^n\sbr[1]{\phi_{\underline{\alpha}_j,\underline{\beta}_j}(\tilde{X}_{i,j})-\phi_{\underline{\alpha}_j,\underline{\beta}_j}(X_{i,j})}}
			\leq
			C \sigma_m \overline{\eta}_n^{1-\frac{1}{m}}
		\end{equation}
		and
		\begin{equation}\label{eq:distcorrupt2}
			\envert[3]{\frac{1}{n}\sum_{i=1}^n\sbr[1]{\phi_{\overline{\alpha}_j,\overline{\beta}_j}(\tilde{X}_{i,j})-\phi_{\overline{\alpha}_j,\overline{\beta}_j}(X_{i,j})}}
			\leq
			C \sigma_m \overline{\eta}_n^{1-\frac{1}{m}}.
		\end{equation}
	\end{lemma}
	\begin{proof}
		Fix~$j\in\cbr[0]{1,\hdots,d}$ and recall the definitions of~$\underline{\alpha}_j,\overline{\alpha}_j,\underline{\beta}_j$ and~$\overline{\beta}_j$ from~\eqref{eq:alphas} and~\eqref{eq:betas}. The case of~$\overline{\eta}_n=0$ is trivial. To establish \eqref{eq:distcorrupt1} for~$\overline{\eta}_n>0$, apply Lemma C.2 in~\cite{Wins1} to~$X_{1,j}$ with~$\eta=\overline{\eta}_n$,~$s_1=c_{1,n}\eps_n$, and~$s_2=1-c_{2,n}\eps_n$, recalling that~$\sigma_{m,j}\leq\sigma_m$ and using that~$\eps_n> \lambda_{1} \overline{\eta}_n>\overline{\eta}_n$ as well as the lower bounds on~$c_{1,n}$ and~$c_{2,n}$ from Lemma~\ref{lem:cControl}, noting that our Assumption~\ref{ass:setting} implies the assumptions \emph{there}. The requirement~$s_1<s_2$ in Lemma C.2 in~\cite{Wins1} follows by~\eqref{eq:epscond} and Lemma~\ref{lem:cControl}. \eqref{eq:distcorrupt2} follows in a similar fashion using~$s_1=c_{2,n}\eps_n$ and~$s_2=1-c_{1,n}\eps_n$ instead.
	\end{proof}

	\begin{lemma}\label{lem:meancontrol}
		There exists a positive constant~$C$ depending only on~$\lambda_1,\lambda_2$, and~$m$, such that if~$\eps_n$ is as in~\eqref{eq:epsfam} satisfying~\eqref{eq:epscond} and Assumption~\ref{ass:setting} is satisfied, then for all~$j=1,\hdots,d$,
		\begin{align}\label{eq:mean1}
			\envert[1]{\E\phi_{\underline{\alpha}_j,\underline{\beta}_j}(X_{1,j})-\mu_j}
			\vee 
			\envert[1]{\E\phi_{\overline{\alpha}_j,\overline{\beta}_j}(X_{1,j})-\mu_j}
			\leq
			C\sigma_m\eps_n^{1-\frac{1}{m}}
		\end{align}	
	\end{lemma}
	\begin{proof}
		Fix~$j\in\cbr[0]{1,\hdots,d}$ and recall the definitions of~$\underline{\alpha}_j,\overline{\alpha}_j,\underline{\beta}_j$ and~$\overline{\beta}_j$ from~\eqref{eq:alphas} and~\eqref{eq:betas}. Since~$\E\phi_{\underline{\alpha}_j,\underline{\beta}_j}(X_{1,j})-\mu_j\leq \E\phi_{\overline{\alpha}_j,\overline{\beta}_j}(X_{1,j})-\mu_j$, \eqref{eq:mean1} follows from i) (C.5) of Lemma C.5 in~\cite{Wins1} applied to $X_{1,j}$ with~$s_1=c_{1,n}\eps_n$,~$s_2=1-c_{2,n}\eps_n$ and ii) (C.6) of Lemma C.5 in~\cite{Wins1} applied to $X_{1,j}$ with~$s_1=c_{2,n}\eps_n$,~$s_2=1-c_{1,n}\eps_n$, recalling that~$\sigma_{m,j}\leq\sigma_m$ (noting that our Assumption~\ref{ass:setting} implies the assumptions \emph{there}, and in either case the requirement~$s_1<s_2$ imposed in Lemma C.5 in~\cite{Wins1} follows from~\eqref{eq:epscond} and Lemma~\ref{lem:cControl}), such that
		\begin{align*}
			&\envert[1]{\E\phi_{\underline{\alpha}_j,\underline{\beta}_j}(X_{1,j})-\mu_j}\vee \envert[1]{ \E\phi_{\overline{\alpha}_j,\overline{\beta}_j}(X_{1,j})-\mu_j}\\
			&\leq  
			2\sigma_m(c_{1,n}\eps_n)^{1-\frac{1}{m}}+\sigma_m\del[2]{1+\sbr[2]{\frac{c_{2,n}\eps_n}{1-c_{2,n}\eps_n}}^{\frac{1}{m}}}(c_{2,n}\eps_n)^{1-\frac{1}{m}}\\
			&	\leq 
			2\sigma_m(c_{1,n}\eps_n)^{1-\frac{1}{m}}+2\sigma_m(c_{2,n}\eps_n)^{1-\frac{1}{m}},
		\end{align*}
		the second inequality following from~$c_{2,n}\eps_n<1/2$ by Lemma~\ref{lem:cControl} and~\eqref{eq:epscond} and~$(0,1)\ni z\mapsto z/(1-z)$ being increasing. The upper bounds on~$c_{1,n}$ and~$c_{2,n}$ from Lemma~\ref{lem:cControl} now yield the conclusion of the lemma.
	\end{proof}
	
	The following result is similar to Lemma~\ref{lem:cControl} and bounds~$c_{1,n}'$ and~$c_{2,n}'$ as defined in~\eqref{eq:c1'} and~\eqref{eq:c2'}, respectively. It will be used in Lemma \ref{lem:quantiles2} below, as well as in the proof of Theorem~\ref{thm:covestimGRAM}. 
	
	\begin{lemma}\label{lem:cControl'}
		For~$c_{1,n}'$ as defined in~\eqref{eq:c1'},~$c_{2,n}'$ as defined in~\eqref{eq:c2'}, and $\eps_n'$ as defined in~\eqref{eq:epsprime} and satisfying~\eqref{eq:epscond'}, it holds that
		\begin{align*}
			0<(1-\lambda_1'^{-1})\exp \del[2]{{-\frac{1}{\lambda_2'(1-\lambda_1'^{-1})}-1}}\leq\ &c_{1,n}'<1,\\
			1<\ &c_{2,n}' \leq 2+\lambda_2'^{-1}+\sqrt{\lambda_2'^{-2}+4\lambda_2'^{-1}},
		\end{align*}
		and
		\begin{align}\label{eq:welldefined'}
			0
			<
			\eps_n'\min(c_{1,n}',c_{2,n}')
			&\leq 
			\eps_n'(c_{1,n}'+c_{2,n}')\notag\\
			&\leq
			2\eps_n' +\frac{\log(d^2n)}{n}+\sqrt{\del[2]{\frac{\log(d^2n)}{n}}^2+4\frac{\log(d^2n)}{n}\eps_n'}.
		\end{align}	
	\end{lemma}
	\begin{proof}
		The proof is almost identical to that of Lemma~\ref{lem:cControl}, but we now apply Lemma~B.3 in~\cite{Wins1} with~$\delta=\frac{6}{d^2n}$ and~$\epsilon = \varepsilon_n'$. \end{proof}
	
	The following lemma, which follows from Lemma B.5 in~\cite{Wins1}, is an analogue to Lemma~\ref{lem:quantiles}, replacing~$\eps_n$ in~\eqref{eq:epsfam} with~$\eps_n'$ in~\eqref{eq:epsprime}.
	\begin{lemma}\label{lem:quantiles2}
		Let~$Z_{1},\hdots, Z_n$ be i.i.d.~real random variables and suppose that the random variables~$\tilde{Z}_1,\hdots,\tilde{Z}_n$ satisfy
		\begin{equation*}
			\envert[1]{\cbr[1]{i\in\cbr[0]{1,\hdots,n}:\tilde{Z}_i\neq Z_i}}\leq \overline{\eta}_n n.
		\end{equation*}
		Suppose that~$\eps_n'$ is as defined in~\eqref{eq:epsprime} and satisfies~\eqref{eq:epscond'}. Then, for~$c_{1,n}'$ as in~\eqref{eq:c1'} and~$c_{2,n}'$ as in~\eqref{eq:c2'}, we have $$0<\eps_n'\min(c_{1,n}',c_{2,n}')\leq \eps_n'(c_{1,n}'+c_{2,n}')<1,$$ and each of~\eqref{eq:LBLQ2}--\eqref{eq:UBUQ2} below holds with probability at least~$1-\frac{1}{d^2n}$:
		\begin{eqnarray}
			\tilde Z_{\lceil \eps'_nn \rceil}^* 
			&\geq &
			Q_{c_{1,n}'\eps'_n}(Z_{1});\label{eq:LBLQ2}\\
			\tilde Z_{\lceil(1-\eps'_n)n\rceil}^*
			&\geq&
			Q_{1-c_{2,n}'\eps'_n}(Z_{1});\label{eq:LBUQ2}\\
			\tilde{Z}_{\lfloor \eps'_n n \rfloor+1}^*
			&\leq&
			Q_{c_{2,n}'\eps'_n}(Z_{1});\label{eq:UBLQ2}\\
			\tilde{Z}_{\lfloor (1-\eps'_n)n\rfloor+1}^*
			&\leq &
			Q_{1-c_{1,n}'\eps'_n}(Z_{1}).\label{eq:UBUQ2}
		\end{eqnarray}	
	\end{lemma}	
	
	\begin{proof}
		The first statement in the present lemma follows from~\eqref{eq:epscond'} along with~\eqref{eq:welldefined'} of Lemma~\ref{lem:cControl'}. For the remaining statements, apply Lemma B.5 in~\cite{Wins1} with~$\delta = \frac{6}{d^2n}$ (which is strictly less than one since we assume throughout that~$d \geq 2$ and~$n \geq 4 $, cf.~the sentence before Theorem~\ref{thm:HDGauss}),~$\epsilon = \varepsilon_n'$, and $\eta$ there being~$\overline{\eta}_n$. To this end, note that
		\begin{enumerate}
			\item $c_1$, and~$c_2$ \emph{there} equal~$c'_{1,n}$ and~$c'_{2,n}$, respectively, for~$\delta=\frac{6}{d^2n}$,
			\item that~$\eps'_n \geq \lambda'_1 \overline{\eta}_n$ by definition and~$\varepsilon'_n c'_{2,n} < 1$ (which has already been verified),
			\item $\varepsilon'_n < 1/2$ follows from~\eqref{eq:epscond'}.
		\end{enumerate}  
	\end{proof}
	
\end{appendix}

\bibliographystyle{ecta} 
\bibliography{ref}		

\end{document}